\documentclass[pdftex]{amsart}
\pdfoutput=1 

\usepackage{amsmath}
\usepackage{amssymb}
\usepackage{mathrsfs}
\usepackage{bussproofs}
\usepackage{enumerate}
\usepackage{rotating}
\usepackage[matrix,arrow,graph,curve]{xy}

\usepackage{graphicx}
\setkeys{Gin}{keepaspectratio=true} 
\usepackage{calc}

\allowdisplaybreaks

\def\fCenter{\mbox{ :: }}

\theoremstyle{plain}
    \newtheorem{theorem}{Theorem}[section]
    \newtheorem{proposition}[theorem]{Proposition}

\theoremstyle{definition}
    
    \newtheorem{example}[theorem]{Example}
\theoremstyle{remark}

\newcommand{\rb}{\raisebox}
\newcommand{\ig}{\includegraphics}
\newcommand{\scale}{0.6}
\newcommand{\boxwidth}{77ex}

\renewcommand{\mathcal}{\mathscr}
\newcommand{\A}{\ensuremath{\mathcal{A}}}
\newcommand{\M}{\ensuremath{\mathcal{M}}}
\renewcommand{\P}{\ensuremath{\mathcal{P}}}
\newcommand{\X}{\ensuremath{\mathcal{X}}}

\newcommand{\op}{\ensuremath{\mathrm{op}}}
\newcommand{\ob}{\ensuremath{\mathrm{ob}}}
\newcommand{\inj}{\ensuremath{\mathrm{inj}}}
\newcommand{\In}{\ensuremath{\mathrm{In}}}
\newcommand{\Out}{\ensuremath{\mathrm{Out}}}
\newcommand{\Chan}{\ensuremath{\mathrm{Chan}}}
\newcommand{\Cont}{\ensuremath{\mathrm{Cont}}}
\newcommand{\Msg}{\ensuremath{\mathbf{Msg}}}
\newcommand{\PMsg}{\ensuremath{\mathbf{PMsg}}}
\newcommand{\Map}{\ensuremath{\mathrm{Map}}}

\newcommand{\vd}{\ensuremath{\vdash}}
\newcommand{\Vd}{\ensuremath{\Vdash}}

\newcommand{\<}{\ensuremath{\langle}}
\renewcommand{\>}{\ensuremath{\rangle}}
\newcommand{\x}{\ensuremath{\times}}
\newcommand{\ox}{\ensuremath{\otimes}}

\newcommand{\parr}{\ensuremath{\oplus}}
\newcommand{\bdot}{\ensuremath{\bullet}}
\newcommand{\wdot}{\ensuremath{\circ}}
\newcommand{\0}{\ensuremath{\mathbf{0}}}
\newcommand{\init}{\ensuremath{\{\,\}}}
\newcommand{\rinit}{\ensuremath{(\,)}}
\newcommand{\ainit}{\ensuremath{\<\;\>}}
\newcommand{\sinit}{\ensuremath{[\;]}}
\newcommand{\cut}[2]{\ensuremath{~_{#1};_{#2} \,}}
\renewcommand{\hat}{\ensuremath{\widetilde}}

\newcommand{\ra}{\ensuremath{\xymatrix@1@=3.5ex{\ar[r]&}}}
\newcommand{\Ra}{\ensuremath{\xymatrix@=3.5ex{\ar@{=>}[r]&}}}
\newcommand{\pc}{\ensuremath{\xymatrix@=3.5ex{\ar@{|=|}[r]&}}}
\renewcommand{\iff}{\ensuremath{\xymatrix@1@=4ex{\ar@{<->}[r]&}}}
\newcommand{\Downpc}{\ensuremath{\xymatrix@!0@R=3.5ex{\ar@{|=|}[d] \\ ~}}}

\newcommand{\case}[4]{\ensuremath{\left\{\begin{array}{l}
    \sigma_1({#1}) \mapsto {#2} \\[1ex]
    \sigma_2({#3}) \mapsto {#4}
    \end{array}\right\}}}

\newcommand{\scase}[4]{\ensuremath{\!\left[\begin{array}{l}
    {#1} \mapsto {#2} \\[1ex]
    {#3} \mapsto {#4}
    \end{array}\right]}}

\begin{document}

\title{The logic of message passing}
\author{J. R. B. Cockett}
\address{Department of Computer Science, University
of Calgary, 2500 University Drive NW, Calgary, Alberta, Canada T2N 1N4.}
\email{robin@cpsc.ucalgary.ca}

\author{Craig Pastro}
\address{Department of Mathematics, Macquarie University,
New South Wales 2109, Australia.}
\email{craig@maths.mq.edu.au}

\date{20 September 2007}

\thanks{The first author gratefully acknowledges the support of NSERC,
Canada while the second gratefully acknowledges the support of an
international Macquarie University Research Scholarship. Parts of this work
were completed while the second author was visiting the University of Calgary
during which he was also supported by the Department of Computer Science
at the University of Calgary, Macquarie International, and the ICS
Postgraduate Research Fund.}


\begin{abstract}
Message passing is a key ingredient of concurrent programming. The purpose
of this paper is to describe the equivalence between the proof theory, the
categorical semantics, and term calculus of message passing.  In order to
achieve this we introduce the categorical notion of a linear actegory and
the related polycategorical notion of a poly-actegory.  Not surprisingly 
the notation used for the term calculus borrows heavily from the 
(synchronous) $\pi$-calculus. The cut elimination procedure for the system 
provides an operational semantics.
\end{abstract}

\maketitle

\tableofcontents

\section{Introduction} 

If programs should be viewed as proofs and types as propositions then for 
what proof system are \emph{concurrent} programs the proofs? Sequential 
programs are connected through the $\lambda$-calculus by the
Curry-Howard-Lambek isomorphism to intuitionistic proofs and cartesian
closed categories. Considering the impact that connection has had on our
understanding of sequential programs, it is reasonable to suppose that the
answer to the question for concurrent programs might have similar and far
reaching consequences. 

On the face of it the question seems beguilingly easy. Indeed, one 
might reasonably be tempted to backward engineer the $\pi$-calculus -- which 
after all holds in the concurrent world an analogous position to the 
$\lambda$-calculus -- to arrive at an answer.  However, a moment's thought 
about the passage between the $\lambda$-calculus and proofs makes one 
realise that the significant missing component, namely the type 
system, has a huge effect.  Adding types to the $\lambda$-calculus 
introduces a program discipline which, for example, is sufficiently 
restrictive to provide a guarantee of termination.  This carries the 
$\lambda$-calculus far from its freewheeling role as a description of 
computability.

Concurrent programs are even more freewheeling than sequential programs 
and, therefore, it is inevitable that collecting them into a proof system 
will introduce a programming discipline which has a similar effect of 
guaranteeing some strong formal properties (for example, being deadlock 
and livelock free). Unlike the $\lambda$-calculus, however, whose 
development was concomitant with the philosophical underpinnings of 
type theory, the $\pi$-calculus was developed in a brave new world 
of computing in which operational behaviours sufficed. Thus, the 
connection to a proof system did not seem inevitable or even necessary.

Despite this it would be a strange world indeed if there was no type 
theoretic underpinning to concurrent programming.  There is now a 
considerable body of literature connecting concurrent semantics with 
linear logic -- or at least the multiplicative fragment of it embodied 
by the logic of the polycategorical cut. While this perhaps should 
not be viewed as providing conclusive evidence that linear logic is 
the correct basis in proof theory for concurrent semantics, the argument
at this stage is 
quite compelling. See for example the work of Abramsky on interaction 
categories~\cite{A2,A3}, of Abramsky and Melli\`es~\cite{AM}, 
of Barber, Gardner, Hasegawa, and Plotkin~\cite{BGHP}, and our own 
work~\cite{CP,P}. In particular, the authors were influenced by the much
earlier paper of Bellin and Scott~\cite{BS} in which this connection was
pursued and which contains further historical commentary.

So how does one model message passing in this formal setting?  Well, of 
course there is a technical answer to that question which covers the 
pages which follow.  However, we should draw the reader's attention to one 
particular aspect of some consequence. We model message passing using a 
two tier logic.  There is a logic for the messages whose proofs should be
thought of as ordinary sequential programs. Then there is a logic 
of message passing which is built on top of the logic of messages.  The 
two logics are really quite distinct: the message logic is concerned with 
what we classically view as computation, while the second logic is concerned 
with manipulating the channels of communication. 

We believe there is a -- perhaps somewhat uncomfortable -- message in this. 
Both functional and imperative programming language designers have introduced 
concurrent features, essentially, by either adding operating system primitives 
to their sequential core or by overloading basically sequential constructs 
(consider the use of monads to obtain IO in Haskell).  However, even the 
briefest perusal of the rules, indicates that the logic concerned with 
managing channels is at least as complicated as the sequential 
programming logic. Furthermore, there are quite significant interactions 
between the two levels.  This suggests that trying to place a boundary 
to programming language design at this point is altogether artificial. 
Thus, we believe programming language designers {\em should\/} be thinking 
in terms of developing integrated two tier languages in order to give high-level 
support for concurrent programming -- operating system primitives manifestly 
fail in this regard.

In this paper we do not pretend that a logic whose proofs are concurrent 
programs is going to be a particularly simple thing.  There are many rules 
involved and consequently many equivalences between proofs. However, there is 
nothing dramatically original about the proof theory we present either. It is
basically the proof theory of (polycategorical) cut (see~\cite{BCST}) with
messages. The aim of this paper is to lay out in some detail how the 
expression of message passing is added to the logic. The resulting 
system is necessarily somewhat more complex as it has to
include the logic of the messages themselves. The logic for message passing
is then built on top of the logic for messages and embodies the interactions
which are necessary between the two levels.

In order to provide a basis which would cover a broad range of semantics we 
decided to use a logic of messages whose categorical semantics is a
distributive monoidal category. That is a monoidal category with coproducts
over which the tensor distributes. Explicitly this means that the natural map
\[
\xymatrix@C=17ex{A \ox X + A \ox Y
    \ar[r]^-{\<1_A \ox \inj_l \,,\, 1_A \ox \inj_r \>} & A \ox (X + Y)}
\]
is an isomorphism. The presence of coproducts in the messages allows us
to indicate how this structure must interact with the message passing
level. This is a rather crucial aspect of the system we present which
allows the contents of messages to determine the interaction which
actually unfolds.

In order to illustrate this point consider the following simple program
which runs a bank machine. A user will insert their card into the bank
machine and provide their personal identification number, \emph{pin}, and
a request for money, $x$. The bank machine will send this information to
the bank which will respond with a transaction identification number,
\emph{tid}, and an amount, $y$, which the bank has permitted the bank
machine to deliver to the user. This amount may be either the amount the
user has requested or, in the case the user had made a request which
cannot be satisfied (e.g., the request will put them over their daily limit
or will have them exceed their balance, etc.), zero. At this point the bank
machine may close the communication with the bank.

The bank machine will then send a request to security to do a check using
the transaction identification, \emph{tid}, supplied from the bank. The
security check will determine, for example, whether the card has been stolen,
or whether it has been used within the last few hours half-way around the
world, etc. The bank machine then receives a response, \emph{srp}, from
security indicating that the card is okay, Accept, or is not, Deny.
If the reply from security does not indicate any problem then close
communication with security and provide the user with the amount $y$. If the
reply from security indicates a problem then the machine will hold on to the
user's card and not provide any amount of money. The program will then
terminate.

We present the program using the syntax which is developed in
Section~\ref{sec-pl}. The program has type

\begin{center}
\begin{minipage}[c]{70ex}
\begin{tabbing}
\quad usr $:$ Request $\circ$ (Response $\bullet \bot$) $\Vdash$
    \= bnk $:$ Request $\circ$ (BResponse $\bullet \bot$), \\
    \> sec $:$ TransID $\circ$ (SResponse
$\bullet \bot$)
\end{tabbing}
\end{minipage}
\end{center}
and thus involves three channels labelled usr, bnk, and sec. The typing is
given as
\begin{itemize}
\item type Request = PIN $*$ Integer
\item type BResponse = TransID $*$ Integer
\item data Response = Dollar Integer $\mid$ TakeCard
\item data SResponse = Accept $\mid$ Deny
\end{itemize}

and the program is:

\begin{center}
\begin{minipage}[c]{70ex}
\begin{tabbing}
get u\=sr $(\emph{pin},x) \Rightarrow$ \\
    \> put b\=nk $(\emph{pin},x)$ \\
    \> get bnk $(\emph{tid},y) \Rightarrow$ \\
        \>\> close b\=nk \\
        \>\> put s\=ec \emph{tid} \\
        \>\> get sec \emph{srp} $\Rightarrow$ \\
            \>\>\> case \= \emph{srp} of \\
                \>\>\>\> $\mid$ Accept \= $\rightarrow$ \= close sec \\
                                      \>\>\>\>\>\> put usr (Dollar $y$) \\
                                      \>\>\>\>\>\> end usr \\
                \>\>\>\> $\mid$ Deny \> $\rightarrow$ close sec \\
                                      \>\>\>\>\>\> put usr TakeCard \\
                                      \>\>\>\>\>\> end usr
\end{tabbing}
\end{minipage}
\end{center}

The point of the example is that it shows how values received from other
processes can not only effect the values subsequently passed but also the
evolution of the communications of the whole process. From a proof theoretic
perspective this means that the coproduct structure must be shared between
the value level (the messages) and the communication level (message passing).
This significantly affects the design of the logic.

It may seem to the reader that we have made a rather esoteric choice of
logic for the messages. The choice is, in fact, minimal in order to
illustrate the interaction between the levels. Intuitively, the reader
should view messages as values of a sequential programming language (as in
the above example). However, in this paper, we have not committed ourselves
to a particular semantics for that sequential world. Thus, these values might
be from a cartesian closed category or, equally, they could be values produced
by a partial recursive function (so embody the possibility of non-termination).
In this latter case, while coproducts are present, products (in the usual
categorical sense) are not present nor is the setting closed. However,
notably, it is an example of the minimal structure we present.

Note that in the above example a rather simple use of the unit $\bot$ is
made. Those familiar with the coherence issues surrounding linearly
distributive and $*$-autonomous categories will be aware that the presence
of units adds significantly to the complexity of determining equality of
maps. Thus, it may seem sensible to avoid these units altogether in a
programming system. However, units have a crucial role as it is their
behavior which allows the proper opening and closing of channels (as was
seen in the above example).

The starting point of our exposition is a description of a logic for the
messages and a term notation to express the proofs of this logic. The term
notation for the proofs is essentially the term logic which Barry
Jay developed for monoidal categories~\cite{Jay}. This logic does not
reflect the obvious symmetry of monoidal categories obtained by reversing
maps and, perhaps for that reason, it did not resonate well with work in
monoidal categories. Here as we wish to contrast the one-sidedness
(multicategorical nature) of message logic with the two-sidedness
(polycategorical nature) of the logic for message passing, it suits our
purpose well.

In more modern terms the message logic is a multicategorical
logic with tensorial representation: we present the term logic from this
perspective (see the work of Abramsky~\cite{A1} and Mackie, Roman, and
Abramsky~\cite{MRA} for similar systems) and add a syntax for coproducts.
To present the message passing logic, we have built on the two-sided logic
presented in~\cite{CP}. In that paper a process reading for the two-sided
terms of additive linear logic was presented. There, in order to present
a two-sided notation for the proofs, we borrowed heavily from the notation
of the $\pi$-calculus. This paper continues this trend by borrowing
notation from the $\pi$-calculus in order to express the proofs associated
with message passing.

Having introduced the logic and a term calculus for its proofs, our next 
aim is to lay out the categorical semantics. We claim that this semantics 
lies in a ``linear actegory''\footnote{The term \emph{actegory} is used to
describe the situation of a monoidal category ``acting'' on a category. They
first appeared (under a different name) in the work of B\'enabou as a simple
example of a bicategory. B. Pareigis developed the theory of actegories
(again under a different name) and showed there usefulness in the
representation theory of monoids and comonoids. The word ``actegory'' was
first suggested at the Australian Category Seminar and first appeared in
print in the thesis of P. McCrudden~\cite{Mc} where they were used to study
categories of representations of coalgebroids.}, by which we mean a linearly
distributive category with a monoidal category acting on it both covariantly
and contravariantly. This structure arrives with a number of coherence
conditions which we have tried to lay out reasonably carefully. We believe
that this description of the categorical semantics of message passing is
novel. 

The final aim of the paper is to connect the term calculus, the proof
theory, and the categorical semantics. In order to make the paper more
accessible, we begin by introducing the proof theory using a sequent
calculus presentation which is annotated to provide a term calculus. However, 
the sequent calculus proof system is equivalent to a natural deduction 
system which in turn is a poly-actegory. In the last sections, 
we move rather freely between the proof theory, the poly-actegorical 
semantics, and the circuit representation of these systems.  To one 
who is not familiar with these techniques these may seem like large 
leaps as we have not tried to provide a detailed justification of it here. 
These techniques are described elsewhere and originate in Lambek's
work~\cite{Lam}. For the linear setting they are described in~\cite{BCST},
where the correctness criterion (i.e., the net condition) 
for the circuit representation of proofs is also discussed: this, although 
present in the current proof system, we barely mention.

The last sections are concerned with establishing the following
three-way equivalence:

\[
\xymatrix@C=5ex@R=5ex{
{\begin{array}{c} \text{proof theory} \\ \hline \text{poly-actegory}
 \end{array}}~~ \ar@{<->}[dr] \ar@{<->}[rr] &&
~~{\begin{array}{c} \text{category theory} \\ \hline \text{linear actegory}
\end{array}} \ar@{<->}[dl]
\\
& {\begin{array}{c} \text{term calculus} \\ \hline \text{message passing}
\end{array}}}
\]

\paragraph*{Outline of the paper.}
In Section~2 we introduce the logic of messages. Section~3 introduces the
logic of message passing which is built atop the logic of messages and
supplies the step from the proof theory to the term calculus. Also, we
introduce the cut elimination process for the logic and so, implicitly, an
operational semantics for the calculus. In Section~4 we introduce the
categorical semantics. In Section~5 we show how to obtain the categorical
semantics from the term calculus. In Section~6 we show how to move from the
categorical semantics to the poly-actegorical semantics and back using
representability. Whence, by the identification of the proof theory 
and the poly-actegorical semantics, we complete the tour of the triangle.

\section{The logic of messages} \label{sec-logic-mc}

In this section a logic for monoidal categories with coproducts is developed.
The logic is presented in a Gentzen sequent style: a sequent takes the form
\[
\Phi \vd A
\]
where the \emph{antecedent} (which we will also call the \emph{context})
of the sequent $\Phi$ is a comma separated list of formulas and the
\emph{succedent} $A$ is a single formula. It is convenient to take the
antecedent to be unordered as this allows the
permutations of the formulas without having to add an explicit exchange rule:
\begin{center}
\AxiomC{$\Phi_1,B,C,\Phi_2 \vd A$}
\LeftLabel{exchange~}
\UnaryInfC{$\Phi_1,C,B,\Phi_2 \vd A$}
\DisplayProof .
\end{center}

The inference rules for this logic are presented in
Figure~\ref{tab-mes-inference}. Notice that the cut rule here is called
``\emph{sub}'' to stand for substitution.

We will consider only the \emph{free logic} built from a multicategory.
This means that we have an arbitrary set
of atoms that will be regarded as the objects of a multicategory, and an
arbitrary set of axioms which will be regarded as the morphisms of a
multicategory. The resulting logic will be denoted by $\Msg$.

\begin{figure}
\begin{center}
\framebox[\boxwidth]{
\begin{tabular}{rlrl}
    axiom &
    \AxiomC{}
    \UnaryInfC{$\Phi \vd A$}
    \DisplayProof
&
    subs & 
    \AxiomC{$\Phi \vd A$}
    \AxiomC{$\Psi_1,A,\Psi_2 \vd B$}
    \BinaryInfC{$\Psi_1,\Phi,\Psi_2 \vd B$}
    \DisplayProof
\bigskip\\
    $*_l$ &
    \AxiomC{$\Phi,A,B \vd C$}
    \UnaryInfC{$\Phi,A * B \vd C$}
    \DisplayProof
&
    $*_r$ & 
    \AxiomC{$\Phi \vd A$}
    \AxiomC{$\Psi \vd B$}
    \BinaryInfC{$\Phi,\Psi \vd A * B$}
    \DisplayProof
\bigskip\\
    $I_l$ &
    \AxiomC{$\Phi \vd A$}
    \UnaryInfC{$\Phi,I \vd A$}
    \DisplayProof
&
    $I_r$ &
    \AxiomC{}
    \UnaryInfC{$\vd I$}
    \DisplayProof
\bigskip\\
    coprod &
    \AxiomC{$\Phi,A \vd C \hspace{-2ex}$}
    \AxiomC{$\Phi,B \vd C$}
    \BinaryInfC{$\Phi,A+B \vd C$}
    \DisplayProof
\quad & \quad
    $\text{inj}_l$ & 
    \AxiomC{$\Phi \vd A$}
    \UnaryInfC{$\Phi \vd A+B$}
    \DisplayProof
\bigskip\\
    \0 &
    \AxiomC{}
    \UnaryInfC{$\Phi,\0 \vd A$}
    \DisplayProof
&
    $\text{inj}_r$ & 
    \AxiomC{$\Phi \vd B$}
    \UnaryInfC{$\Phi \vd A+B$}
    \DisplayProof
\end{tabular}
}\end{center}
\caption{Inference rules for $\Msg$} \label{tab-mes-inference}
\end{figure}

\subsection{A term calculus for $\Msg$} \label{sec-mon-tc}

\begin{figure}
\begin{center}
\framebox[\boxwidth]{
\begin{tabular}{cc}
\!\!\!\!\!\!\!\!\!
\AxiomC{axiom $f$}
\UnaryInfC{$x_1\!:\!A_1,\ldots,x_n\!:\!A_n \!\vd\! f(x_1,\ldots,x_n)\!:\!B$}
\DisplayProof
    & 
\AxiomC{$\Phi \vd \! f : \! A$ \hspace{-2ex}}
\AxiomC{$\Psi_1,w:A,\Psi_2 \vd g \! : \! B$}
\BinaryInfC{$\Psi_1,\Phi,\Psi_2 \vd (w \mapsto g)\,f:B$}
\DisplayProof
\!\!\!\!\!\!\!\!
    \bigskip\\ 
\AxiomC{$\Phi,x:A,y:B \vd f:C$}
\UnaryInfC{$\Phi,(x,y):A * B \vd f:C$}
\DisplayProof
    &
\AxiomC{$\Phi \vd f:A$}
\AxiomC{$\Psi \vd g:B$}
\BinaryInfC{$\Phi,\Psi \vd (f,g):A * B$}
\DisplayProof
    \bigskip\\
\AxiomC{$\Phi \vd f:A$}
\UnaryInfC{$\Phi,(\,):I \vd f:A$}
\DisplayProof
    &
\AxiomC{}
\UnaryInfC{$\vd (\,):I$}
\DisplayProof
    \bigskip\\ 
\!\!\!\!\!\!
\AxiomC{$\Phi,x:A \vd f:C$ \hspace{-2ex}}
\AxiomC{$\Phi,y:B \vd g:C$}
\BinaryInfC{$\Phi,z:A+B \vd \case{x}{f}{y}{g}z: C$}
\DisplayProof
    &
\AxiomC{$\Phi \vd f:A$}
\UnaryInfC{$\Phi \vd \sigma_1(f):A+B$}
\DisplayProof
    \bigskip\\ 
\AxiomC{}
\UnaryInfC{$\Phi,z:\0 \vd \init z: A$}
\DisplayProof
    &
\AxiomC{$\Phi \vd f:B$}
\UnaryInfC{$\Phi \vd \sigma_2(f):A+B$}
\DisplayProof
\end{tabular}
}\end{center}
\caption{Term formation rules for $\Msg$} \label{tab-mes-ann}
\end{figure}

We now introduce a term calculus for this logic. The idea is that, given a
derivable sequent, to annotate the formulas on the left of the turnstile
(``\,$\vdash$\,'') with ``patterns'' made up of variables ($x,y,z,\ldots$), and
the formula on the right of the turnstile with a term ($f,g,\ldots$) which
together describe a derivation of the sequent.

Given a derivable annotated sequent $\Phi \vd A$ its annotation and
corresponding term are constructed inductively (top-down) from the
derivation. The description is given in Figure~\ref{tab-mes-ann}.

For the identity derivation on atoms $x:A \vd A$, instead of $1_A(x)$ we
will simply write $x$. That is, for atoms, the term formation rule is given
by
\begin{center}
\AxiomC{atomic $A$}
\UnaryInfC{$x:A \vd x:A$}
\DisplayProof~.
\end{center}
In order to avoid variable name clashes, an assumption that will be made is
that whenever two or more annotated sequents are involved in a derivation
(i.e., a sub, $*_r$, or coprod rule) no two will contain a variable name in
common unless mentioned explicitly.

Notice that different derivations of the same sequent will be described by the
same annotation. For example, notice that the derivations
\begin{center}
\AxiomC{}
\UnaryInfC{$x:A \vd x:A$}
\AxiomC{}
\UnaryInfC{$\vd (\,):I$}
\UnaryInfC{$(\,):I \vd (\,):I$}
\BinaryInfC{$x:A,(\,):I \vd (x,(\,)):A * I$}
\UnaryInfC{$(x,(\,)):A * I \vd (x,(\,)):A * I$}
\DisplayProof
\quad \text{and} \quad
\AxiomC{}
\UnaryInfC{$x:A \vd x:A$}
\AxiomC{}
\UnaryInfC{$\vd (\,):I$}
\BinaryInfC{$x:A \vd (x,(\,)):A * I$}
\UnaryInfC{$x:A,(\,):I \vd (x,(\,)):A * I$}
\UnaryInfC{$(x,(\,)):A * I \vd (x,(\,)):A * I$}
\DisplayProof
\end{center}
are both described by the annotation
\[
(x,(\,)):A * I \vd (x,(\,)):A * I
\]
and are therefore implicitly identified in the term calculus.

Here are two derivations of the same sequent in which the terms describing
the derivation differ.
\begin{enumerate}
\item \quad
\AxiomC{}
\UnaryInfC{$x:A \vd x:A$}
\AxiomC{}
\UnaryInfC{$y:A \vd y:A$}
\BinaryInfC{$z:A+A \vd \case{x}{x}{y}{y}z:A$}
\UnaryInfC{$z:A+A \vd \sigma_1(\case{x}{x}{y}{y}z):A+B$}
\DisplayProof
\bigskip

\item \quad
\AxiomC{}
\UnaryInfC{$x:A \vd x:A$}
\UnaryInfC{$x:A \vd \sigma_1(x):A+B$}
\AxiomC{}
\UnaryInfC{$y:A \vd y:A$}
\UnaryInfC{$y:A \vd \sigma_1(y):A+B$}
\BinaryInfC{$z:A+A \vd \case{x}{\sigma_1(x)}{y}{\sigma_1(y)}z:A+B$}
\DisplayProof
\end{enumerate}
It will be seen in Section~\ref{sec-mon-equiv} that these two terms must
be identified.

\subsection{Cut elimination for $\Msg$}

In this section the cut (``sub'') elimination rewrites are described. Recall
that, unless explicitly mentioned, a term may not contain a variable name in
common.

Notice that the $*_l$ and $I_l$ rules have no effect on the terms. In these
cases the term is actually encoding an implicit cut elimination step. For
example, both the left-hand and right-hand derivations in the cut
elimination step
\begin{center}
\AxiomC{$\Phi \vd A$}
\AxiomC{$\Psi_1,A,\Psi_2,B,C \vd $D}
\UnaryInfC{$\Psi_1,A,\Psi_2,B*C \vd D$}
\BinaryInfC{$\Psi_1,\Phi,\Psi_2,B*C \vd D$}
\DisplayProof
~$\Ra$~
\AxiomC{$\Phi \vd A$}
\AxiomC{$\Psi_1,A,\Psi_2,B,C \vd D$}
\BinaryInfC{$\Psi_1,\Phi,\Psi_2,B,C \vd D$}
\UnaryInfC{$\Psi_1,\Phi,\Psi_2,B*C \vd D$}
\DisplayProof
\end{center}
are represented by the same term.

In what follows the notation $f[y/x]$ will mean ``in $f$ substitute $y$ for
all occurrences of $x$''. Also to recover the type of a term $\Phi \vd f:A$
we denote $\Cont(f) = \Phi$ and $\Out(f) = A$.

The cut elimination rewrites are as follows.
\begin{align*}
\text{[id-sequent]} &\quad
(w \mapsto f)\,x \Ra f[x/w]
\\[2ex]
\text{[sequent-id]} &\quad
(w \mapsto w)\,g \Ra g
\\[2ex]
\text{[sequent-$\,*_r$]} &\quad
(w \mapsto (g_1,g_2))f \Ra
    \begin{cases}
    \big((w \mapsto g_1)f \,,\, g_2\big) & w \in \Cont(g_1) \smallskip\\
    \big(g_1 \,,\, (w \mapsto g_2)f \big)  & w \in \Cont(g_2)
    \end{cases}
\\[2ex]
\text{[sequent-coprod]} &\quad
\Big(w \mapsto \case{x}{g_1}{y}{g_2}z\, \Big)f
    \Ra \case{x}{(w \mapsto g_1)f}{y}{(w \mapsto g_2)f}z
\\[2ex]
\text{[coprod-sequent]} &\quad
(w \mapsto g) \left(\case{x}{f_1}{y}{f_2}z \right)
\Ra \case{x}{(w \mapsto g)f_1}{y}{(w \mapsto g)f_2}z
\\[2ex]
\text{[sequent-\0]} &\quad
(w \mapsto \init z)\,f \Ra \init z
\\[2ex]
\text{[\0-sequent]} &\quad
(w \mapsto g)(\init z) \Ra \init z
\\[2ex]
\text{[sequent-inj$_l$]} &\quad
(w \mapsto \sigma_1(g))f \Ra \sigma_1((w \mapsto g)f)
\\[2ex]
\text{[sequent-inj$_r$]} &\quad
(w \mapsto \sigma_2(g))f \Ra \sigma_2((w \mapsto g)f)
\\[2ex]
\text{[$*_r\,$-$\,*_l$]} &\quad
((x,y) \mapsto g)(f_1,f_2) \Ra (x \mapsto (y \mapsto g)f_2\,)f_1
\\[2ex]
\text{[$I_r\,$-$\,I_l$]} &\quad
((\,) \mapsto g)\, (\,) \Ra g
\\[2ex]
\text{[inj$_l$-coprod]} &\quad
\Big(z \mapsto \case{x}{g_1}{y}{g_2}z\,\Big) \sigma_1(f) \Ra (x \mapsto g_1)f
\\[2ex]
\text{[inj$_r$-coprod]} &\quad
\Big(z \mapsto \case{x}{g_1}{y}{g_2}z\,\Big)\sigma_2(f) \Ra (y \mapsto g_2)f
\end{align*}
The cut elimination procedure accounts for all the ways in which a cut can
move above a compound formula which is introduced either on the left or on
the right. Of course the cut elimination procedure will get stuck on the
atomic cuts (composition) in the multicategory. However, it is easy to check
that, if composition terminates in the underlying multicategory, this
process will terminate. Indeed, in terms with only primitive function
symbols (no axioms) which only involve primitive types (no atoms) the cuts
can be completely eliminated.

\subsection{Equations in $\Msg$} \label{sec-mon-equiv}

In order to ensure that the cut elimination procedure is confluent
identities (for which we use the notation ``\pc'') between cut eliminated
terms need to be introduced. Firstly, if $f$, $g$, and $h$ are axioms,
equations are needed describing the associative law and interchange law.
For associativity suppose
\[
\Phi \vd f:A, \quad \Psi_1,x:A,\Psi_2 \vd g:B, \quad \text{and} \quad
\Psi'_1,y:B,\Psi'_2 \vd h:B.
\]
for which the identity
\[
(y \mapsto h)((x \mapsto g)f) \pc (x \mapsto (y \mapsto h)g)f
\]
describing associativity must be added. Similarly, for the interchange law
suppose
\[
\Phi \vd f:A, \quad \Phi' \vd g:B, \quad \text{and} \quad
\Psi_1,x:A,\Psi_2,y:B,\Psi_3 \vd h:C.
\]
The identity 
\[
(y \mapsto (x \mapsto h)f)g \pc (x \mapsto (y \mapsto h)g)f
\]
describes the interchange law.

We now move on to examining the compound terms. Here is an example which
shows how such an identity arises. There are two ways to cut eliminate the
top term in the following diagram.
\[
\xymatrix@!0@C=20ex@R=12ex{
& (w \mapsto (g_1,g_2))\left(\case{x_1}{f_1}{x_2}{f_2}z\right)
\ar@{=>}[dr] \ar@{=>}[dl]_{w \in g_1~} \\
    \left((w \mapsto g_1) (\case{x_1}{f_1}{x_2}{f_2}z) ~,~ g_2\right)
\ar@{=>}[d] &&
    \case{x_1}{(w \mapsto (g_1,g_2))f_1}{x_2}{(w \mapsto (g_1,g_2))f_2}z
\ar@{=>}[d] \\
    \left(\case{x_1}{(w \mapsto g_1)f_1}{x_2}{(w \mapsto g_1)f_2}z
    ~,~ g_2\right)
\ar@{|=|}[rr] &&
    \case{x_1}{((w \mapsto g_1)f_1,g_2)}
         {x_2}{((w \mapsto g_1)f_2,g_2)}z
}
\]
and this forces the identity at the bottom of the diagram. Similarly, the
identification of the two terms at the end of Section~\ref{sec-mon-tc}
arises from the following diagram.
\[
\xymatrix@!0@C=20ex@R=12ex{
& {(w \mapsto \sigma_1(f))\left(\case{x_1}{g_1}{x_2}{g_2}z\right)}
\ar@{=>}[dr] \ar@{=>}[dl] \\
    \case{x}{(w \mapsto \sigma_1(f))g_1}{y}{(w \mapsto \sigma_1(f))g_2}z
\ar@{=>}[d] &&
    \sigma_1\left((w \mapsto f)\left(\case{x_1}{g_1}{x_2}{g_2}z\right)\right)
\ar@{=>}[d] \\
    \case{x}{\sigma_1((w \mapsto f)g_1)}{y}{\sigma_1((w \mapsto f)g_2)}z\
\ar@{|=|}[rr] &&
    ~\sigma_1\left(\case{x_1}{(w \mapsto f)g_1}{x_2}{(w \mapsto f)g_2}z\right)
}
\]

The list of identities which are introduced into the system in this manner
are presented in Figure~\ref{fig-id-mon}.

\begin{figure}
\begin{tabular}{c}
\\ \hline
\hspace{74ex}
\end{tabular}
\begin{enumerate}
\item $(y \mapsto h)((x \mapsto g)f) \pc (x \mapsto (y \mapsto h)g)f$
    \qquad $f,g,h$ atomic
    \medskip
\item $(y \mapsto (x \mapsto h)f)g \pc (x \mapsto (y \mapsto h)g)f$
    \qquad $f,g,h$ atomic
    \medskip
\item $\left(\case{x}{f_1}{y}{f_2}z,g \right) \pc
    \case{x}{(f_1,g)}{y}{(f_2,g)}z$
    \medskip
\item $\left(f,\case{x}{g_1}{y}{g_2}z \right) \pc
    \case{x}{(f,g_1)}{y}{(f,g_2)}z$
    \medskip
\item $\case{x_1}{\case{y_1}{f_1}{y_2}{f_2}z_2}
     {x_2}{\case{y_1}{g_1}{y_2}{g_2}z_2}z_1 \pc
            \case{y_1}{\case{x_1}{f_1}{x_2}{g_1}z_1}
            {y_2}{\case{x_1}{f_2}{x_2}{g_2}z_1}z_2$
    \medskip
\item $\sigma_1\left(\case{x}{f}{y}{g}z\right) \pc
    \case{x}{\sigma_1(f)}{y}{\sigma_1(g)}z$
    \medskip
\item $\sigma_2\left(\case{x}{f}{y}{g}z\right) \pc
    \case{x}{\sigma_2(f)}{y}{\sigma_2(g)}z$
    \medskip
\item $(f,\init z) \pc \init z$
    \medskip
\item $(\init z,g) \pc \init z$
    \medskip
\item $ \case{x}{\init z}{y}{\init z} \pc \init z$
    \medskip
\item $ \sigma_1(\init z) \pc \init z$
    \medskip
\item $ \sigma_2(\init z) \pc \init z$
\end{enumerate}
\begin{tabular}{c}
\\ \hline
\hspace{74ex}
\end{tabular}
\caption{Identities in $\Msg$}
\label{fig-id-mon}
\end{figure}

\begin{example}
The cut elimination procedure allows us to prove the distributive law. This
involves proving that the composite of
\[
w:A * x:(B+C) \xymatrix@C=27ex{\ar[r]_{
    \case{x_1}{\sigma_1((w,x_1))}{x_2}{\sigma_2((w,x_2))}x} &}
    A*B + A*C
\]
and
\[
z:((y_1,y_2):A *B+ (z_1,z_2):A*C) \xymatrix@C=31ex{\ar[r]_{
\case{(y_1,y_2)}{(y_1,\sigma_1(y_2))}{(z_1,z_2)}{(z_1,\sigma_2(z_2))}z} &}
    A*(B+C),
\]
and its reverse, are the identity. The above composite gives:
\newpage
\begin{gather*}
\left(z \mapsto \case{(y_1,y_2)}{(y_1,\sigma_1(y_2))}
                     {(z_1,z_2)}{(z_1,\sigma_2(z_2))}z\right)
    \left(\case{x_1}{\sigma_1((w,x_1))}{x_2}{\sigma_2((w,x_2))}x \right)
\\[1ex] \Downarrow [\textup{coprod-sequent}] \\[1ex]
    \case{x_1}{\left(z \mapsto
            \case{(y_1,y_2)}{(y_1,\sigma_1(y_2))}
                 {(z_1,z_2)}{(z_1,\sigma_2(z_2))}z\right)\sigma_1((w,x_1))}
         {x_2}{\left(z \mapsto
            \case{(y_1,y_2)}{(y_1,\sigma_1(y_2))}
                 {(z_1,z_2)}{(z_1,\sigma_2(z_2))}z\right)\sigma_2((w,x_2))}x
\\[1ex] \Downarrow [\textup{inj-coprod}] \\[1ex]
    \case{x_1}{\big((y_1,y_2) \mapsto (y_1,\sigma_1(y_2))\big)(w,x_1)}
         {x_2}{\big((z_1,z_2) \mapsto (z_1,\sigma_2(z_2))\big)(w,x_2)}x
\\[1ex] \Downarrow [*_r\text{-}*_l] \\[1ex]
    \case{x_1}{(y_1 \mapsto (y_2 \mapsto (y_1,\sigma_1(y_2)))x_1)w}
         {x_2}{(z_1 \mapsto (z_2 \mapsto (z_1,\sigma_2(z_2)))x_2)w}x
\\[1ex] \Downarrow [\textup{id-sequent}] \\[1ex]
    \case{x_1}{(y_2 \mapsto (w,\sigma_1(y_2)))x_1}
         {x_2}{(z_2 \mapsto (w,\sigma_2(z_2)))x_2}x
\\[1ex] \Downarrow [\textup{id-sequent}] \\[1ex]
    \case{x_1}{(w,\sigma_1(x_1))}{x_2}{(w,\sigma_2(x_2))}x \\
\end{gather*}
where this last is (one form of) the identity map of $A*(B+C)$.
The other way around also works (giving the identity on $A*B + A*C$) and is
left for the interested reader to familiarise themselves with this calculus.
\end{example}

\section{The logic of message passing} \label{sec-logic-mp}

This section introduces the sequent rules for the message passing logic
which will be denoted by $\PMsg$. As the message passing logic is built on
top of the logic of messages -- which in this case we are
taking to be $\Msg$ -- the logic involves inference rules whose premisses
are inferences of both systems. In order to help the reader keep this
straight we shall use two different entailment symbols: ``\,$\vdash$\,''
for the messages themselves and ``\,$\Vdash$\,'' for the message passing
logic. For example, the inference rule $\bdot_l$ has the form
\[
\AxiomC{$\pi$}
\UnaryInfC{$\Phi \vd A$}
\AxiomC{$\Pi$}
\UnaryInfC{$\Psi \mid \Gamma, X \Vd \Delta$}
\BinaryInfC{$\Phi,\Psi \mid \Gamma, A \bdot X \Vd \Delta$}
\DisplayProof
\]
where $\pi$ denotes a derivation in $\Msg$ and $\Pi$ a derivation of $\PMsg$.

Semantically the message passing logic builds a linearly distributive
category (which, when linear adjoints are present, is just a $*$-autonomous 
category~\cite{CKS}) from the underlying message logic. This linearly
distributive category, as we shall see, is part of a linear actegory. The
term calculus, which we construct in the next section, then becomes a very
basic language for concurrent programs which can pass as values the messages
provided by our message logic. 

A sequent of $\PMsg$ has three components: the message type context $\Phi$,
its input message passing types $\Gamma$, and its output message passing
types $\Delta$ which together define a sequent of $\PMsg$:
\[
\Phi \mid \Gamma \Vdash \Delta.
\]
We shall treat all three components as unordered lists.  

The axioms for the message passing logic should be regarded as being maps in a
poly-actegory. A poly-actegory (see Section~\ref{sec-rep}) is a
symmetric polycategory whose components have certain inputs from a
multicategory. Intuitively one may think of such an axiom as a process
between certain channels which is parameterised by certain values (from
the sequential world). The usual associativity and interchange laws hold
for this (multi and) polycomposition.

The inference rules for $\PMsg$ are presented in Figure~\ref{message-infrules}.
Our convention will be to denote formulas from the message logic $\Msg$
using uppercase letters from the beginning of the alphabet $A,B,\ldots$ and
unordered lists of these formulas using $\Phi$ and $\Psi$. Formulas from the
message passing logic $\PMsg$ will be denoted using uppercase letters from
the end of the alphabet $X,Y,\ldots$ and unordered lists of these formulas
using $\Gamma$ and $\Delta$.

Most of the rules of this calculus are just the standard ones for two-sided
multiplicative linear logic. The main novel aspects are the ``action'' rules
$\circ_l$, $\circ_r$, $\bdot_l$, and $\bdot_r$, which allow messages to be
bound to channels of interaction in the message passing logic $\PMsg$.
However, also notice that the effect of a sum of messages can be derived
from how they are passed in $\PMsg$.

\begin{figure}
\begin{center}
\begin{tabular}{c}
\\ \hline
\hspace{74ex}
\end{tabular}
\begin{tabular}{rl}
cut \hspace{-2ex} & \hspace{-2ex} 
\AxiomC{$\Phi \mid \Gamma_1 \Vd \Delta_1,X$}
\AxiomC{$\hspace{-13pt} \Psi \mid X,\Gamma_2 \Vd \Delta_2$}
\BinaryInfC{$\Phi,\Psi \mid \Gamma_1,\Gamma_2 \Vd \Delta_1,\Delta_2$}
\DisplayProof
    \bigskip
\end{tabular}
\begin{tabular}{rlrl}
\hspace{-2ex} atom id \hspace{-2ex} & \hspace{-2ex}
\AxiomC{}
\UnaryInfC{$\emptyset \mid X \vd X$}
\DisplayProof
    &
\hspace{-3ex} axiom \hspace{-2ex} & \hspace{-2ex}
\AxiomC{}
\UnaryInfC{$\Phi \mid \Gamma \Vd \Delta$}
\DisplayProof
    \bigskip\\
$\ox_l$ \hspace{-2ex} & \hspace{-2ex}
\AxiomC{$\Phi \mid \Gamma,X,Y \Vd \Delta$}
\UnaryInfC{$\Phi \mid \Gamma,X \ox Y \Vd \Delta$}
\DisplayProof
    &
$\ox_r$ \hspace{-2ex} & \hspace{-2ex}
\AxiomC{$\Phi\! \mid\! \Gamma_1 \Vd \Delta_1,X$}
\AxiomC{$\hspace{-13pt} \Psi\! \mid\! \Gamma_2 \Vd Y,\Delta_2$}
\BinaryInfC{$\Phi,\Psi\! \mid\! \Gamma_1,\Gamma_2 \Vd \Delta_1,X \ox Y,\Delta_2$}
\DisplayProof
    \bigskip\\
$\parr_l$ \hspace{-2ex} & \hspace{-2ex}
\AxiomC{$\Phi\! \mid\! \Gamma_1,X \Vd \Delta_1$}
\AxiomC{$\hspace{-13pt} \Psi\! \mid\! Y,\Gamma_2 \Vd \Delta_2$}
\BinaryInfC{$\Phi,\Psi\! \mid\! \Gamma_1,X \parr Y,\Gamma_2 \Vd \Delta_1,\Delta_2$}
\DisplayProof
    \hspace{-2.5ex} & \hspace{-2.5ex}
$\parr_r$ \hspace{-2ex} & \hspace{-2ex}
\AxiomC{$\Phi \mid \Gamma \Vd X,Y,\Delta$}
\UnaryInfC{$\Phi \mid \Gamma \Vd X \parr Y,\Delta$}
\DisplayProof
    \bigskip\\
$\top_l$ \hspace{-2ex} & \hspace{-2ex}
\AxiomC{$\Phi \mid \Gamma \Vd \Delta$}
\UnaryInfC{$\Phi \mid \Gamma,\top \Vd \Delta$}
\DisplayProof
    &
$\top_r$ \hspace{-2ex} & \hspace{-2ex}
\AxiomC{}
\UnaryInfC{$\emptyset \mid~ \Vd \top$}
\DisplayProof
    \bigskip\\
$\bot_l$ \hspace{-2ex} & \hspace{-2ex}
\AxiomC{}
\UnaryInfC{$\emptyset \mid \bot \Vd $}
\DisplayProof
    &
$\bot_r$ \hspace{-2ex} & \hspace{-2ex}
\AxiomC{$\Phi \mid \Gamma \Vd \Delta$}
\UnaryInfC{$\Phi \mid \Gamma \Vd \bot,\Delta$}
\DisplayProof
    \bigskip\\
$\wdot_l$ \hspace{-2ex} & \hspace{-2ex}
\AxiomC{$\Phi,A \mid \Gamma, X \Vd \Delta$}
\UnaryInfC{$\Phi \mid \Gamma, A \wdot X \Vd \Delta$}
\DisplayProof
    &
$\wdot_r$ \hspace{-2ex} & \hspace{-2ex}
\AxiomC{$\Phi \vd A$}
\AxiomC{$\Psi \mid \Gamma \Vd X,\Delta$}
\BinaryInfC{$\Phi,\Psi \mid \Gamma \Vd A \wdot X,\Delta$}
\DisplayProof
    \bigskip\\
$\bdot_l$ \hspace{-2ex} & \hspace{-2ex}
\AxiomC{$\Phi \vd A$}
\AxiomC{$\Psi \mid \Gamma, X \Vd \Delta$}
\BinaryInfC{$\Phi,\Psi \mid \Gamma, A \bdot X \Vd \Delta$}
\DisplayProof
    &
$\bdot_r$ \hspace{-2ex} & \hspace{-2ex}
\AxiomC{$\Phi,A \mid \Gamma \Vd X,\Delta$}
\UnaryInfC{$\Phi \mid \Gamma \Vd A \bdot X,\Delta$}
\DisplayProof
    \bigskip\\
$*$ \hspace{-2ex} & \hspace{-2ex}
\AxiomC{$\Phi,A,B \mid \Gamma \Vd \Delta$}
\UnaryInfC{$\Phi,A*B \mid \Gamma \Vd \Delta$}
\DisplayProof
   ~&~ 
$I$ \hspace{-2ex} & \hspace{-2ex}
\AxiomC{$\Phi \mid \Gamma \Vd \Delta$}
\UnaryInfC{$\Phi,I \mid \Gamma \Vd \Delta$}
\DisplayProof
\bigskip\\
coprod \hspace{-2ex} & \hspace{-2ex}
\AxiomC{$\Phi,A \mid \Gamma \Vd \Delta \hspace{-1ex}$}
\AxiomC{$\Phi,B \mid \Gamma \Vd \Delta \hspace{-1ex}$}
\BinaryInfC{$\Phi,A+B \mid \Gamma \Vd \Delta$}
\DisplayProof
    \quad & \quad
$\0$ \hspace{-2ex} & \hspace{-2ex}
\AxiomC{}
\UnaryInfC{$\Phi,\0 \mid \Gamma \Vd \Delta$}
\DisplayProof
\end{tabular}
    \bigskip\\
\begin{tabular}{rl}
subs \hspace{-2ex} & \hspace{-2ex}
\AxiomC{$\Psi \vd A$}
\AxiomC{$\Phi,A \mid \Gamma \Vd \Delta$}
\BinaryInfC{$\Phi,\Psi \mid \Gamma \Vd \Delta$}
\DisplayProof
\end{tabular}

\begin{tabular}{c}
\hspace{74ex}
\\ \hline
\end{tabular}
\end{center}
\caption{Inference rules for $\PMsg$} \label{message-infrules}
\end{figure}

\subsection{Term calculus for $\PMsg$} \label{sec-mess-term}

The term calculus we now introduce for message passing should be
thought of as a very basic language for concurrency which permits
point-to-point interactions along channels. In order to refer to the
individual formulas of a sequent, which are to be thought of as the
channels through which the process embodied by the sequent interacts, we
shall label them with ``channel names'' using lowercase Greek letters. For
example
\[
(x,y):A*B,z:C ~\mid~ \alpha:W \ox X \Vdash \beta:D \bdot Y,~ \gamma:Z.
\]

In the $\pi$-calculus much emphasis is lain on how these channel names are
propagated. In particular the ability to pass channel names as messages
introduces scope extrusion and the necessity for channel relabeling. The
calculus we present is not as free ranging and does not allow the passing
of channel names as messages. In particular, here we distinguish sharply
between the world of messages and the mechanisms for message passing. It
might, therefore, be supposed that the ability to pass channel names is
completely absent from this system. However, this is not the case. Although
we have not chosen to concentrate on these issues (or abilities), in fact
they are already present in the message passing calculus. In the
$\pi$-calculus message passing is the \emph{only} mechanism present and so
passing of channel names has to be achieved by passing them as messages. In
the current calculus, significantly, there are other mechanisms present, in
particular, one can bundle channels together (using the tensor or the par)
and, thus, one can pass simultaneously on a channel multiple channel names
along which the receiving process can subsequently interact. The issue of
scope extrusion is actually handled in the cut elimination procedure which,
in effect, also defines the operational semantics of the system.
Particularly relevant in this regard is the $\ox_r-\ox_l$ cut elimination
step which shows how a process can use channel names which are passed to it.

Modern process calculi are also concerned with the issue of "mobility":
this means both computation carried out on mobile devices (i.e., networks
that have a dynamic topology), and mobile computation (i.e., executable
code that is able to move around a network). For example, the ambient
calculus~\cite{CG} of L. Cardelli and A. Gordon was introduced to addresses
these issues: ambients being conceptual locations in which computation can
occur. The calculus we have presented is not intended to address these
issues and, indeed, is completely neutral on its ``ambient'' implementation.
It is, of course, a pertinent issue of how to model ambient calculi
categorically and proof theoretically: this may provide a useful mathematical
and logical insight into these calculi.

As discussed above, the cut elimination procedure forces the behavior of
scope in our calculus and, thus, the channel renaming which is necessary.
As we are using a two-sided calculus it is possible to have completely
separate name spaces for input and output channels. The ``plugging
together'' of processes on a channel -- which is a cut --- binds an output
channel of one process $s$ to an input channel of another $t$ and, thus,
may be denoted by an infix syntax:
\[
s ~_\alpha;_\beta t
\]
where we assume $s$ and $t$ have distinct output channel names and distinct
input channel names. We will also allow the use of the simpler
\[
s ~;_\alpha t
\]
where $\alpha$ is both an output channel bound in $s$ and an input channel
bound in $t$. Once a channel name is bound it can be renamed to
any unused name and, indeed, this may be required to reassociate cuts:
\[
(s ~_\alpha;_{\alpha'} t) ~_\beta;_{\beta'} u
    = s ~_\alpha;_{\alpha'} (t ~_\beta;_{\beta'} u).
\]
Here this equality is only valid without renaming if $\alpha' \not\in
\In(u)$ and $\beta \not\in \Out(s)$, although with renaming this
associativity is always valid.

We now describe the term formation rules. The terms presented here make use
of the self-dual nature of the logic. That is, term formation rules for dual
inference rules (e.g., $\ox_l$ and $\oplus_r$) will be identical. This
agrees with the process reading of the rules when there is no distinction
made between being an input channel and an output channel.

The notation ``\,::\,'' is used to denote the term-type membership relation,
e.g., $s:: \Phi \mid \Gamma \Vd \Delta$ means that $s$ is a term
of type $\Phi \mid \Gamma \Vd \Delta$. The lengthy syntax will not permit us
to present these rules in a table, and so we do so as a list. However, a
summary of the term formation rules is provided in Figure~\ref{tab-summary}.
\bigskip

\begin{tabular}{ll}
[cut] &
    \AxiomC{$s :: \Phi \mid \Gamma_1 \Vd \Delta_1,\alpha:X$}
    \AxiomC{$t :: \Psi \mid \beta:X,\Gamma_2 \Vd \Delta_2$}
    \BinaryInfC{$s \cut{\alpha}{\beta} t ::
                \Phi,\Psi \mid \Gamma_1,\Gamma_2 \Vd \Delta_1,\Delta_2$}
    \DisplayProof
\end{tabular}
\bigskip

\begin{tabular}{ll}
[atomic identity] &
\AxiomC{}
\UnaryInfC{$\alpha =_X \beta :: \emptyset \mid \alpha:X \Vd \beta:Y$}
\DisplayProof
\end{tabular}
\bigskip

\begin{tabular}{ll}
[axiom] &
\AxiomC{$s(\Phi)[\Gamma;\Delta]$ primitive process}
\UnaryInfC{$s(x_1,\ldots,x_r)[\alpha_1,\ldots,\alpha_m;
                                    \beta_1,\ldots,\beta_n] 
                 :: x:\Phi \mid \alpha:\Gamma \Vd \beta:\Delta$}
\DisplayProof \\ ~ & ~ \\
& where $\alpha:\Gamma$ stands for $\alpha_1:A_1,\ldots,\alpha_m:A_m$ etc.
\end{tabular}
\bigskip

\begin{tabular}{ll}
[$\ox_l$ and $\parr_r$] &
    \Axiom$s \fCenter \Phi \mid \Gamma,\alpha_1:X,\alpha_2:Y \Vd \Delta$
    \UnaryInf$\alpha\<\alpha_1,\alpha_2\> \cdot s \fCenter
            \Phi \mid \Gamma,\alpha:X \ox Y \Vd \Delta$
    \DisplayProof
    \quad and dually \bigskip\\
&
    \Axiom$s \fCenter \Phi \mid \Gamma \Vd \alpha_1:X,\alpha_2:Y,\Delta$
    \UnaryInf$\alpha\<\alpha_1,\alpha_2\> \cdot s \fCenter
            \Phi \mid \Gamma \Vd \alpha:X \parr Y,\Delta$
    \DisplayProof
\end{tabular}
\bigskip

\begin{tabular}{ll}
[$\parr_l$ and $\ox_r$] &
    \AxiomC{$s_1::\Phi \mid \Gamma_1,\alpha_1:X \Vd \Delta_1$}
    \AxiomC{$s_2::\Psi \mid \alpha_2:Y,\Gamma_2 \Vd \Delta_2$}
    \BinaryInfC{$\alpha \scase{\alpha_1}{s_1}{\alpha_2}{s_2} ::
      \Phi,\Psi \mid \Gamma_1,\alpha:X \parr Y,\Gamma_2 \Vd \Delta_1,\Delta_2$}
    \DisplayProof
    \bigskip\\
and dually &
    \AxiomC{$s_1::\Phi \mid \Gamma_1 \Vd \Delta_1,\alpha_1:X$}
    \AxiomC{$s_2::\Psi \mid \Gamma_2 \Vd \alpha_2:Y,\Delta_2$}
    \BinaryInfC{$\alpha \scase{\alpha_1}{s_1}{\alpha_2}{s_2} ::
       \Phi,\Psi \mid \Gamma_1,\Gamma_2 \Vd \Delta_1,\alpha:X \ox Y,\Delta_2$}
    \DisplayProof
\end{tabular}
\bigskip

\begin{tabular}{ll}
[$\top_l$ and $\bot_r$] &
    \Axiom$s \fCenter \Phi \mid \Gamma \Vd \Delta$
    \UnaryInf$\alpha\ainit \cdot s \fCenter
             \Phi \mid \Gamma,\alpha:\top \Vd \Delta$
    \DisplayProof \quad and dually \bigskip\\
&
    \Axiom$s \fCenter \Phi \mid \Gamma \Vd \Delta$
    \UnaryInf$\alpha\ainit \cdot s
             \fCenter \Phi \mid \Gamma \Vd \alpha:\bot,\Delta$
    \DisplayProof
\end{tabular}
\bigskip

\begin{tabular}{ll}
[$\bot_l$ and $\top_r$] &
    \AxiomC{}
    \UnaryInfC{$\alpha\sinit :: \emptyset \mid \alpha:\bot \Vd $}
    \DisplayProof
    \quad and dually \bigskip\\
&
    \AxiomC{}
    \UnaryInfC{$\alpha\sinit :: \emptyset \mid~ \Vd \alpha:\top $}
    \DisplayProof
\end{tabular}
\bigskip

\begin{tabular}{ll}
[$\wdot_l$ and $\bdot_r$] &
    \Axiom$s \fCenter x:A,\Phi \mid \Gamma, \alpha:X \Vd \Delta$
    \UnaryInf$\alpha\<x\> \cdot s \fCenter
         \Phi \mid \Gamma,\alpha:A \wdot X \Vd \Delta$
    \DisplayProof \quad and dually
    \bigskip\\
&
    \Axiom$s \fCenter x:A,\Phi \mid \Gamma \Vd \alpha:X,\Delta$
    \UnaryInf$\alpha\<x\> \cdot s \fCenter
         \Phi \mid \Gamma \Vd \alpha:A \bdot X,\Delta$
    \DisplayProof
\end{tabular}
\bigskip

A process reading of the terms for the [$\wdot_l$ and $\bdot_r$] inferences
may be thought of as follows: read $x$ on channel $\alpha$ and bind it in
the process $s$.
\bigskip

\begin{tabular}{ll}
[$\bdot_l$ and $\wdot_r$] &
    \AxiomC{$\Phi \vd f:A$}
    \AxiomC{$s :: \Psi \mid \Gamma, \alpha:X \Vd \Delta$}
    \BinaryInfC{$\alpha[f] \cdot s::
             \Phi,\Psi \mid \Gamma, \alpha:A \bdot X \Vd \Delta$}
    \DisplayProof \quad and dually
    \bigskip\\
&
    \AxiomC{$\Phi \vd f:A$}
    \AxiomC{$s::\Psi \mid \Gamma \Vd \alpha:X,\Delta$}
    \BinaryInfC{$\alpha[f] \cdot s::
             \Phi,\Psi \mid \Gamma \Vd \alpha:A \wdot X,\Delta$}
    \DisplayProof
\end{tabular}
\bigskip

A process reading of the terms for the [$\bdot_l$ and $\wdot_r$]
inferences may be thought of as follows: output $f$ on channel $\alpha$ and
continue with the process $s$.
\bigskip

\begin{tabular}{rl}
    [$*$] & 
     \Axiom$s \fCenter \Phi, x:A,y:B\mid \Gamma \Vd \Delta$
     \UnaryInf$s \fCenter \Phi, (x,y):A * B \mid \Gamma \Vd \Delta$
     \DisplayProof
\\[4ex]
    [$I$] &
    \Axiom$s \fCenter \Phi \mid \Gamma \Vd \Delta$
    \UnaryInf$s \fCenter \Phi,\rinit:I \mid \Gamma \Vd \Delta$
    \DisplayProof
\\[4ex]
    [coprod] &
    \AxiomC{$s:: \Phi,x:A \mid \Gamma \Vd \Delta$}
    \AxiomC{$t:: \Phi,y:B \mid \Gamma \Vd \Delta$}
    \BinaryInfC{$\case{x}{s}{y}{t}z :: \Phi,z:A+B \mid \Gamma \Vd \Delta$}
    \DisplayProof
\\[7ex]
    [$\0$] &
    \AxiomC{}
    \UnaryInfC{$\init z :: \Phi,z:\0 \mid \Gamma \Vd \Delta$}
    \DisplayProof
\\[4ex]
    [substitution] &
    \AxiomC{$\Phi \vd f:A$}
    \AxiomC{$s :: \Psi,x:A \mid \Gamma \Vd \Delta$}
    \BinaryInfC{$(x \mapsto s)f :: \Phi,\Psi \mid \Gamma \Vd \Delta$}
    \DisplayProof
\end{tabular}

\begin{figure}
\begin{center}
\framebox[\boxwidth]{
\begin{tabular}{rlrl}
cut &
\AxiomC{$s$}
\AxiomC{$t$}
\BinaryInfC{$s \cut{\alpha}{\beta} t$}
\DisplayProof
    \bigskip & & \\
atomic id & 
\AxiomC{}
\UnaryInfC{$\alpha =_X \beta$}
\DisplayProof
&
axiom &
\AxiomC{}
\UnaryInfC{$s(\Phi)[\Gamma;\Delta]$}
\DisplayProof
    \bigskip\\
$\ox_l$ and $\parr_r$
&
\AxiomC{$s$}
\UnaryInfC{$\alpha\<\alpha_1,\alpha_2\> \cdot s$}
\DisplayProof
\quad & \quad
$\parr_l$ and $\ox_r$ &
\AxiomC{$s_1$}
\AxiomC{$s_2$}
\BinaryInfC{$\alpha\scase{\alpha_1}{s_1}{\alpha_2}{s_2}$}
\DisplayProof
    \bigskip\\
$\top_l$ and $\bot_r$ &
\AxiomC{$s$}
\UnaryInfC{$\alpha\ainit \cdot s$}
\DisplayProof
&
$\bot_l$ and $\top_r$ &
\AxiomC{}
\UnaryInfC{$\alpha\sinit$}
\DisplayProof
    \bigskip\\
$\wdot_l$ and $\bdot_r$ &
\AxiomC{$s$}
\UnaryInfC{$\alpha\<x\> \cdot s$}
\DisplayProof
&
$\bdot_l$ and $\wdot_r$ &
\AxiomC{$f$}
\AxiomC{$s$}
\BinaryInfC{$\alpha[f] \cdot s$}
\DisplayProof
    \bigskip\\
$*$ &
\AxiomC{$s$}
\UnaryInfC{$s$}
\DisplayProof
&
$I$ &
\AxiomC{$s$}
\UnaryInfC{$s$}
\DisplayProof
    \bigskip\\
coprod &
\AxiomC{$s$}
\AxiomC{$t$}
\BinaryInfC{$\case{x}{s}{y}{t}z$}
\DisplayProof
&
$\0$ &
\AxiomC{}
\UnaryInfC{$\init z$}
\DisplayProof
    \bigskip\\
subs &
\AxiomC{$f$}
\AxiomC{$s$}
\BinaryInfC{$(x \mapsto s)f$}
\DisplayProof
&
\end{tabular}
}\end{center}
\caption{Summary of the term formation rules for $\PMsg$}\label{tab-summary}
\end{figure}

\subsection{A programming syntax term calculus}\label{sec-pl}

The term calculi for $\PMsg$ presented in the previous section is quite
useful for manipulations of the logic, but it does not illustrate well the
programming view of those proofs. To illustrate this relationship, and to
connect with the example in the introduction, we present a programming
syntax for $\PMsg$ similar to the syntax presented in~\cite{CP}.

Coproducts in programming languages are usually introduced as
(non-recursive) datatypes and so it is useful to show how this may be
incorporated into the message logic. A non-recursive datatype is defined as
\begin{center}
data $F(A_1,\ldots,A_n)$ = $C_1 T_1 \mid \cdots \mid C_r T_r$
\end{center}
where $T_i$ is a type expression in variables $A_1,\ldots,A_n$. Note that,
in contrast to the message logic in which coproducts were defined using a
binary rule and a nullary rule, here we are defining coproducts indexed by an
arbitrary finite set of constructors. This means in particular that this
definition also captures nullary and unary coproducts.

Two rules are needed to introduce the syntax associated with this datatype:
\begin{center}
\begin{tabular}{rl}
    coprod &
    \AxiomC{$\Phi,x_1:T_1 \vd f_1:B$}
    \AxiomC{$\cdots$}
    \AxiomC{$\Phi,x_r:T_r \vd f_r:B$}
    \TrinaryInfC{
    \begin{minipage}[t]{10ex}
    \begin{tabbing}
    $\Phi,z:F(A_1,\ldots,A_n) \vd$ case $z$ as
                \= $\mid$ \= $C_1 x_1 \rightarrow f_1 ~~:B$ \\
                        \>\> $\cdots$ \\
                \> $\mid$ $C_r x_r \rightarrow f_r$
    \end{tabbing}
    \end{minipage}}
    \DisplayProof
\\
    construct &
    \AxiomC{$\Phi \vd f:T_i$}
    \UnaryInfC{$\Phi \vd C_i f:F(A_1,\ldots,A_n)$}
    \DisplayProof~.
\end{tabular}
\end{center}

A convenient programming syntax for substitution is given by:
\begin{center}
\begin{tabular}{rl}
    subs &
    \AxiomC{$\Phi \vd \! f : \! A$ \hspace{-2ex}}
    \AxiomC{$\Psi_1,w:A,\Psi_2 \vd g \! : \! B$}
    \BinaryInfC{$\Psi_1,\Phi,\Psi_2 \vd$ $g$ where $w = f:B$}
    \DisplayProof~.
\end{tabular}
\end{center}

A programming syntax for $\PMsg$ is given in Figure~\ref{pmsg-prog-syn}.
The only novel aspect here is the syntax for the action rules and the rules
for the interaction between the two logics.

\begin{figure}
\begin{center}
\framebox[78ex][c]{
\begin{tabular}{rlrl}
cut &
\AxiomC{$s$}
\AxiomC{$t$}
\BinaryInfC{on $\alpha=\beta$ plug $s$ to $t$}
\DisplayProof
    \bigskip
& & \\
atomic id & 
\AxiomC{}
\UnaryInfC{$\alpha =_X \beta$}
\DisplayProof
&
axiom &
\AxiomC{}
\UnaryInfC{$s(\Phi)[\Gamma;\Delta]$}
\DisplayProof
    \bigskip\\
$\ox_l$ and $\oplus_r$
&
\AxiomC{$s$}
\UnaryInfC{split $\alpha$ as $\alpha_1,\alpha_2$; $s$}
\DisplayProof
&
$\oplus_l$ and $\ox_r$ &
\AxiomC{$s_1$}
\AxiomC{$s_2$}
\BinaryInfC{
    \begin{minipage}[t]{10ex}
    \begin{tabbing}
    fork $\alpha$ as \= $\mid \alpha_1 \rightarrow s_1$ \\
                     \> $\mid \alpha_2 \rightarrow s_2$
    \end{tabbing}
    \end{minipage}
}
\DisplayProof
    \bigskip\\
$\top_l$ and $\bot_r$ &
\AxiomC{$s$}
\UnaryInfC{close $\alpha$; $s$}
\DisplayProof
&
$\bot_l$ and $\top_r$ &
\AxiomC{}
\UnaryInfC{end $\alpha$}
\DisplayProof
    \bigskip\\
$\wdot_l$ and $\bdot_r$ &
\AxiomC{$s$}
\UnaryInfC{get $\alpha$ $x \Rightarrow s$}
\DisplayProof
&
$\bdot_l$ and $\wdot_r$ &
\AxiomC{$f$}
\AxiomC{$s$}
\BinaryInfC{put $\alpha$ $f$; $s$}
\DisplayProof
    \bigskip\\
$*$ &
\AxiomC{$s$}
\UnaryInfC{$s$}
\DisplayProof
&
$I$ &
\AxiomC{$s$}
\UnaryInfC{$s$}
\DisplayProof
    \bigskip\\
coprod &
\AxiomC{$s_1$}
\AxiomC{$\cdots$}
\AxiomC{$s_r$}
\TrinaryInfC{
    \begin{minipage}[t]{10ex}
    \begin{tabbing}
    case $z$ of \= $\mid$ \= $C_1 x_1 \rightarrow s_1$ \\
                        \>\> $\cdots$ \\
                \> $\mid$ $C_r x_r \rightarrow s_r$
    \end{tabbing}
    \end{minipage}
}
\DisplayProof
    \hspace{-8ex}
    \bigskip\\
subs &
\AxiomC{$f$}
\AxiomC{$s$}
\BinaryInfC{$s$ where $x=f$}
\DisplayProof
&
\end{tabular}
}\end{center}
\caption{Programming syntax for $\PMsg$}
\label{pmsg-prog-syn}
\end{figure}

\subsection{Cut elimination for $\PMsg$}

In this section the cut elimination rewrites are described. Recall that
unless explicitly mentioned no two subterms of a term may contain a variable
or channel name in common. Also recall that the notation $s[\alpha/\beta]$
means ``in $s$ substitute $\alpha$ for all occurrences of $\beta$''.

To recover the type of a term
\[
s :: \Phi \mid \Gamma \Vd \Delta
\]
denote $\Cont(s) = \Phi$, $\In(s) = \Gamma$, $\Out(s) = \Delta$, and
$\Chan(s) = \Gamma \cup \Delta$. These are respectively the \emph{context},
\emph{input channels}, \emph{output channels}, and \emph{channels} of the
term $s$. Similarly for $\Phi \vd f:A$ in $\Msg$ we have $\Cont(f)
= \Phi$ and $\Out(f) = A$.

Since there are two types of cuts, substitution and cut, the rewrites will
be split into two sections.

\subsubsection{Cut rewrites}

A cut in the message passing logic has the form
\begin{center}
\AxiomC{$\Pi$}
\UnaryInfC{$\Phi \mid \Gamma_1 \Vd \Delta_1,X$}
\AxiomC{$\Pi'$}
\UnaryInfC{$\Psi \mid X,\Gamma_2 \Vd \Delta_1$}
\BinaryInfC{$\Phi,\Psi \mid \Gamma_1,\Gamma_2 \Vd \Delta_1,\Delta_1$}
\DisplayProof
\end{center}
where both are derivations in $\PMsg$. The rewrites to eliminate cuts are
as follows.
\begin{align*}
\text{[id-sequent]} &\quad
    \alpha =_X \beta \cut{\beta}{\gamma} t \Ra t[\alpha/\gamma]
    \\[1.2ex]
\text{[sequent-id]} &\quad
    s \cut{\gamma}{\alpha} \alpha =_X \beta \Ra s[\beta/\gamma]
    \\[1.2ex]
\text{[$\ox_l/\oplus_r$-seq]} &\quad
    \alpha\<\alpha_1,\alpha_2\> \cdot s \cut{\beta}{\gamma} t \Ra
    \alpha\<\alpha_1,\alpha_2\> \cdot (s \cut{\beta}{\gamma} t)
    \\[1.2ex]
\text{[seq-$\ox_l/\oplus_r$]} &\quad
    s \cut{\beta}{\gamma} \alpha\<\alpha_1,\alpha_2\> \cdot t \Ra
    \alpha\<\alpha_1,\alpha_2\> \cdot (s \cut{\beta}{\gamma} t)
    \\[1.2ex]
\text{[$\oplus_l/\ox_r$-seq]} &\quad
    \alpha\scase{\alpha_1}{s_1}{\alpha_2}{s_2} \cut{\beta}{\gamma} t \Ra
        \begin{cases}
            \alpha\scase{\alpha_1}{s_1 \cut{\beta}{\gamma} t}
                        {\alpha_2}{s_2} & \text{if } \beta \in \Out(s_1)
            \bigskip\\
            \alpha\scase{\alpha_1}{s_1}
                        {\alpha_2}{s_2 \cut{\beta}{\gamma} t}
                            & \text{if } \beta \in \Out(s_2) \\
        \end{cases}
    \\[1.2ex]
\text{[seq-$\oplus_l/\ox_r$]} &\quad
    s \cut{\beta}{\gamma} \alpha\scase{\alpha_1}{t_1}{\alpha_2}{t_2} \Ra 
        \begin{cases}
            \alpha\scase{\alpha_1}{s \cut{\beta}{\gamma} t_1}
                        {\alpha_2}{t_2} & \text{if } \gamma \in \In(t_1)
            \bigskip\\
            \alpha\scase{\alpha_1}{t_1}
                        {\alpha_2}{s \cut{\beta}{\gamma} t_2}
                            & \text{if } \gamma \in \In(t_2) \\
        \end{cases}
    \\[1.2ex]
\text{[$\top_l/\bot_r$-seq]} &\quad
    \alpha\ainit \cdot s \cut{\beta}{\gamma} t \Ra
    \alpha\ainit \cdot (s \cut{\beta}{\gamma} t)
    \\[1.2ex]
\text{[seq-$\top_l/\bot_r$]} &\quad
    s \cut{\beta}{\gamma} \alpha\ainit \cdot t \Ra
    \alpha\ainit \cdot (s \cut{\beta}{\gamma} t)
    \\[1.2ex]
\text{[$\wdot_l/\bdot_r$-sequent]} &\quad
    \alpha\<x\> \cdot s \cut{\beta}{\gamma} t \Ra
    \alpha\<x\>\cdot (s \cut{\beta}{\gamma} t)
    \\[1.2ex]
\text{[sequent-$\wdot_l/\bdot_r$]} &\quad
    s \cut{\beta}{\gamma} \alpha\<x\> \cdot t \Ra
    \alpha\<x\>\cdot (s \cut{\beta}{\gamma} t)
    \\[1.2ex]
\text{[$\bdot_l/\wdot_r$-sequent]} &\quad
    \alpha[f] \cdot s \cut{\beta}{\gamma} t \Ra
    \alpha[f] \cdot (s \cut{\beta}{\gamma} t)
    \\[1.2ex]
\text{[sequent-$\bdot_l/\wdot_r$]} &\quad
    s \cut{\beta}{\gamma} \alpha[f] \cdot t \Ra
    \alpha[f] \cdot (s \cut{\beta}{\gamma} t)
    \\[1.2ex]
\text{[coprod-seq]} &\quad
    \case{x}{s_1}{y}{s_2}z \cut{\beta}{\gamma} t \Ra
    \case{x}{s_1 \cut{\beta}{\gamma} t}{y}{s_2 \cut{\beta}{\gamma} t}z
    \\[1.2ex]
\text{[seq-coprod]} &\quad
    s \cut{\beta}{\gamma} \case{x}{t_1}{y}{t_2}z \Ra
    \case{x}{s \cut{\beta}{\gamma} t_1}{y}{s \cut{\beta}{\gamma} t_2}z
    \\[1.2ex]
\text{[$\0$-sequent]} &\quad
    \init z \cut{\beta}{\gamma} t \Ra \init z
    \\[1.2ex]
\text{[sequent-$\0$]} &\quad
    s  \cut{\beta}{\gamma} \init z \Ra \init z
    \\[1.2ex]
\text{[$\ox_r$-$\ox_l$]} &\quad
    \alpha\scase{\alpha_1}{s_1}{\alpha_2}{s_2} \cut{\alpha}{\beta}
    \beta\<\beta_1,\beta_2\> \cdot t \Ra
    s_1 \cut{\alpha_1}{\beta_1} (s_2 \cut{\alpha_2}{\beta_2} t)
    \\[1.2ex]
\text{[$\oplus_r$-$\oplus_l$]} &\quad
    \alpha\<\alpha_1,\alpha_2\> \cdot s \cut{\alpha}{\beta}
    \beta\scase{\beta_1}{t_1}{\beta_2}{t_2} \Ra
    (s \cut{\alpha_1}{\beta_1} t_1) \cut{\alpha_2}{\beta_2} t_2
    \\[1.2ex]
\text{[$\top_r$-$\top_l$]} &\quad
    \alpha\sinit \cut{\alpha}{\beta} \beta\ainit \cdot t \Ra t
    \\[1.2ex]
\text{[$\bot_r$-$\bot_l$]} &\quad
    \alpha\ainit \cdot s \cut{\alpha}{\beta} \beta\sinit \Ra s
    \\[1.2ex]
\text{[$\wdot_r$-$\wdot_l$]} &\quad
    \alpha[f] \cdot s \cut{\alpha}{\beta} \beta\<x\> \cdot t \Ra
    s \cut{\alpha}{\beta} (x \mapsto t)f
    \\[1.2ex]
\text{[$\bdot_r$-$\bdot_l$]} &\quad
    \alpha\<x\> \cdot s \cut{\alpha}{\beta} \beta[f] \cdot t \Ra
    (x \mapsto s)f \cut{\alpha}{\beta} t
\end{align*}
The last two ``action'' cut elimination steps (the [$\wdot_r$-$\wdot_l$]
and [$\bdot_r$-$\bdot_l$] rewrites) are a bit novel so it will be nice to
see an explicit sequent cut elimination. This is the $\wdot_r$-$\wdot_l$
cut elimination step.
\[
\begin{array}{c}
  \AxiomC{$\pi$}
  \UnaryInfC{$\Phi_1 \vd f:A$}
  \AxiomC{$\Pi$} 
  \UnaryInfC{$\Phi_2 \mid \Gamma_1 \Vd X,\Delta_1$} 
  \BinaryInfC{$\Phi_1,\Phi_2 \mid \Gamma_1 \Vd A \circ X,\Delta_1$} 
  \AxiomC{$\Pi'$} 
  \UnaryInfC{$\Psi,A \mid \Gamma_2,X \Vd \Delta_2$}
  \UnaryInfC{$\Psi \mid \Gamma_2,A \circ X \Vd \Delta_2$}
  \BinaryInfC{$\Phi_1,\Phi_2,\Psi \mid \Gamma_1,\Gamma_2 \Vd \Delta_1,\Delta_2$}
  \DisplayProof
\smallskip\\ \xymatrix@!0@R=3.5ex{\ar@{=>}[d] \\ &} \\
  \AxiomC{$\Pi$} 
  \UnaryInfC{$\Phi_2 \mid \Gamma_1 \Vd X,\Delta_1$} 
  \AxiomC{$\pi$}
  \UnaryInfC{$\Phi_1 \vd f:A$}
  \AxiomC{$\Pi'$} 
  \UnaryInfC{$\Psi,A \mid \Gamma_2,X \Vd \Delta_2$}
  \BinaryInfC{$\Phi_1,\Psi \mid \Gamma_2,X \Vd \Delta_2$}
  \BinaryInfC{$\Phi_1,\Phi_2,\Psi \mid \Gamma_1,\Gamma_2 \Vd \Delta_1,\Delta_2$}
  \DisplayProof
\end{array}
\]
Notice that there are basically two kinds of interactions which are modeled
by the cut. The first is that one or the other of the terms is not active
on the channel of the cut. In this case the reduction is to pass the other
process into where it is active on
that channel. Notice in particular how a message gets evaluated when it is
passed. The other possibility is when the cut is on a channel for which, on
both sides, the type formation is the leading component of the term, i.e.,
both terms lead with activity on that channel (the last six rewrites above).
In this case we get a reduction which breaks down the type.

\subsubsection{Substitution rewrites}

A substitution in the message passing logic has the form
\begin{center}
\AxiomC{$\pi$}
\UnaryInfC{$\Psi \vd A$}
\AxiomC{$\Pi$}
\UnaryInfC{$\Phi,A \mid \Gamma \Vd \Delta$}
\BinaryInfC{$\Phi,\Psi \mid \Gamma \Vd \Delta$}
\DisplayProof
\end{center}
where $\pi$ is a derivation in $\Msg$ and $\Pi$ a derivation in $\PMsg$.
There are eleven substitution rewrites as follows.
\begin{align*}
\text{[subs-$\ox_l/\oplus_r$]} &\quad
    (w \mapsto \alpha\<\alpha_1,\alpha_2\> \cdot s)f \Ra
    \alpha\<\alpha_1,\alpha_2\> \cdot (w \mapsto s)f
    \\[1.2ex]
\text{[subs-$\oplus_l/\ox_r$]} &\quad
    \left(w \mapsto \alpha\scase{\alpha_1}{s}{\alpha_2}{t}\right)f \Ra
    \begin{cases}
        \alpha\scase{\alpha_1}{(w \mapsto s)f}{\alpha_2}{t} & \text{if }
            w \in \Cont(s) \medskip\\
        \alpha\scase{\alpha_1}{s}{\alpha_2}{(w \mapsto t)f} & \text{if }
            w \in \Cont(t) \\
    \end{cases}
    \\[1.2ex]
\text{[subs-$\top_l/\top_r$]} &\quad
    (w \mapsto \alpha\ainit \cdot s)f \Ra
    \alpha\ainit \cdot (w \mapsto s)f
    \\[1.2ex]
\text{[subs-$\wdot_l/\bdot_r$]} &\quad
    (w \mapsto \alpha\<x\> \cdot s)f \Ra
    \alpha\<x\> \cdot (w \mapsto s)f
    \\[1.2ex]
\text{[subs-$\bdot_l/\wdot_r$]} &\quad
    (w \mapsto \alpha[g] \cdot s)f \Ra
    \begin{cases}
        \alpha[g] \cdot (w \mapsto s)f & \text{if } w \in \Cont(s) \\
        \alpha[(w \mapsto g)f] \cdot s & \text{if } w \in \Cont(g) \\
    \end{cases}
    \\[1.2ex]
\text{[subs-coprod]} &\quad
    \left(w \mapsto \case{x}{s}{y}{t}z\right)f \Ra
    \case{x}{(w \mapsto s)f}{y}{(w \mapsto t)f}z
    \\[1.2ex]
\text{[subs-$\0$]} &\quad
    (w \mapsto \init z)f \Ra \init z
    \\[1.2ex]
\text{[$*_r$-$*$]} &\quad
    ((x,y) \mapsto s)(f,g) \Ra (x \mapsto (y \mapsto s)g)f
    \\[1.2ex]
\text{[$I_r$-$I$]} &\quad
    (\rinit \mapsto s)\rinit \Ra s
    \\[1.2ex]
\text{[inj$_l$-coprod]} &\quad
    \left(z \mapsto \case{x}{s}{y}{t}z\right)\sigma_1(f) \Ra (x \mapsto s)f
    \\[1.2ex]
\text{[inj$_r$-coprod]} &\quad
    \left(z \mapsto \case{x}{s}{y}{t}z\right)\sigma_2(f) \Ra (y \mapsto t)f
    \\[1.2ex]
\end{align*}
Presented as a derivation the [inj$_l$-coprod] rewrite is:
\[
\begin{array}{c}
\AxiomC{$\pi$}
\UnaryInfC{$\Phi \vd f:X$}
\UnaryInfC{$\Phi \vd \sigma_1(f):X+Y$}
\AxiomC{$\Pi$}
\UnaryInfC{$x:X,\Psi \mid \Gamma \Vd_s \Delta$}
\AxiomC{$\Pi'$}
\UnaryInfC{$y:Y,\Psi \mid \Gamma \Vd_t \Delta$}
\BinaryInfC{$z:X+Y,\Psi \mid \Gamma \Vd \Delta$}
\BinaryInfC{$\Phi,\Psi \mid \Gamma \Vd \Delta$}
\DisplayProof
\medskip\\ \Downarrow \\
\AxiomC{$\pi$}
\UnaryInfC{$\Phi \vd f:X$}
\AxiomC{$\Pi$}
\UnaryInfC{$x:X,\Psi \mid \Gamma \Vd_s \Delta$}
\BinaryInfC{$\Phi,\Psi \mid \Gamma \Vd \Delta$}
\DisplayProof~.
\end{array}
\]

\subsection{Equations in $\PMsg$} \label{sec-msg-equiv}

As in $\Msg$, in order that the cut elimination procedure is confluent
identities between the cut eliminated terms must also be added here.
Slightly unusual is the fact that this system has two distinct cut rules.
So that one is not given preference over the other, in addition to the usual
rewrites (which do not involve either of the cut rules), an identity is
required which allows the interchange of these two cut rules:
\[
\begin{array}{c}
\AxiomC{$\Pi$}
\UnaryInfC{$\Phi_2 \mid \Gamma_1 \Vd X,\Delta_1$} 
\AxiomC{$\pi$}
\UnaryInfC{$\Phi_1 \vd A$}
\AxiomC{$\Pi'$}
\UnaryInfC{$\Psi,A \mid \Gamma_2,X \Vd \Delta_2$}
\BinaryInfC{$\Phi_1,\Psi \mid \Gamma_2,X \Vd \Delta_2$}
\BinaryInfC{$\Phi_1,\Phi_2,\Psi \mid \Gamma_1,\Gamma_2 \Vd \Delta_1,\Delta_2$}
\DisplayProof
\\ \Downpc  \smallskip\\
\AxiomC{$\pi$}
\UnaryInfC{$\Phi_1 \vd A$}
\AxiomC{$\Pi$}
\UnaryInfC{$\Phi_2 \mid \Gamma_1 \Vd X,\Delta_1$} 
\AxiomC{$\Pi'$}
\UnaryInfC{$\Psi,A \mid \Gamma_2,X \Vd \Delta_2$}
\BinaryInfC{$\Phi_2,\Psi,A \mid \Gamma_1,\Gamma_2 \Vd \Delta_1,\Delta_2$}
\BinaryInfC{$\Phi_1,\Phi_2,\Psi \mid \Gamma_1,\Gamma_2 \Vd \Delta_1,\Delta_2$}
\DisplayProof
\end{array}
\]
This is called the [cut-subs] interchange below.

In presenting the identities we shall always assume that both sides make
sense. For example, the equation
\[
\alpha\<x\> \cdot \beta[f] \cdot s \pc \beta[f] \cdot \alpha\<x\> \cdot s
\]
only makes sense if $x$ does not occur in $\Cont(f)$.

Many of the identities below simply express that actions on one channel
should be independent -- as far as is possible given the bindings -- from
actions on another. Thus, although there may seem to be a lot of equations
their genus is quite simple.

\noindent $\bullet$ The two cut rewrites.
\medskip

\begin{tabular}{rl}
\qquad [subs-cut] & $(x \mapsto s)f \cut{\beta}{\gamma} t \pc
    (x \mapsto s \cut{\beta}{\gamma} t)f$
    \\[1.2ex]
\qquad [cut-subs] & $s \cut{\beta}{\gamma} (x \mapsto t)f \pc
    (x \mapsto s \cut{\beta}{\gamma} t)f$
\end{tabular}
\medskip

\noindent $\bullet$ Rewrites involving $\ox_l$ or $\oplus_r$.
\medskip

\begin{tabular}{rl}
[$\ox_l/\oplus_r$-$\ox_l/\oplus_r$] &
$\alpha\<\alpha_1,\alpha_2\> \cdot \beta\<\beta_1,\beta_2\> \cdot s \pc
    \beta\<\beta_1,\beta_2\> \cdot \alpha\<\alpha_1,\alpha_2\> \cdot s$
\end{tabular}
\medskip

\begin{tabular}{l}
[$\ox_l/\oplus_r$-$\oplus_l/\ox_r$] \\
\quad $\alpha\<\alpha_1,\alpha_2\> \cdot \scase{\beta_1}{s}{\beta_2}{t} \pc
    \begin{cases}
    \beta\scase{\beta_1}{\alpha\<\alpha_1,\alpha_2\> \cdot s}{\beta_2}{t}
    & \text{if } \alpha_1,\alpha_2 \in \Chan(s) \medskip\\
    \beta\scase{\beta_1}{s}{\beta_2}{\alpha\<\alpha_1,\alpha_2\> \cdot t}
    & \text{if } \alpha_1,\alpha_2 \in \Chan(t)
    \end{cases}$
\end{tabular}
\medskip

\begin{tabular}{rl}
[$\ox_l/\oplus_r$-$\top_l/\bot_r$] &
    $\alpha\<\alpha_1,\alpha_2\> \cdot \beta\ainit \cdot s \pc
    \beta\ainit \cdot \alpha\<\alpha_1,\alpha_2\> \cdot s$
    \\[2ex]
[$\ox_l/\oplus_r$-$\wdot_l/\bdot_r$] &
    $\alpha\<\alpha_1,\alpha_2\> \cdot \beta\<x\> \cdot f \pc
    \beta\<x\> \cdot \alpha\<\alpha_1,\alpha_2\> \cdot f$
    \\[2ex]
[$\ox_l/\oplus_r$-$\bdot_l/\wdot_r$] &
    $\alpha\<\alpha_1,\alpha_2\> \cdot \beta[f] \cdot s \pc
    \beta[f] \cdot \alpha\<\alpha_1,\alpha_2\> \cdot s$
    \\[2ex]
[$\ox_l/\oplus_r$-coprod] &
    $\alpha\<\alpha_1,\alpha_2\> \cdot \case{x}{s}{y}{t}z \pc
    \case{x}{\alpha\<\alpha_1,\alpha_2\> \cdot s}
         {y}{\alpha\<\alpha_1,\alpha_2\> \cdot t}z$
    \\[5ex]
[$\ox_l/\oplus_r$-$\0$] &
    $\alpha\<\alpha_1,\alpha_2\> \cdot \init z \pc \init z$
\end{tabular}
\medskip

\noindent $\bullet$ Rewrites involving $\oplus_l$ or $\ox_r$. We only
describe an interchange occurring on the left-hand (the $\alpha_1$) branch.
Similar rewrites are needed for the right-hand (the $\alpha_2$) branch, but
are omitted here as they are easy to infer.

\begin{align*}
\text{[$\oplus_l/\ox_r$-$\oplus_l/\ox_r$]} \quad &
    \alpha\scase{\alpha_1}{\beta\scase{\beta_1}{s_1}{\beta_2}{s_2}}
                 {\alpha_2}{t} \pc
     \beta\scase{\beta_1}{\alpha\scase{\alpha_1}{s_1}{\alpha_2}{t}}
                {\beta_2}{s_2}
    \\[1ex]
\text{[$\oplus_l/\ox_r$-$\top_l/\bot_r$]} \quad &
    \alpha\scase{\alpha_1}{\beta\ainit \cdot s}{\alpha_2}{t} \pc
        \beta\ainit \cdot \alpha\scase{\alpha_1}{s}{\alpha_2}{t}
    \\[1ex]
\text{[$\oplus_l/\ox_r$-$\wdot_l/\bdot_r$]} \quad &
    \alpha\scase{\alpha_1}{\beta\<x\> \cdot s}{\alpha_2}{t} \pc
        \beta\<x\> \cdot \alpha\scase{\alpha_1}{s}{\alpha_2}{t}
    \\[1ex]
\text{[$\oplus_l/\ox_r$-$\bdot_l/\wdot_r$]} \quad &
    \alpha\scase{\alpha_1}{\beta[f] \cdot s}{\alpha_2}{t} \pc
     \beta[f] \cdot \alpha\scase{\alpha_1}{s}{\alpha_2}{t}
    \\[1ex]
\text{[$\oplus_l/\ox_r$-coprod]} \quad &
    \alpha\scase{\alpha_1}{\case{x}{s_1}{y}{s_2}z}{\alpha_2}{t} \pc
     \case{x}{\alpha\scase{\alpha_1}{s_1}{\alpha_2}{t}}
          {y}{\alpha\scase{\alpha_1}{s_2}{\alpha_2}{t}}z
    \\[1ex]
\text{[$\oplus_l/\ox_r$-$\0$]} \quad &
    \alpha\scase{\alpha_1}{\init z}{\alpha_2}{t} \pc \init z
\end{align*}
\noindent $\bullet$ Rewrites involving $\top_l$ or $\bot_r$.
\begin{align*}
\text{[$\top_l/\bot_r$-$\top_l/\bot_r$]} \quad &
    \alpha\ainit \cdot \beta\ainit \cdot s \pc
        \beta\ainit \cdot \alpha\ainit \cdot s
    \\[1ex]
\text{[$\top_l/\bot_r$-$\wdot_l/\bdot_r$]} \quad &
    \alpha\ainit \cdot \beta\<x\> \cdot s \pc
        \beta\<x\> \cdot \alpha\ainit \cdot s
    \\[1ex]
\text{[$\top_l/\bot_r$-$\bdot_l/\wdot_r$]} \quad &
    \alpha\ainit \cdot \beta[f] \cdot s \pc
        \beta[f] \cdot \alpha\ainit \cdot s
    \\[1ex]
\text{[$\top_l/\bot_r$-coprod]} \quad &
    \alpha\ainit \cdot \beta\case{x}{s}{y}{t}z \pc
        \case{x}{\beta\ainit \cdot s}{y}{\beta\ainit \cdot t}z \qquad\qquad
    \\[1ex]
\text{[$\top_l/\bot_r$-$\0$]} \quad &
    \alpha\ainit \cdot \beta\init \cdot s \pc \beta\init
\end{align*}
\noindent $\bullet$ Rewrites involving $\wdot_l$ or $\bdot_r$.
\begin{align*}
\text{[$\wdot_l/\bdot_r$-$\wdot_l/\bdot_r$]} \quad &
    \alpha\<x\> \cdot \beta\<y\> \cdot s \pc
    \beta\<y\> \cdot \alpha\<x\> \cdot s
    \\[1ex]
\text{[$\wdot_l/\bdot_r$-$\bdot_l/\wdot_r$]} \quad &
    \alpha\<x\> \cdot \beta[f] \cdot s \pc \beta[f] \cdot \alpha\<x\> \cdot s
    \\[1ex]
\text{[$\wdot_l/\bdot_r$-coprod]} \quad &
    \alpha\<x\> \cdot \case{y_1}{s}{y_2}{t}z \pc
    \case{y_1}{\alpha\<x\> \cdot s}{y_2}{\alpha\<x\> \cdot t}z \qquad
    \\[1ex]
\text{[$\wdot_l/\bdot_r$-$\0$]} \quad &
    \alpha\<x\> \cdot \init z \pc \init z
\end{align*}
\noindent $\bullet$ Rewrites involving $\bdot_l$ or $\wdot_r$. In this case
we must also investigate the message terms which are passed by the $\bdot_l$
and $\wdot_r$ rules as they may also lead to interchanges with the coproduct or
$\0$ rules. Thus, the two distinct interchanges for each of the coproduct and
$\0$ rule.
\begin{align*}
\text{[$\bdot_l/\wdot_r$-$\bdot_l/\wdot_r$]} \quad &
    \alpha[f] \cdot \beta[g] \cdot s \pc \beta[g] \cdot \alpha[f] \cdot s
    \\[1ex]
\text{[$\bdot_l/\wdot_r$-coprod]} \quad &
    \alpha[f] \cdot \case{x}{s}{y}{t}z \pc
    \case{x}{\alpha[f] \cdot s}{y}{\alpha[f] \cdot t}z
    \\[1ex]
\text{[$\bdot_l/\wdot_r$-coprod]} \quad &
    \alpha\left[\case{x}{f}{y}{g}z\right] \cdot s \pc
    \case{x}{\alpha[f] \cdot s}{y}{\alpha[g] \cdot s}z \qquad
    \\[1ex]
\text{[$\bdot_l/\wdot_r$-$\0$]} \quad &
    \alpha[f] \cdot \init z \pc \init z
    \\[1ex]
\text{[$\bdot_l/\wdot_r$-$\0$]} \quad &
    \alpha[\init z] \cdot s \pc \init z
\end{align*}
\noindent $\bullet$ Rewrites involving the coproduct.

\begin{tabular}{l}
[coprod-coprod] \\[1.2ex] \quad
   $\case{x_1}{\case{y_1}{s_1}{y_2}{t_1}w}{x_2}{\case{y_1}{s_2}{y_2}{t_2}w}z \pc
    \case{y_1}{\case{x_1}{s_1}{x_2}{s_2}z}{y_2}{\case{x_1}{t_1}{x_2}{t_2}z}w$
\end{tabular}

\begin{tabular}{rl}
[coprod-$\0$] &
   $\case{x}{\init z}{y}{\init z}w \pc \init z$
\end{tabular}

We shall demonstrate how these identities are used in Section~\ref{semantics}
where we relate this term logic to the categorical semantics.

\section{Linear actegories}

Given a (symmetric) monoidal category $\A$ we introduce the notion of a
linear $\A$-actegory, which is a linearly distributive category $\X$
equipped with two functors
\[
\wdot: \xymatrix{\A \x \X \ar[r] & \X} \qquad \text{and} \qquad
\bdot:\xymatrix{\A^\op \x \X \ar[r] & \X},
\]
the ``actions'' of $\A$ on $\X$. These must satisfy a number of coherence 
conditions which are described below.

Our aim, in the following section, is to show that these form the categorical 
semantics for the proof theory of $\PMsg$.

\subsection{Linearly distributive categories}

A \emph{linearly distributive category} is a category $\X$ equipped with a
``tensor'' $\ox:\X \x \X \ra \X$ with unit $\top$ and coherent natural
isomorphisms
\[
a_\ox:(X \ox Y) \ox Z \ra X \ox (Y \ox Z), \quad
l_\ox:\top \ox X \ra X, \quad
r_\ox:X \ox \top \ra X,
\]
and a ``par'' $\parr:\X \x \X \ra \X$ with unit $\bot$ and coherent natural
isomorphism
\[
a_\parr:X \parr (Y \parr Z) \ra (X \parr Y) \parr Z, \quad
l_\parr:X \ra \bot \parr X, \quad
r_\parr:X \ra X \parr \bot
\]
(the odd choice of direction is used to maximise the symmetry below),
and two linear distributions
\[
d^\ox_\parr : X \ox (Y \parr Z) \ra (X \ox Y) \parr Z \quad \text{and} \quad
d^\parr_\ox : (Y \parr Z) \ox X \ra Y \parr (Z \ox X)
\]
relating the two structures. This data must satisfy several coherence
conditions (see~\cite{CS}).

If both the tensor and the par are symmetric (with $c_\ox$ and $c_\parr$)
and several other coherence conditions are satisfied (again see~\cite{CS})
then it is called a \emph{symmetric} linear distributive category. In this
case there are two induced ``permuting'' linear distributions
\[
d^{\ox'}_\parr \!:\! X \ox (Y \parr Z) \ra Y \parr (X \ox Z)
    \quad \text{and} \quad
d^{\parr'}_\ox \!:\! (Y \parr Z) \ox X \ra (Y \ox X) \parr Z.
\]

\begin{example} \quad
\begin{enumerate}
\item Any distributive lattice is a linearly distributive category with
the objects being the elements and the maps being comparisons: $\ox$ is
the meet and $\parr$ is the join.

\item Any monoidal category gives rise to a ``compact'' (i.e., $\ox = \parr$)
linearly distributive category. When both $\ox$ and $\parr$ are interpreted
by the same tensor the linear distribution becomes associativity.

\item Any $*$-autonomous category is a linearly distributive with $A \parr
B := (B^* \ox A^*)$.

\item The category of sets with $\ox = \x$ and $\oplus$ given as follows
(due to J\"urgen Koslowski):
\[
A \parr B  = \begin{cases}
    B & A = \emptyset \\
    A & B = \emptyset \\
    1 & \text{otherwise.}
\end{cases}
\]
This is an example of a non-compact, non-posetal, non-$*$-autonomous
linearly distributive category.
\end{enumerate}
\end{example}

\subsection{Linear actegories}

Let $\A = (\A,*,I,a_*,l_*,r_*,c_*)$ be a symmetric monoidal category.
A \emph{(symmetric) linear $\A$-actegory} consists of the following data.

\begin{itemize}
\item A symmetric linearly distributive category $\X$ (as above),

\item Functors
\[
\wdot:\A \x \X \ra \X \quad \text{and} \quad \bdot:\A^\op \x \X \ra \X,
\]
such that $\wdot$ is the left parameterised left adjoint of $\bdot$, i.e.,
for all $A \in \A$, $A \wdot - \dashv A \bdot -$. The unit and counit of
this adjunction (natural in $A \in \A$ and $X \in \X$) are denoted
respectively by
\[
n_{A,X}:X \ra A \bdot (A \wdot X) \quad \text{and} \quad 
e_{A,X}:A \wdot (A \bdot X) \ra X.
\]

\item For all $A,B \in \A$ and $X,Y \in \X$ natural isomorphisms in $\X$
\begin{align*}
u_\wdot &: I \wdot X \ra X, \\
u_\bdot &: X \ra I \bdot X, \\
a^*_\wdot &:(A * B) \wdot X \ra A \wdot (B \wdot X), \\
a^*_\bdot &: A \bdot (B \bdot X) \ra (A * B) \bdot X, \\
a^\wdot_\ox &: A \wdot (X \ox Y) \ra (A\wdot X) \ox Y, \\
a^\bdot_\parr &: (A \bdot X) \parr Y \ra A \bdot (X \parr Y).
\end{align*}

\item For all $A,B \in \A$ and $X,Y \in \X$ natural morphisms in $\X$
\begin{align*}
d^\wdot_\parr &: A \wdot (X \parr Y) \ra (A \wdot X) \parr Y, \\
d^\bdot_\ox &: (A \bdot X) \ox Y \ra A \bdot (X \ox Y), \\
d^\wdot_\bdot &: A \wdot (B \bdot X) \ra B \bdot (A \wdot X).
\end{align*}
\end{itemize}
The symmetries of $*$, $\ox$, and $\parr$ induce the following permuting
morphisms (or isomorphisms):
\begin{align*}
a^{*'}_\wdot &: (A * B) \wdot X \ra B \wdot (A \wdot X), \\
a^{*'}_\bdot &: B \bdot (A \bdot X) \ra (A*B) \bdot X, \\
a^\wdot_{\ox'} &: A \wdot (X \ox Y) \ra X \ox (A \wdot Y), \\
a^\bdot_{\parr'} &: X \parr (A \bdot Y) \ra A \bdot (X \parr Y), \\
d^\wdot_{\parr'} &: A \wdot (X \parr Y) \ra X \parr (A \wdot Y), \\
d^\bdot_{\ox'} &: X \ox (A \bdot Y) \ra A \bdot (X \ox Y).
\end{align*}
This data must satisfy several coherence conditions which we shall discuss
shortly. Firstly we try to give some intuition behind the notation that has
been chosen. The ``$a$'' maps are (invertible) associativity isomorphisms
and the ``$d$'' maps are (non-invertible) linear distributions. The direction
of the maps have been chosen (when there is a choice) to maximise the amount
of symmetry and so that the $\wdot$ is pushed in a bracket and the $\bdot$
is pulled out of a bracket. This choice of notation may allow us to leave
off the subscripts and let the types disambiguate the maps (which is not
however done here).

The symmetries of this data are as follows:

[$\op'$] Reverse the arrows and swap $\ox$ and $\parr$, $\top$ and $\bot$,
and $\wdot$ and $\bdot$. This gives the following assignment of generating
maps
\[
\begin{array}{cccc}
a_\ox \iff a_\parr
    ~&~ l_\ox \iff l_\parr
    ~&~ r_\ox \iff r_\parr
    ~&~ c_\ox \iff c^{-1}_\parr \\
u_\wdot \iff u_\bdot
    ~&~ a^*_\wdot \iff a^*_\bdot
    ~&~ a^\wdot_\ox \iff a^\bdot_\parr \\
d^\ox_\parr \iff d^\parr_\ox
    ~&~ d^\wdot_\parr \iff d^\bdot_\ox 
    ~&~ d^\wdot_\bdot \iff d^\wdot_\bdot \\
n \iff e
\end{array}
\]
with the remainder of the morphisms unchanged. 

[$*'$] Reverse the $*$ (i.e., $A *' B = B * A$); this assigns
\[
\begin{array}{ccccc}
a_* \iff a^{-1}_*
    ~&~ l_* \iff r_*
    ~&~ c_* \iff c^{-1}_*
    ~&~ a^*_\wdot \iff a^{*'}_\wdot
    ~&~ a^*_\bdot \iff a^{*'}_\bdot 
\end{array}
\]
with the remainder unchanged.

[$\ox'$] Reverse the $\ox$; this assigns
\[
\begin{array}{cccc}
a_\ox \iff a^{-1}_\ox
    ~&~ l_\ox \iff r_\ox
    ~&~ c_\ox \iff c^{-1}_\ox \\
d^\ox_\parr \iff d^{\parr'}_\ox
    ~&~ d^\parr_\ox \iff d^{\ox'}_\parr
    ~&~ a^\wdot_\ox \iff a^\wdot_{\ox'}
    ~&~ d^\bdot_\ox \iff d^\bdot_{\ox'}
\end{array}
\]
with the remainder unchanged.

[$\parr'$] Reverse the $\parr$; this assigns
\[
\begin{array}{cccc}
a_\parr \iff a^{-1}_\parr
    ~&~ l_\parr \iff r_\parr
    ~&~ c_\parr \iff c^{-1}_\parr \\
d^\ox_\parr \iff d^{\ox'}_\parr
    ~&~ d^\parr_\ox \iff d^{\parr'}_\ox
    ~&~ a^\bdot_\parr \iff a^\bdot_{\parr'}
    ~&~ d^\wdot_\parr \iff d^\wdot_{\parr'}
\end{array}
\]
with the remainder unchanged.

[$*'$] Reverse the $*$; this assigns
\[
\begin{array}{cccc}
a_* \iff a_*^{-1}
   ~&~ l_* \iff l_*^{-1} 
   ~&~ r_* \iff r_*^{-1} 
   ~&~ c_* \iff c_*^{-1}
\end{array}
\]
with the remainder unchanged.

[$\odot'$] There are four remaining symmetries obtained by reversing any
combination of two or more of $*$, $\ox$, and $\parr$. The assignments are
evident.

The notion of a linear $\A$-actegory is preserved by these symmetries. It is 
important to notice that the first symmetry is the most significant as it 
indicates a fundamental relationship between {\em different\/} functorial 
operations. 

The coherence conditions for a linear $\A$-actegory are now described.

[Symmetries.] The two diagrams below linking the symmetries and the
associativities must commute.
\[
\xygraph{{A \wdot (X \ox Y)} (
    :[r(2.5)] {(A \wdot X) \ox Y} ^-{a^\wdot_\ox}
    :[d(1.2)] {Y \ox (A \wdot X)}="e" ^-{c_\ox}
    )
    :[d(1.2)] {A \wdot (Y \ox X)} _-{A \wdot c_\ox}
    :"e" ^-{a^\wdot_{\ox'}}}
\quad
\xygraph{{A \wdot (X \parr Y)} (
    :[r(2.5)] {(A \wdot X) \parr Y} ^-{d^\wdot_\parr}
    :[d(1.2)] {Y \parr (A \wdot X)}="e" ^-{c_\parr}
    )
    :[d(1.2)] {A \wdot (Y \parr X)} _-{A \wdot c_\parr}
    :"e" ^-{a^\wdot_{\parr'}}}
\]
These diagrams and the results of applying the symmetries to them yield
the following equations:
\begin{align*}
a^\wdot_\ox ; c_\ox &= (A \wdot c_\ox) ; a^\wdot_{\ox'} &
    c_\ox ; a^\bdot_\parr &= a^\bdot_{\parr'} ; (A \bdot c_\parr)
    \tag{1} \\[1ex]
d^\wdot_\parr ; c_\parr &= (A \wdot c_\parr) ; a^\wdot_{\parr'} &
    c_\parr ; d^\bdot_\ox &= a^\bdot_{\ox'} ; (A \bdot c_\ox)
    \tag{2}
\end{align*}
[Unit and associativity.] The following diagrams linking the unit and
associatively morphisms must commute. 
\[
\begin{array}{cc}
    \xygraph{{(A * I) \wdot X}
        (:[r(2.3)] {A \wdot (I \wdot X)} ^-{a^*_\wdot}
         :[d] {A \wdot X}="e" ^-{A \wdot u_\wdot}
        )
        :"e" _-{r_* \wdot X}}
\qquad & \qquad
    \xygraph{{(I * A) \wdot X}
        (:[r(2.3)] {I \wdot (A \wdot X)} ^-{a^*_\wdot}
         :[d] {A \wdot X}="e" ^-{u_\wdot}
        )
        :"e" _-{l_* \wdot X}}
\bigskip\\
    \xygraph{{A \wdot (X \ox \top)}
        (:[r(2.3)] {(A \wdot X) \ox \top} ^-{a^\wdot_\ox}
         :[d] {A \wdot X}="e" ^-{r_\ox}
        )
        :"e" _-{A \wdot r_\ox}}
\qquad & \qquad
    \xygraph{{A \wdot X}
        (:[r(2.3)] {A \wdot (X \parr \bot)} ^-{A \wdot r_\parr}
         :[d] {(A \wdot X) \parr \bot}="e" ^-{d^\wdot_\parr}
        )
        :"e" _-{r_\parr}}
\end{array}
\]
These result in the following equations:
\begin{align*}
a^*_\wdot ; (A \wdot u_\wdot) &= r_* \wdot X
    & (A \bdot u_\bdot); a^*_\bdot &= r^{-1}_* \bdot X
    \tag{3} \\[1ex]
a^*_\wdot ; u_\wdot &= l_* \wdot X
    & u_\bdot ; a^*_\bdot &= l^{-1}_* \bdot X
    \tag{4} \\[1ex]
a^\wdot_\ox ; r_\ox &= A \wdot r_\ox
    & r_\parr ; a^\bdot_\parr &= A \bdot r_\parr
    \tag{5} \\[1ex]
a^\wdot_{\ox'} ; l_\ox &= A \wdot l_\ox
    & l_\parr ; a^\bdot_{\parr'} &= A \bdot l_\parr
\\[1ex]
(A \wdot r_\parr) ; d^\wdot_\parr &= r_\parr
    & d^\bdot_\ox ; (A \bdot r_\ox) &= r_\ox
    \tag{6} \\[1ex]
(A \wdot l_\parr) ; d^\wdot_{\parr'} &= l_\parr
    & d^\bdot_{\ox'} ; (A \bdot l_\ox) &= l_\ox
\end{align*}
[Unit and distributivity.] The following diagrams linking the unit and
distributivity morphisms must commute. 
\[
\begin{array}{cc}
\xygraph{{I \wdot (X \parr Y)}
   (:[r(2.3)] {(I \wdot X) \parr Y} ^-{d^\wdot_\parr}
    :[d] {X \parr Y}="e" ^-{u_\wdot \parr Y}
   )
    :"e" _-{u_\wdot}}
\qquad & \qquad
\xygraph{{I \wdot (X \ox Y)}
   (:[r(2.3)] {(I \wdot X) \ox Y} ^-{a^\wdot_\ox}
    :[d] {X \ox Y}="e" ^-{u_\wdot \ox Y}
   )
    :"e" _-{u_\wdot}}
\bigskip
\\
\xygraph{{I \wdot (A \bdot X)}
   (:[r(2.3)] {A \bdot (I \wdot X)} ^-{d^\wdot_\bdot}
    :[d] {A \bdot X}="e" ^-{A \bdot u_\wdot}
   )
    :"e" _-{u_\wdot}}
\qquad & \qquad
\xygraph{{A \wdot X}
   (:[r(2.3)] {A \wdot (I \bdot X)} ^-{A \wdot u_\bdot}
    :[d] {I \bdot (A \wdot X)}="e" ^-{d^\wdot_\bdot}
   )
    :"e" _-{u_\bdot}}
\end{array}
\]
These result in the following equations:
\begin{align*}
d^\wdot_\parr ; (u_\wdot \parr Y) &= u_\wdot
    & (u_\bdot \ox Y) ; d^\bdot_\ox &= u_\bdot
\tag{7} \\[1ex]
d^\wdot_{\parr'} ; (Y \parr u_\wdot) &= u_\wdot
    & (Y \ox u_\bdot) ; d^\bdot_{\ox'} &= u_\bdot
\\[1ex]
a^\wdot_\ox ; (u_\wdot \ox Y) &= u_\wdot
    & (u_\bdot \parr Y) ; a^\bdot_\parr &= u_\bdot
\tag{8} \\[1ex]
a^\wdot_{\ox'} ; (Y \ox u_\wdot) &= u_\wdot
    & (Y \parr u_\bdot) ; a^\bdot_{\parr'} &= u_\bdot
\\[1ex]
d^\wdot_\bdot ; (A \bdot u_\wdot) &= u_\wdot
    & (A \wdot u_\bdot) ; d^\wdot_\bdot &= u_\bdot
\\[1ex]
d^\wdot_\bdot ; (A \bdot u_\wdot) &= u_\wdot
    & (A \wdot u_\bdot) ; d^\wdot_\bdot &= u_\bdot
\tag{9} \\[1ex]
A \wdot u_\bdot ; d^\wdot_\bdot &= u_\bdot
    & d^\wdot_\bdot ; A \bdot u_\wdot &= u_\wdot
\tag{10}
\end{align*}
[Associativity.] The following diagrams are ones linking the various
associativity morphisms. There are four pentagon shaped diagrams which
must commute.
\[
\xygraph{
{((A * B) * C) \wdot X} (
    :[u(1.7)r(1.4)] {(A * (B * C)) \wdot X} ^-{a_* \wdot X}
    :[dr] {A \wdot ((B * C) \wdot X)} ^-{a^*_\wdot}
    :[d(1.4)] {A \wdot (B \wdot (C \wdot X))}="e" ^-{A \wdot a^*_\wdot}
    )
    :[d(1.7)r(1.4)] {(A * B) \wdot (C \wdot X)} _-{a^*_\wdot}
    :"e" _-{a^*_\wdot}}
\
\xygraph{ {(A * B) \wdot (C \bdot X)} (
    :[u(1.7)r(1.4)] {A \wdot (B \wdot (C \bdot X))} ^-{a^*_\wdot}
    :[dr] {A \wdot (C \bdot (B \wdot X))} ^-{A \wdot d^\wdot_\bdot}
    :[d(1.4)] {C \bdot (A \wdot (B \wdot X))}="e" ^-{d^\wdot_\bdot}
    )
    :[d(1.7)r(1.4)] {C \bdot ((A * B) \wdot X)} _-{d^\wdot_\bdot}
    :"e" _-{C \bdot a^*_\wdot}}
\]
\[
\xygraph{
{(A * B) \wdot (X \ox Y)} (
    :[u(1.7)r(1.4)] {A \wdot (B \wdot (X \ox Y))} ^-{a^*_\wdot}
    :[dr] {A \wdot ((B \wdot X) \ox Y)} ^-{A \wdot a^\wdot_\ox}
    :[d(1.4)] {(A \wdot (B \wdot X)) \ox Y}="e" ^-{a^\wdot_\ox}
    )
    :[d(1.7)r(1.4)] {((A * B) \wdot X) \ox Y} _-{a^\wdot_\ox}
    :"e" _-{a^*_\wdot \ox Y}}
\
\xygraph{
{A \wdot ((X \ox Y) \ox Z)} (
    :[u(1.7)r(1.4)] {(A \wdot (X \ox Y)) \ox Z} ^-{{a^\wdot_\ox}}
    :[dr] {((A \wdot X) \ox Y) \ox Z} ^-{a^\wdot_\ox \ox Z}
    :[d(1.4)] {(A \wdot X) \ox (Y \ox Z)}="e" ^-{a_\ox}
    )
    :[d(1.7)r(1.4)] {A \wdot (X \ox (Y \ox Z))} _-{A \wdot a_\ox}
    :"e" _-{a^\wdot_\ox}}
\]
These result in the following equations:
\[
\begin{array}{lcc}
(11) & (a_* \wdot X) ; a^*_\wdot ; (A \wdot a^*_\wdot) = a^*_\wdot ; a^*_\wdot
& (A \bdot a^*_\bdot) ; a^*_\bdot ; (a^{-1}_* \bdot X) = a^*_\bdot ; a^*_\bdot 
\\[1ex]
& (a^{-1}_* \wdot X) ; a^{*'}_\wdot ; (A \wdot a^{*'}_\wdot)
    = a^{*'}_\wdot ; a^{*'}_\wdot
    & (A \bdot a^{*'}_\bdot) ; a^{*'}_\bdot ; (a_* \bdot X)
    = a^{*'}_\bdot ;a^{*'}_\bdot 
\\[1ex]
(12) & a^*_\wdot ; (A \wdot d^\wdot_\bdot) ; d^\wdot_\bdot
    = d^\wdot_\bdot ; (C \bdot a^*_\wdot)
    & d^\wdot_\bdot ; (A \bdot d^\wdot_\bdot) ; a^*_\wdot
    = C \wdot a^*_\bdot ; d^\bdot_\wdot
\\[1ex]
& a^{*'}_\wdot ; (A \wdot d^\wdot_\bdot) ; d^\wdot_\bdot
    = d^\wdot_\bdot ; (C \bdot a^{*'}_\wdot)
    & d^\wdot_\bdot ; (A \bdot d^\wdot_\bdot) ; a^{*'}_\wdot
    = C \wdot a^{*'}_\bdot ; d^\bdot_\wdot
\\[1ex]
(13) & a^*_\wdot ; (A \wdot a^\wdot_\ox) ; (a^\wdot_\ox)
    = (a^\wdot_\ox) ; (a^*_\wdot \ox Y)
    & a^\bdot_\parr ; (A \bdot a^\bdot_\parr) ; a^*_\bdot
    = (a^*_\bdot \parr Y) ; a^\bdot_\parr
\\[1ex]
& a^{*'}_\wdot ; (A \wdot a^\wdot_\ox) ; (a^\wdot_\ox)
    = (a^\wdot_\ox) ; (a^{*'}_\wdot \ox Y)
    & a^\bdot_\parr ; (A \bdot a^\bdot_\parr) ; a^{*'}_\bdot
    = (a^{*'}_\bdot \parr Y) ; a^\bdot_\parr
\\[1ex]
& a^*_\wdot ; (A \wdot a^\wdot_{\ox'}) ; (a^\wdot_{\ox'})
    = (a^\wdot_{\ox'}) ; (Y \ox a^*_\wdot)
    & a^\bdot_{\parr'} ; (A \bdot a^\bdot_{\parr'}) ; a^*_\bdot
    = (Y \parr a^*_\bdot) ; a^\bdot_{\parr'}
\\[1ex]
& a^{*'}_\wdot ; (A \wdot a^\wdot_{\ox'}) ; (a^\wdot_{\ox'})
    = (a^\wdot_{\ox'}) ; (Y \ox a^{*'}_\wdot)
    ~&~ a^\bdot_{\parr'} ; (A \bdot a^\bdot_{\parr'}) ; a^{*'}_\bdot
    = (Y \parr a^{*'}_\bdot) ; a^\bdot_{\parr'}
\\[1ex]
(14) & a^\wdot_\ox ; (a^\wdot_\ox \ox Z) ; a_\ox = (A \wdot a_\ox) ; a^\wdot_\ox
    & a_\parr ; (a^\bdot_\parr \ox Z) ; a^\bdot_\parr =
    a^\bdot_\parr ; (A \bdot a_\parr)
\\[1ex]
& a^\wdot_{\ox'} ; (Z \ox a^\wdot_{\ox'}) ; a^{-1}_\ox =
    (A \wdot a^{-1}_\ox) ; a^\wdot_{\ox'}
    ~&~ a^{-1}_\parr ; (Z \ox a^\bdot_{\parr'}) ; a^\bdot_{\parr'} =
    a^\bdot_{\parr'} ; (A \bdot a^{-1}_\parr)
\end{array}
\]
There is additionally a hexagon shaped diagram which must commute
\[
\xygraph{{(A*B) \wdot (X \ox Y)} (
    :[r(1.4)d] {B\! \wdot\! (A\! \wdot\! (X\! \ox\! Y))} ^-{a^{*'}_\wdot}
    :[d]  {B\! \wdot\! ((A\! \wdot\! X)\! \ox\! Y)} ^-{B \wdot a^\wdot_\ox}
    :[l(1.4)d] {(A \wdot X) \ox (B \wdot Y)}="e" ^-{a^\wdot_{\ox'}}
    )
    :[l(1.4)d] {A\! \wdot\! (B\! \wdot\! (X\! \ox\! Y))} _-{a^*_\wdot}
    :[d]  {A\! \wdot\! (X\! \ox\! (B\! \wdot\! Y))} _-{A \wdot a^\wdot_{\ox'}}
    :"e" _-{a^\wdot_\ox}}
\]
which results in the equations:
\begin{align*}
a^*_\wdot ; (A\! \wdot\! a^\wdot_{\ox'}) ; a^\wdot_\ox &=
    a^{*'}_\wdot ; (B\! \wdot\! a^\wdot_\ox) ; a^\wdot_{\ox'}
    & a^\bdot_\parr ; (A\! \bdot\! a^\bdot_{\parr'}) ; a^*_\bdot &=
    a^{*'}_\bdot ; (B\! \bdot\! a^\bdot_\parr) ; a^\bdot_{\parr'}
\tag{15} \\[1ex]
a^*_\wdot ; (A\! \wdot\! a^\wdot_\ox) ; a^\wdot_{\ox'}
    &= a^{*'}_\wdot ; (B\! \wdot\! a^\wdot_{\ox'}) ; a^\wdot_\ox & 
a^\bdot_{\parr'} ; (A\! \bdot\! a^\bdot_\parr) ; a^*_\bdot
    &= a^{*'}_\bdot ; (B\! \bdot\! a^\bdot_{\parr'}) ; a^\bdot_\parr
\end{align*}
[Distributivity and associativity.] Diagrams made up of the
distributivity morphisms. There are five pentagon shaped diagrams which must
commute.
\[
\xygraph{ {(A * B) \wdot (X \parr Y)} (
    :[u(1.7)r(1.4)] {A \wdot (B \wdot (X \parr Y))} ^-{a^*_\wdot}
    :[dr] {A \wdot ((B \wdot X) \parr Y)} ^-{A \wdot d^\wdot_\parr}
    :[d(1.4)] {(A \wdot (B \wdot X)) \parr Y}="e" ^-{d^\wdot_\parr}
    )
    :[d(1.7)r(1.4)] {((A * B) \wdot X) \parr Y} _-{d^\wdot_\parr}
    :"e" _-{a^*_\wdot \parr Y}}
\
\xygraph{ {A \wdot ((B \bdot X) \ox Y)} (
    :[u(1.7)r(1.4)] {A \wdot (B \bdot (X \ox Y))} ^-{A \wdot d^\bdot_\ox}
    :[dr] {B \bdot (A \wdot (X \ox Y))} ^-{d^\wdot_\bdot}
    :[d(1.4)]  {B \bdot (X \ox (A \wdot Y))}="e" ^-{B \bdot a^\wdot_{\ox'}}
    )
    :[d(1.7)r(1.4)] {(B \bdot X) \ox (A \wdot Y)} _-{a^\wdot_{\ox'}}
    :"e"  _-{d^\bdot_\ox}}
\]
\[
\xygraph{ {A \wdot (X \ox (Y \parr Z))} (
    :[u(1.7)r(1.4)] {A \wdot ((X \ox Y) \parr Z)} ^-{A \wdot d^\ox_\parr}
    :[dr] {(A \wdot (X \ox Y)) \parr Z} ^-{d^\wdot_\parr}
    :[d(1.4)] {((A \wdot X) \ox Y) \parr Z}="e" ^-{a^\wdot_\ox \parr Z} 
    )
    :[d(1.7)r(1.4)] {(A \wdot X) \ox (Y \parr Z)} _-{a^\wdot_\ox}
    :"e" _-{d^\ox_\parr}}
\
\xygraph{ {A \wdot ((Y \parr Z) \ox X)} (
    :[u(1.7)r(1.4)] {(A \wdot (Y \parr Z)) \ox X} ^-{a^\wdot_\ox}
    :[dr] {((A \wdot Y) \parr Z) \ox X} ^-{d^\wdot_\parr \ox X}
    :[d(1.4)] {(A \wdot Y) \parr (Z \ox X)}="e" ^-{d^\parr_\ox} 
    )
    :[d(1.7)r(1.4)] {A \wdot (Y \parr (Z \ox X))} _-{A \wdot d^\parr_\ox}
    :"e" _-{d^\wdot_\parr}}
\]
\[
\xygraph{ {A \wdot (X \parr (Y \parr Z))} (
    :[u(1.7)r(1.4)] {A \wdot ((X \parr Y) \parr Z)} ^-{A \wdot a_\parr}
    :[dr] {(A \wdot (X \parr Y)) \parr Z} ^-{d^\wdot_\parr}
    :[d(1.4)] {((A \wdot X) \parr Y) \parr Z}="e" ^-{d^\wdot_\parr \parr Z} 
    )
    :[d(1.7)r(1.4)] {(A \wdot X) \parr (Y \parr Z)} _-{d^\wdot_\parr}
    :"e" _-{a_\parr}}
\]
These result in the following equations:
\[
\begin{array}{lcc}
(16) & a^*_\wdot ; (A \wdot d^\wdot_\parr) ; d^\wdot_\parr
    = d^\wdot_\parr ; (a^*_\wdot \parr Y)
    & d^\bdot_\ox ; (A \bdot d^\bdot_\ox) ; a^*_\bdot
    = (a^*_\bdot \ox Y) ; d^\bdot_\ox
\\[1ex]
& a^{*'}_\wdot ; (A \wdot d^\wdot_\parr) ; d^\wdot_\parr
    = d^\wdot_\parr ; (a^{*'}_\wdot \parr Y)
    & d^\bdot_\ox ; (A \bdot d^\bdot_\ox) ; a^{*'}_\bdot
    = (a^{*'}_\bdot \ox Y) ; d^\bdot_\ox
\\[1ex]
& a^*_\wdot ; (A \wdot d^\wdot_{\parr'}) ; d^\wdot_{\parr'}
    = d^\wdot_{\parr'} ; (Y \parr a^*_\wdot)
    & d^\bdot_{\ox'} ; (A \bdot d^\bdot_{\ox'}) ; a^*_\bdot
    = (Y \ox a^*_\bdot) ; d^\bdot_{\ox'}
\\[1ex]
& a^{*'}_\wdot ; (A \wdot d^\wdot_{\parr'}) ; d^\wdot_{\parr'}
    = d^\wdot_{\parr'} ; (Y \parr a^{*'}_\wdot)
    & d^\bdot_{\ox'} ; (A \bdot d^\bdot_{\ox'}) ; a^{*'}_\bdot
    = (Y \ox a^{*'}_\bdot) ; d^\bdot_{\ox'}
\\[1ex]
(17) & (A \wdot d^\bdot_\ox) ; d^\wdot_\bdot ; (B \bdot a^\wdot_{\ox'})
    = a^\wdot_{\ox'} ; d^\bdot_\ox
    & (B \wdot a^\bdot_{\parr'}) ; d^\wdot_\bdot ; (A \bdot d^\wdot_\parr)
    = d^\wdot_\parr ; a^\bdot_{\parr'}
\\[1ex]
& (A \wdot d^\bdot_{\ox'}) ; d^\wdot_\bdot ; (B \bdot a^\wdot_\ox)
    = a^\wdot_\ox ; d^\bdot_{\ox'}
    & (B \wdot a^\bdot_\parr) ; d^\wdot_\bdot ; (A \bdot d^\wdot_{\parr'})
    = d^\wdot_{\parr'} ; a^\bdot_\parr
\\[1ex]
(18) & (A \wdot d^\ox_\parr) ; d^\wdot_\parr ; (a^\wdot_\ox \parr Z)
    = a^\wdot_\ox ; d^\ox_\parr
    & (a^\bdot_\parr \ox Z) ; d^\bdot_\ox ; (A \bdot d^\parr_\ox)
    = d^\parr_\ox ; a^\bdot_\parr 
\\[1ex]
& (A \wdot d^{\parr'}_\ox) ; d^\wdot_\parr ; (a^\wdot_{\ox'} \parr Z)
    = a^\wdot_{\ox'} ; d^{\parr'}_\ox
    & (a^\bdot_{\parr'} \ox Z) ; d^\bdot_\ox ; (A \bdot d^{\parr'}_\ox)
    = d^{\parr'}_\ox ; a^\bdot_{\parr'} 
\\[1ex]
& (A \wdot d^{\ox'}_\parr) ; d^\wdot_{\parr'} ; (Z \parr a^\wdot_\ox)
    = a^\wdot_\ox ; d^{\ox'}_\parr
    & (Z \ox a^\bdot_\parr) ; d^\bdot_{\ox'} ; (A \bdot d^{\ox'}_\parr)
    = d^{\ox'}_\parr ; a^\bdot_\parr 
\\[1ex]
& (A \wdot d^\parr_\ox) ; d^\wdot_{\parr'} ; (Z \parr a^\wdot_{\ox'})
    = a^\wdot_{\ox'} ; d^\parr_\ox
    & (Z \ox a^\bdot_{\parr'}) ; d^\bdot_{\ox'} ; (A \bdot d^\ox_\parr)
    = d^\ox_\parr ; a^\bdot_{\parr'}
\\[1ex]
(19) & a^\wdot_\ox ; (d^\wdot_\parr \ox X) ; d^\parr_\ox
    = (A \wdot d^\parr_\ox) ; d^\wdot_\parr
    & d^\ox_\parr ; (d^\bdot_\ox \parr X) ; a^\bdot_\parr 
    = d^\bdot_\ox ; (A \bdot d^\ox_\parr)
\\[1ex]
& a^\wdot_{\ox'} ; (X \ox d^\wdot_\parr) ; d^{\ox'}_\parr
    = (A \wdot d^{\ox'}_\parr) ; d^\wdot_\parr
    & d^{\ox'}_\parr ; (X \parr d^\bdot_\ox) ; a^\bdot_{\parr'} 
    = d^\bdot_\ox ; (A \bdot d^{\ox'}_\parr)
\\[1ex]
& a^\wdot_\ox ; (d^\wdot_{\parr'} \ox X) ; d^{\parr'}_\ox
    = (A \wdot d^{\parr'}_\ox) ; d^\wdot_{\parr'}
    & d^{\parr'}_\ox ; (d^\bdot_{\ox'} \parr X) ; a^\bdot_\parr 
    = d^\bdot_{\ox'} ; (A \bdot d^{\parr'}_\ox)
\\[1ex]
& a^\wdot_{\ox'} ; (X \ox d^\wdot_{\parr'}) ; d^\ox_\parr
    = (A \wdot d^\ox_\parr) ; d^\wdot_{\parr'}
    & d^\parr_\ox ; (X \parr d^\bdot_{\ox'}) ; a^\bdot_{\parr'}
    = d^\bdot_{\ox'} ; (A \bdot d^\parr_\ox)
\\[1ex]
(20) & (A \wdot a_\parr) ; d^\wdot_\parr ; (d^\wdot_\parr \parr Z)
    = d^\wdot_\parr ; a_\parr
    & (d^\bdot_\ox \ox Z) ; d^\bdot_\ox ; (A \bdot a_\ox)
    = a_\ox ; d^\bdot_\ox 
\\[1ex]
& (A \wdot a^{-1}_\parr) ; d^\wdot_{\parr'} ; (Z \parr d^\wdot_{\parr'})
    = d^\wdot_{\parr'} ; a^{-1}_\parr
    ~&~ (Z \ox d^\bdot_{\ox'}) ; d^\bdot_{\ox'} ; (A \bdot a^{-1}_\ox)
    = a^{-1}_\ox ; d^\bdot_{\ox'} 
\end{array}
\]
There are additionally three hexagon shaped diagram which must commute.
The first is
\[
\xygraph{{(A*B) \wdot (X \parr Y)} (
    :[r(1.4)d] {B\! \wdot\! (A\! \wdot\! (X\! \parr\! Y))} ^-{a^{*'}_\wdot}
    :[d(1.2)]  {B\! \wdot\! ((A\! \wdot\! X)\! \parr\! Y)} ^-{B \wdot d^\wdot_\parr}
    :[l(1.4)d] {(A \wdot X) \parr (B \wdot Y)}="e" ^-{d^\wdot_{\parr'}}
    )
    :[l(1.4)d] {A\! \wdot\! (B\! \wdot\! (X\! \parr\! Y))} _-{a^*_\wdot}
    :[d(1.2)]  {A\! \wdot\! (X\! \parr\! (B\! \wdot\! Y))} _-{A \wdot d^\wdot_{\parr'}}
    :"e" _-{d^\wdot_\parr}}
\]
which results in the equations:
\begin{align*}
a^*_\wdot ; (A\! \wdot\! d^\wdot_{\parr'}) ; d^\wdot_\parr &=
    a^{*'}_\wdot ; (B\! \wdot\! d^\wdot_\parr) ; d^\wdot_{\parr'}
    & d^\bdot_\ox ; (A\! \bdot\! d^\bdot_{\ox'}) ; a^*_\bdot &=
    d^\bdot_{\ox'} ; (B\! \bdot\! d^\bdot_\ox) ; a^{*'}_\bdot
\tag{21} \\[1ex]
a^*_\wdot ; (A\! \wdot\! d^\wdot_\parr) ; d^\wdot_{\parr'} &=
    a^{*'}_\wdot ; (B\! \wdot\! d^\wdot_{\parr'}) ; d^\wdot_\parr
    & d^\bdot_{\ox'} ; (A\! \bdot\! d^\bdot_\ox) ; a^*_\bdot &=
    d^\bdot_\ox ; (B\! \bdot\! d^\bdot_{\ox'}) ; a^{*'}_\bdot
\end{align*}
The remaining two are
\[
\xygraph{{A \wdot ((B \bdot X) \ox Y)} (
    :[r(1.4)d] {(A\! \wdot\! (B\! \bdot\! X))\! \ox\! Y} ^-{a^\wdot_\ox}
    :[d(1.2)]  {(B\! \bdot\! (A\! \wdot\! X))\! \ox\! Y} ^-{d^\wdot_\bdot \ox Y}
    :[l(1.4)d] {B \bdot ((A \wdot X) \ox Y)}="e" ^-{d^\bdot_\ox}
    )
    :[l(1.4)d] {A\! \wdot\! (B\! \bdot\! (X\! \ox\! Y))} _-{A \wdot d^\bdot_\ox}
    :[d(1.2)]  {B\! \bdot\! (A\! \wdot\! (X\! \ox\! Y))} _-{d^\wdot_\bdot}
    :"e" _-{B \bdot a^\wdot_\ox}}
\quad
\xygraph{{A \wdot ((B \bdot X) \parr Y)} (
    :[r(1.4)d] {(A\! \wdot\! (B\! \bdot\! X))\! \parr\! Y} ^-{d^\wdot_\parr}
    :[d(1.2)]  {(B\! \bdot\! (A\! \wdot\! X))\! \parr\! Y} ^-{d^\wdot_\bdot \parr Y}
    :[l(1.4)d] {B \bdot ((A \wdot X) \parr Y)}="e" ^-{a^\bdot_\parr}
    )
    :[l(1.4)d] {A\! \wdot\! (B\! \bdot\! (X\! \parr\! Y))} _-{A \wdot
a^\bdot_\parr}
    :[d(1.2)]  {B\! \bdot\! (A\! \wdot\! (X\! \parr\! Y))} _-{d^\wdot_\bdot}
    :"e" _-{B \bdot d^\wdot_\parr}}
\]
and these result in the equations:
\begin{align*}
a^\wdot_\ox ; (d^\wdot_\bdot \ox Y) ; d^\bdot_\ox
    &= (A \wdot d^\bdot_\ox) ; d^\wdot_\bdot ; (B \bdot a^\wdot_\ox)
\tag{22} \\[1ex]
d^\wdot_\parr ; (d^\wdot_\bdot \parr Y) ; a^\bdot_\parr
    &= (B \wdot a^\bdot_\parr) ; d^\wdot_\bdot ; (A \bdot d^\wdot_\parr)
\\[1ex]
a^\wdot_{\ox'} ; (Y \ox d^\wdot_\bdot) ; d^\bdot_{\ox'}
    &= (A \wdot d^\bdot_{\ox'}) ; d^\wdot_\bdot ; (B \bdot a^\wdot_{\ox'})
\\[1ex]
d^\wdot_{\parr'} ; (Y \parr d^\wdot_\bdot) ; a^\bdot_{\parr'}
    &= (B \wdot a^\bdot_{\parr'}) ; d^\wdot_\bdot ; (A \bdot d^\wdot_{\parr'})
\\[1ex]
d^\wdot_\parr ; (d^\wdot_\bdot \parr Y) ; a^\bdot_\parr
    &= (A \wdot a^\bdot_\parr) ; d^\wdot_\bdot ; (B \bdot d^\wdot_\parr)
\tag{23} \\[1ex]
a^\wdot_\ox ; (d^\wdot_\bdot \ox Y) ; d^\bdot_\ox
    &= (B \wdot d^\bdot_\ox) ; d^\wdot_\bdot ; (A \bdot a^\wdot_\ox) \\[1ex]
d^\wdot_{\parr'} ; (Y \parr d^\wdot_\bdot) ; a^\bdot_{\parr'}
    &= (A \wdot a^\bdot_{\parr'}); d^\wdot_\bdot; (B \bdot d^\wdot_{\parr'})
\\[1ex]
a^\wdot_{\ox'} ; (Y \ox d^\wdot_\bdot) ; d^\bdot_{\ox'}
    &= (B \wdot d^\bdot_{\ox'}) ; d^\wdot_\bdot ; (A \bdot a^\wdot_{\ox'})
\end{align*}
[$n$ and $e$.] Finally, there are diagrams linking the unit and counit
of the $A \wdot - \dashv A \bdot -$ adjunction with $I$ and the
associativity and distributivity morphisms.
\[
    \vcenter{\xymatrix{X \ar[r]^-{n_{I,X}} \ar[dr]_-{u_\bdot} 
              & I \bdot (I \circ X) \ar[d]^-{I \bdot u_\circ} \\
              & I \bdot X}}
\quad
    \vcenter{\xymatrix{X \ar[r]^-{n_{A*B,X}} \ar[d]_{n_{A,X}}
     & (A * B) \bdot ((A * B) \wdot X) \ar[dd]^{(A * B) \bdot a^{*'}_\wdot} \\
             A \bdot (A \circ X) \ar[d]_{A \bdot n_{B,A \wdot X}} \\
             A \bdot (B \bdot (B \circ (A \wdot X)) \ar[r]^-{a_\bdot^*} 
             & (A * B) \bdot (B \circ (A \wdot X))}}
\]
\[
\xygraph{{X \ox Y}
   (:[ur(1.2)] {A \bdot (A \wdot (X \ox Y))} ^-{n_{A,X \ox Y}}
    :[dr(1.2)] {A \bdot ((A \wdot X) \ox Y)}="e" ^-{A \bdot a^\wdot_\ox}
   )
    :[dr(1.2)] {(A \bdot (A \wdot X)) \ox Y} _-{n_{A,X} \ox Y}
    :"e" _-{d^\bdot_\ox}}
\qquad
\xygraph{{X \parr Y}
   (:[ur(1.2)] {A \bdot (A \wdot (X \parr Y))} ^-{n_{A,X \parr Y}}
    :[dr(1.2)] {A \bdot ((A \wdot X) \parr Y)}="e" ^-{A \bdot d^\wdot_\parr}
   )
    :[dr(1.2)] {(A \bdot (A \wdot X)) \parr Y} _-{n_{A,X} \parr Y}
    :"e" _-{a^\ox_\bdot}}
\]
With the symmetries this gives the equations:
\[
\begin{array}{lcc}
(24) & n_{I,X} ; I \bdot u_\wdot = u_\bdot & e_{I,X} ; u_\bdot = u_\wdot
    \\[1ex]
(25) & n_{A*B,X} ; (A*B) \bdot a^{*'}_\wdot =
            n_{A,X} ; A \bdot n_{B,A \wdot X} ; a_\bdot^*
     ~&~ a^{*'}_\bdot ; e_{A*B,X} =
            a^*_\wdot ; A \wdot e_{B,A\bdot X} ; e_{A,X}
    \\[1ex]
(26) & n_{X \ox Y} ; (A \bdot a^\wdot_\ox) = (n_X \ox Y) ; d^\bdot_\ox
    & (A \wdot a^\bdot_\parr) ; e_{X \parr Y} = d^\wdot_\parr ; (e_X \parr Y)
    \\[1ex]
     & n_{Y \ox X} ; (A \bdot a^\wdot_{\ox'}) = (Y \ox n_X) ; d^\bdot_{\ox'}
     & (A \wdot a^\bdot_{\parr'}) ; e_{Y \parr X} = d^\wdot_{\parr'} ;
       (Y \parr e_X) \\[1ex]
(27) & n_{X \parr Y} ; (A \bdot d^\wdot_\parr) = (n_X \parr Y) ; a^\bdot_\parr 
    & (A \wdot d^\bdot_\ox) ; e_{X \ox Y} = a^\wdot_\ox ; (e_X \ox Y)
    \\[1ex]
     & n_{Y \parr X};(A \bdot d^\wdot_{\parr'}) = (Y \parr n_X);a^\bdot_{\parr'}
     & (A \wdot d^\bdot_{\ox'}) ; e_{Y \ox X} = a^\wdot_{\ox'} ; (Y \ox e_X)
\end{array}
\]
In the next section we will explore these diagrams in more detail.

\begin{example} \quad
\begin{enumerate}
\item Given any monoidal closed category, regarding it as a compact linearly
distributive category it acts on itself via
\[
A \circ B = A \ox B \qquad \text{and} \qquad A \bdot B = A \multimap B.
\]
This fails in general to give a linear actegory as the isomorphism
$a^\bdot_\parr: (A \bdot X) \parr Y \ra A \bdot (X \parr Y)$ is absent.
However if one restricts the action to the compact objects
(i.e., those objects for which the natural map $(A \multimap I) \ox B \ra A
\multimap B$ is an isomorphism) then this becomes a linear actegory with
actions as above.

\item Compact closed categories are, of course, a source of examples of
the above. They are (compact) $*$-autonomous categories and so a bridge to the
next example. They are important as not only are they the foundation for
what Girard calls the ``geometry of interaction'', but also for a family of
compact closed categories (essentially span categories) which were studied
by Abramsky~\cite{A2,A3} under the name of ``interaction categories''.
These used explicit ideas from process calculus to give a categorical
semantics for processes. Acting on themselves these give examples of linear
actegories. Abramsky, Gay, and Nagarajan also considered adding
``specifications'' to these categories~\cite{AGN} which made them (mix)
$*$-autonomous categories and, thus, less degenerate models of linear
actegories. 

\item Given any linearly distributive category $\X$, the objects which have
linear adjoints (complements) form a $*$-autonomous subcategory $\Map(\X)$.
There is an obvious action of $\Map(\X)$ on $\X$ defined by $A \circ X =
A \ox X$ and $A \bdot X = A^* \parr X$ (much as in the first example). In
particular, a $*$-autonomous category acting in the obvious manner on itself
is a linear actegory.

\item None of the above examples illustrate well the separation of messages 
from the message passing. However, using Benton's approach~\cite{B}
to models of linear logic (with exponentials), which links the
intuitionistic terms to the linear terms by a monoidal adjunction, gives
an important model of a linear actegory where such separation is displayed.

This model was used by Barber et al.~\cite{BGHP} to provide a semantics 
for Milner's action calculus. As the action calculus was developed, in part, 
to provide a semantic framework for systems such as the $\pi$-calculus this 
suggests that there is a very close connection between the work of Barber
et al.~\cite{BGHP} and what is being proposed here. This is indeed the case, 
however, there are also some important differences. A model of their logic
is an example of a linear actegory only when their monoidal closed category
is actually a $*$-autonomous category. In this regard linear actegories
demand \emph{more} structure. On the other hand, not every linear actegory
arises in this manner. One very obvious reason is that in their setting the
messages and the processes are still very closely linked by the adjoint
and it is thus possible to turn process code into message code. While often
in practice this may be a desirable feature, it may also be something that
one does not want to allow. In our models we do not assume such a connection
exists.

If $\X$ is a $*$-autonomous category with an exponential comonad $!:\X \ra
\X$ then, as described by Benton~\cite{B}, the comonad induces a monoidal 
adjunction $V \vdash W: \X \ra \X_!$, where essentially $W =\; !$
and $\X_!$ is a cartesian closed category. The action $\circ:\X_! \x \X
\ra \X$ is
\[
X \circ Y = W(X) \ox Y
\]
and as $W(X \x Y) \cong W(X) \ox W(Y)$ this automatically gives an action.
The action $\bullet:(\X_!)^\op \x \X \ra \X$ is given by
\[
X \bullet Y = W(X)^* \oplus Y.
\]
That this is a linear actegory is now a straightforward, if lengthy, exercise.
\end{enumerate}
\end{example}

\section{Categorical semantics} \label{semantics}

We shall say that a linear $\A$-actegory is $\A$-additive in case the
monoidal category $\A$ is a distributive monoidal category (i.e., it has
coproducts over which the tensor distributes) and the covariant action
preserves these coproducts while the contravariant action turns them into
products.

Our aim is to show that the proof theory of the message passing logic, as 
represented by the terms, forms a linear additive actegory in the above
sense. To achieve this we shall show in this section how one may collect
the proof theory for message passing into a linear additive actegory 
(completeness).  In the next section we show that given any interpretation 
of the axioms into such a linear actegory one can extend the interpretation to 
the whole message logic (soundness).  

To preserve the sanity of reader and writer alike we shall present a recipe 
for these processes exemplifying only some of the details.  These matters 
have already been well-explored for the message fragment and the linearly 
distributive fragment. Accordingly, when it comes to the details we shall 
focus on the actions. We begin by proving completeness. This involves 
showing: 
\begin{enumerate}[(a)]
\item The proofs of the message passing logic sequents with empty context and 
one input and output type 
\[
\emptyset \mid X \Vd Y 
\]
form a linearly distributive category.  

\item The proofs of the message logic with one input (and necessarily one
output)
\[
x:A \vd s:B
\]
form a distributive monoidal category.

\item The two required actions of a linear actegory are present together
with the coherence maps and that they satisfy all the required coherences
of the previous section.
\end{enumerate}

Parts (a) and (b) have been established in~\cite{CKS} and~\cite{Her}
respectively, so we shall concentrate our efforts on (c). To begin this
proof we need to establish the functorial nature of the actions. This, in
turn, will lead into the naturality of the coherence maps and establishing
the coherences.

\subsection{The actions are functors}

Here is the definition of the action $\circ:\A \x \X \ra \X$ built from an
arbitrary monoidal map $f:A \ra B$ and a process $s:X \ra Y$.
\begin{center}
\AxiomC{$x:A \vd f:B \quad$}
\AxiomC{$\quad s :: \alpha:X \Vd \beta:Y$}
\BinaryInf$\beta[f] \cdot s ~\fCenter~ x:A \mid \alpha:X \Vd \beta:B \circ Y$
\UnaryInf$\alpha\<x\> \cdot \beta[f] \cdot s ~\fCenter~ \emptyset \mid
          \alpha:A \circ X \Vd \beta:B \circ Y$
\DisplayProof
\end{center}
In order to show that this is a functor we must show that it preserves
composition. As composition amounts to a cut, this amounts to showing that
cuts inside the functor can be equivalently expressed as a cut outside. Here
is the calculation in reverse.
\begin{align*}
\alpha\<x\> \cdot \beta[f] \cdot s \cut{\beta}{\gamma} \gamma\<y\> \cdot
        \delta[g] \cdot t 
 &~\Ra~ \alpha\<x\> \cdot (\beta[f] \cdot s \cut{\beta}{\gamma} \gamma\<y\>
       \cdot \delta[g] \cdot t) \\
 &~\Ra~ \alpha\<x\> \cdot (y\mapsto s \cut{\beta}{\gamma} \delta[g] \cdot t)f \\
 &~\Ra~ \alpha\<x\>\cdot (y\mapsto\delta[g]\cdot (s \cut{\beta}{\gamma} t))f \\
 &~\Ra~ \alpha\<x\>\cdot \delta[(y \mapsto g)f]\cdot (s \cut{\beta}{\gamma} t)
\end{align*}

Thus $\circ$ preserves composition. As we shall see shortly, by definition,
it preserves identities.

The symmetry which is embodied in the term logic means that the action
$\bdot:\A^\op \x \X \ra \X$ has an identical term though the arrangement
of the types is different.
\begin{center}
\AxiomC{$x:A \vd f:B \quad$}
\AxiomC{$\quad s:: \alpha: X \Vd \beta: Y$}
\BinaryInf$\alpha[f] \cdot s ~\fCenter~ x:A \mid \alpha:B \bdot X \Vd \beta:Y$
\UnaryInf$\beta\<x\> \cdot \alpha[f] \cdot s ~\fCenter~ \emptyset \mid
          \alpha: B \bdot X \Vd \beta:A \bdot Y$
\DisplayProof
\end{center}
This means the preservation of composition (and identities) for $\bdot$ is
essentially the same proof.

\subsection{Identities}

We shall denote the identity map on a type by
\[
\alpha =_X \beta :: \alpha:X \ra \beta:X.
\]
However, its definition as a term depends on the type $X$. If the type is
primitive then this identity is built-in and defined to behave in the
correct manner.  However, if $X$ is not primitive we must provide an
inductive definition.
\begin{align*}
\alpha =_X \beta &\quad \text{is} \quad
    \alpha =_X \beta \qquad \text{for $X$ a primitive type}
    \\[1ex]
\alpha =_{X \ox Y} \beta &\quad \text{is} \quad
    \alpha\<\alpha_1,\alpha_2\> \cdot \beta
    \scase{\beta_1}{\alpha_1 =_X \beta_1}{\beta_2}{\alpha_2 =_Y \beta_2}
    \\[1ex]
\alpha =_{X \parr Y} \beta &\quad \text{is} \quad
    \beta\<\beta_1,\beta_2\> \cdot \alpha
    \scase{\alpha_1}{\alpha_1 =_X \beta_1}{\alpha_2}{\alpha_2 =_Y \beta_2}
    \\[1ex]
\alpha =_\top \beta &\quad \text{is} \quad
    \alpha\ainit \cdot \beta\sinit
    \\[1ex]
\alpha =_\bot \beta &\quad \text{is} \quad
    \beta\ainit \cdot \alpha\sinit
    \\[1ex]
\alpha =_{A \circ X} \beta  &\quad \text{is} \quad
             \alpha\<x\>\cdot\beta[\iota(x)]\cdot \alpha =_X \beta
    \\[1ex]
\alpha =_{A \bdot X} \beta &\quad \text{is} \quad
            \beta\<x\>\cdot\alpha[\iota(x)]\cdot \alpha =_X \beta
\end{align*}
Here $x:A \vd \iota(x):A$ is the identity map in the message logic for the
type $A$ (where $x$ here stands for a pattern in general). Notice also that
this definition confirms that the functors (above) preserve identities.
Clearly we also need an inductive definition of the identities in the
message logic:
\begin{align*}
x:A \vd \iota(x):A & ~~\text{is}~~
    x:A \vd x:A \qquad \text{for $A$ a primitive type}
    \\[1ex]
(x,y)\!:\!A*B \vd \iota(x,y)\!:\!A*B & ~~\text{is}~~ (x,y)\!:\!A*B \vd
    (\iota(x),\iota(y))\!:\!A*B
    \\[1ex]
z\!:\!A+B \vd \iota(z)\!:\!A+B & ~~\text{is}~~
    z\!:\!A+B \vd \case{x}{\sigma_0(\iota(x))}{y}{\sigma_1(\iota(y))}z:A+B
\end{align*}
It remains to do an inductive proof that these terms do act as identity maps.
We shall focus on the step for the covariant action to give a feel of how
this inductive proof (which is straightforward) plays out. Consider
\[
\alpha\<x\>\cdot\beta[\iota(x)]\cdot \alpha =_X \beta ~~;_\beta~ s.
\]
We do an induction on the structure of the term $s$. There are two basic
cases: either the leading structure of $s$ interacts along the channel
$\beta$ or it does not. If it does not then this allows us to push the
identity term inside the leading structure and to invoke the inductive
hypothesis. This leaves the case in which there is interaction on
$\beta$ and this means $s$ must be of the form $\beta\<y\>\cdot s'$
and we have
\begin{align*}
\alpha\<x\> \cdot \beta[\iota(x)] \cdot \alpha =_X \beta ~;_\beta
        \beta\<y\>\cdot s'
    &~\Ra~ \alpha\<x\> \cdot (y \mapsto \alpha =_X \beta ~;_\beta s')
        \iota(x) \\
    &~\Ra~ \alpha\<y\> \cdot \alpha =_X \beta ~;_\beta s' \\
    &~\Ra~ \alpha\<y\> \cdot s'
\end{align*}
where we invoke the induction hypothesis to obtain the last step.

\subsection{Associativity and interchange}

The associativity of composition, as represented by a cut, must also be
proven. The proof is again an inductive argument concerning the cut
elimination process.  Here we highlight some of the inductive steps of this
argument which involve the action $\circ:\A \x \X \ra \X$. We consider the
cases determined by terms which have their leading action on this type. This
means the term is either of the form $\alpha\<x\> \cdot s$ or $\alpha[f]
\cdot s$. We consider the first case in more detail. In a sequence of three
cuts such a term can occur of course in three positions. For whichever
position it occurs in it is shown that the inductive hypothesis may be used
on a combination of smaller terms to show that the original term is
associative.

Consider the case when it is in the first position and the cut does not
occur on $\alpha$. The following diagram of cut elimination rewrites, where
the bottom equality uses the inductive hypothesis proves associativity in
this case.
\[
\xymatrix{
(\alpha\<x\> \cdot s ~;_\beta t) ~;_\gamma u \ar@{=>}[d] &
    \alpha\<x\> \cdot s ~;_\beta (t ~;_\gamma u) \ar@{=>}[dd] \\
(\alpha\<x\> \cdot (s ~;_\beta t)) ~;_\gamma u \ar@{=>}[d] & \\
\alpha\<x\> \cdot ((s ~;_\beta t) ~;_\gamma u)~ \ar@{|=|}[r] &
    ~\alpha\<x\> \cdot (s ~;_\beta (t ~;_\gamma u))}
\]
A very similar argument holds for the middle position provided $\alpha$ is
not the interacting channel determined by the leftmost cut.
\[
\xymatrix{
(s ~;_\beta \alpha\<x\> \cdot t) ~;_\gamma u \ar@{=>}[d]
    & s ~;_\beta (\alpha\<x\> \cdot t ~;_\gamma u) \ar@{=>}[dd] \\
(\alpha\<x\> \cdot (s ~;_\beta t)) ~;_\gamma u \ar@{=>}[d] \\
\alpha\<x\> \cdot ((s ~;_\beta t) ~;_\gamma u) \ar@{|=|}[r] &
    \alpha\<x\> \cdot (s ~;_\beta (t ~;_\gamma u))}
\]
Now if $\alpha$ is the channel of the leftmost cut then either the leftmost
terms leading action is on that channel or not. If it is not, we can move
the cut inside the action. Now provided the second cut (on $\gamma$) is
also not on that channel this can be moved inside the action and we can
then invoke the inductive hypothesis. Fortunately, due to the way cut
works, it is impossible for the outer terms to share a channel, thus the
second cut is guaranteed to be independent of this action.

This leaves the case when the leftmost terms leading action is on the
channel. The following diagram then proves this case.
\[
\xymatrix{
(\beta[f] \cdot s ~;_\beta \beta\<x\> \cdot t) ~;_\gamma u \ar@{=>}[d]
  & \beta[f] \cdot s ~;_\beta (\beta\<x\> \cdot t ~;_\gamma u) \ar@{=>}[d] \\
(s ~;_\beta (x \mapsto t)f) ~;_\gamma u \ar@{|=|}[d]
  & \beta[f] \cdot s ~;_\beta \beta\<x\> \cdot (t ~;_\gamma u) \ar@{=>}[d] \\
(x \mapsto s ~;_\beta t)f ~;_\gamma u \ar@{|=|}[d]
  & s ~;_\beta (x \mapsto (t ~;_\gamma u))f \ar@{|=|}[d] \\
(x \mapsto (s ~;_\beta t) ~;_\gamma u)f \ar@{|=|}[r]
  & (x \mapsto s ~;_\beta (t ~;_\gamma u))f  \\
}
\]
Similar arguments are now easily inferred for the term in the last position.

\subsection{The natural transformations}

At this stage we have demonstrated that we have the basic functorial data
for a linear actegory, namely a (distributive) monoidal category acting on
a linearly distributive category in a covariant and contravariant way. It
remains to show that the coherent transformations are present and satisfy
the required conditions.  

To accomplish this task we shall indicate the definition of the natural
transformations and illustrate how one establishes their naturality. Here
is how one natural transformation from each of the major symmetry classes
is defined in the term logic:
\begin{align*}
r_\oplus &:= \beta\<\beta_1,\beta_2\>\cdot\beta_2 \ainit\cdot \alpha =_X \beta_1
    \\[1ex]
a_\ox &:= \alpha\<\alpha_1,\alpha_2\> \cdot \alpha_1\<\alpha_{11},\alpha_{12}\>
    \cdot \beta\scase{\beta_1}{\alpha_{11} =_X \beta_1}
               {\beta_2}{\beta_2\scase{\beta_{21}}{\alpha_{12} =_Y \beta_{21}}
                                     {\beta_{22}}{\alpha_2 =_Z \beta_{22}}}
    \\[1ex]
l_\ox &:= \alpha\<\alpha_1,\alpha_2\> \cdot \alpha_1\ainit \cdot
         \alpha_2 =_X \beta
    \\[1ex]
r_\ox &:= \alpha\<\alpha_1,\alpha_2\> \cdot \alpha_2\ainit \cdot
    \alpha_1 =_X \beta
    \\[1ex]
c_\ox &:= \alpha\<\alpha_1,\alpha_2\> \cdot \beta
    \scase{\beta_1}{\alpha_2 =_Y \beta_1}{\beta_2}{\alpha_1 =_X \beta_2}
    \\[1ex]
u_\wdot &:= \alpha\<\rinit\> \cdot \alpha =_X \beta
    \\[1ex]
a^*_\wdot &:= \alpha\<(x,y)\> \cdot \beta[x] \cdot \beta[y] \cdot
    \alpha =_X \beta
    \\[1ex]
a^\wdot_\ox &:= \alpha\<x\> \cdot \alpha\<\alpha_1,\alpha_2\> \cdot 
    \beta\scase{\beta_1}{\beta_1[x] \cdot \alpha_1 =_X \beta_1}
               {\beta_2}{\alpha_2 =_X \beta_2}
    \\[1ex]
d^\ox_\parr &:= \beta\<\beta_1,\beta_2\> \cdot \alpha\< \alpha_1,\alpha_2 \>
    \cdot \alpha_2\scase{\alpha_{21}}{\beta
            \scase{\beta_1}{\alpha_1=_X \beta_{11}}
                  {\beta_2}{\alpha_{21}=_Y \beta_{12}}
              }{\alpha_{22}}{\alpha_{22} =_Z \beta_2}
    \\[1ex]
d^\wdot_\oplus &:= \beta\<\beta_1,\beta_2\> \cdot \alpha\<x\> \cdot
    \beta_1[x] \cdot
    \alpha\scase{\alpha_1}{\alpha_1 =_X \beta_1}
                {\alpha_2}{\alpha_2 =_X \beta_2}
    \\[1ex]
d^\wdot_\bdot &:= \alpha\<x\> \cdot \beta\<y\> \cdot \beta[x] \cdot \alpha[y]
    \cdot \alpha_1 =_X \beta_1
    \\[1ex]
n &:= \beta\<x\>\cdot\beta[x]\cdot \alpha =_X \beta
\end{align*}

The natural transformations labeled with $a$ should be isomorphisms.
In particular $a^\wdot_\ox$ should be an isomorphism. It is not difficult to
check that its inverse is given by
\[
(a^{\wdot}_\ox)^{-1} := \alpha\<\alpha_1,\alpha_2\> \cdot \alpha_1\<x\> \cdot \beta[x] \cdot
    \beta\scase{\beta_1}{\alpha_1 =_X \beta_1}
               {\beta_2}{\alpha_2 =_X \beta_2}~.
\]
To show that these transformations are natural we must check that
naturality conditions. We shall demonstrate what is involved for
$a^\wdot_\ox$, i.e., given $f:A \ra B$, $s:W \ra Y$, and $t:X \ra Z$, that
the following categorical diagram commutes.
\[
\xymatrix@R=6ex{
    A \circ (W \ox X) \ar[d]_{f \circ (s \ox t)}\ar[r]^-{a^\wdot_\ox} 
    & (A \circ W) \ox X \ar[d]^{(f \circ s) \ox t} \\
    B \circ (Y \ox Z) \ar[r]^-{a^\wdot_\ox} & (B \circ Y) \ox Z}
\]
We shall translate the upper route into a term and show that (without looking
into $f$, $s$, or $t$) we can manipulate it into a form which is equivalent to
the lower route. When we translate the terms $s$ and $t$ we shall indicate
that it runs from channel $\alpha$ to $\beta$ (i.e.,
$s::\alpha:W \ra \beta:Y$) by labeling it $s[\alpha;\beta]$. The top route
then gives
\begin{align*}
& \alpha\<x\>  \cdot \alpha\<\alpha_1,\alpha_2\> \cdot
    \begin{array}[t]{l}
        \beta\scase{\beta_1}{\beta_1[x]\cdot\alpha_1 =_W \beta_1}
                {\beta_2}{\alpha_2 =_X \beta_2}~~;_\beta
        \\ \qquad\qquad
        \beta\<\beta_1,\beta_2\> \cdot
        \gamma\scase{\gamma_1}{\beta_1\<x\>\cdot\gamma_1[f(x)]\cdot
        s[\beta_1;\gamma_1]}{\gamma_2}{t[\beta_2;\gamma_2]}
    \end{array}
    \\[1ex]
&\Ra~ \alpha\<x\>\cdot\alpha\<\alpha_1,\alpha_2\> \cdot
    (\beta_1[x] \cdot \alpha_1 =_W \beta_1
    \begin{array}[t]{l}
     ~;_{\beta_1} \alpha_2 =_X \beta_2 ~;_{\beta_2} \\
    \gamma\scase{\gamma_1}{\beta_1\<x\>\cdot\gamma_1[f(x)]
    \cdot s[\beta_1;\gamma_1]}{\gamma_2}{t[\beta_2;\gamma_2]}
    \end{array} \\[1ex]
&\Ra~ \alpha\<x\>\cdot\alpha\<\alpha_1,\alpha_2\> \cdot \gamma
     \scase{\gamma_1}{\beta_1[x] \cdot \alpha_1 =_W \beta_1 ~;_{\beta_1}
        \beta_1\<x\>\cdot\gamma_1[f(x)]\cdot s[\beta_1;\gamma_1]
           }{\gamma_2}{\alpha_2 =_X \beta_2 ~;_{\beta_2}t[\beta_2;\gamma_2]}
\\[1ex]
&\Ra~ \alpha\<x\>\cdot\alpha\<\alpha_1,\alpha_2\> \cdot \gamma
    \scase{\gamma_1}{\gamma_1[f(x)]\cdot s[\alpha_1;\gamma_1]}
          {\gamma_2}{t[\alpha_2;\gamma_2]}~,
\end{align*}
and the bottom route
\begin{align*}
& \alpha\<x\> \cdot \beta[f(x)] \cdot \alpha\<\alpha_1,\alpha_2\> \cdot \beta
    \begin{array}[t]{l}
    \scase{\beta_1}{s[\alpha_1;\beta_1]}{\beta_2}{t[\alpha_2;\beta_2]} 
    ~;_\beta
    \\ \qquad
    \beta\<y\> \cdot \beta\<\beta_1,\beta_2\> \cdot
    \gamma\scase{\gamma_1}{\gamma_1[y]\cdot \beta_1 =_Y \gamma_1}
                {\gamma_2}{\beta_2 =_Z \gamma_2}
    \end{array}
    \\[1ex]
&\Ra~ \alpha\<x\> \cdot \bigg(y \mapsto \alpha\<\alpha_1,\alpha_2\> \cdot 
    \begin{array}[t]{l}
        \beta\scase{\beta_1}{s[\alpha_1;\beta_1]}
                   {\beta_2}{t[\alpha_2;\beta_2]} ~;_\beta \\ \qquad
        \beta\<\beta_1,\beta_2\> \cdot \gamma
        \scase{\gamma_1}{\gamma_1[y] \cdot \beta_1 =_Y \gamma_1}
              {\gamma_2}{\beta_2 =_Z \gamma_2}
        \bigg) f(x)
    \end{array}
    \\[1ex]
&\Ra~ \alpha\<x\>\cdot \bigg(\alpha\<\alpha_1,\alpha_2\> \cdot 
    \begin{array}[t]{l}
        \beta\scase{\beta_1}{s[\alpha_1;\beta_1]}
                   {\beta_2}{t[\alpha_2;\beta_2]} ~;_\beta \\ \qquad
        \beta\<\beta_1,\beta_2\>\cdot \gamma
        \scase{\gamma_1}{\gamma_1[f(x)]\cdot \beta_1 =_Y \gamma_1}
              {\gamma_2}{\beta_2 =_Z \gamma_2}
        \bigg)
    \end{array}
    \\[1ex]
&\Ra~ \alpha\<x\> \!\cdot\! \alpha\<\alpha_1,\alpha_2\> \!\cdot\! \bigg(\!
    s[\alpha_1;\beta_1] ~;_{\beta_1} \bigg(t[\alpha_2;\beta_2] ~;_{\beta_2}
    \gamma \scase{\gamma_1}{\gamma_1[f(x)]\cdot \beta_1 =_Y \gamma_1}
                 {\gamma_2}{\beta_2 =_Z \gamma_2}\!\bigg)\!\bigg)
    \\[1ex]
&\Ra~ \alpha\<x\>\cdot\alpha\<\alpha_1,\alpha_2\> \cdot \gamma
    \scase{\gamma_1}{\gamma_1[f(x)]\cdot s[\alpha_1;\beta_1] ~;_{\beta_1}
                     \beta_1 =_Y \gamma_1}
          {\gamma_2}{t[\alpha_2;\beta_2] ~;_{\beta_2} \beta_2 =_Z \gamma_2}
    \\[1ex]
&\Ra~ \alpha\<x\>\cdot\alpha\<\alpha_1,\alpha_2\> \cdot 
        \gamma\scase{\gamma_1}{\gamma_1[f(x)]\cdot s[\alpha_1;\gamma_1]
           }{\gamma_2}{t[\alpha_2;\gamma_2]}
\end{align*}
which therefore proves the naturality of $a^\wdot_\ox$.

\subsection{Completeness}

In order to establish completeness of the logic it is now necessary to check
that all the coherence diagrams commute. This is a lengthy exercise most
of which we will leave to the reader!

We shall explicitly check the triangle equalities for the parameterised
adjunction.  Because of the symmetry it actually suffices to check just one:
\[
A \wdot n_X ; e_{A \wdot X} = 1_{A \wdot X}: A \wdot X \ra A \wdot X.
\]
Here is the explicit calculation:
\begin{align*}
A \wdot  n_X & ~;_\beta e_{A \wdot X}
\\[1ex]
    =~~ & \alpha\<w\> \cdot \beta[w] \cdot \beta\<x\> \cdot \beta[x] \cdot
    \alpha =_X \beta ~~;_\beta~ \beta\<y\> \cdot \beta[y] \cdot \beta\<z\>
    \cdot \gamma[z] \cdot \beta =_X \gamma
\\[1ex]
    \Ra & \alpha\<w\>\cdot \big(\beta\<x\> \cdot \beta[x] \cdot \alpha =_X \beta
    ~~;_\beta~ \beta[w] \cdot \beta\<z\> \cdot \gamma[z] \cdot \beta =_X
    \gamma \big)
\\[1ex]
    \Ra & \alpha\<w\> \cdot \big(\beta[w] \cdot \alpha =_X \beta
    ~~;_\beta~ \beta\<z\> \cdot \gamma[z] \cdot \beta =_X \gamma \big)
\\[1ex]
    \Ra & \alpha\<w\> \cdot \big(\alpha =_X \beta
    ~~;_\beta~ \gamma[w] \cdot \beta =_X \gamma \big)
\\[1ex]
    \Ra & \alpha\<w\> \cdot \gamma[w] \cdot \alpha =_X \gamma = 1_{A \wdot X}.
\end{align*}

The other aspect of this setting we have not discussed is the coproducts. We 
expect that $(A+B) \wdot X$ is the coproduct of $A \wdot X$ and $B \wdot X$.  
If we have proofs 
\[
s_1:: \emptyset \mid \Gamma,\alpha:A \wdot X \Vd \Delta
\qquad \text{and} \qquad
s_2:: \emptyset \mid \Gamma,\alpha:B \wdot X \Vd \Delta
\]
it is not hard to see (especially using the representability results
discussed below) that one can construct a proof
\[
s:: \emptyset \mid \Gamma,\alpha:(A+B) \wdot X \Vd \Delta.
\]
However if one has a proof $s$ as above how does one see it as a cotuple of
proofs?

Consider the effect of substituting on $\alpha$ the identity map of
$(A + B) \wdot X$ into $s$ (where the first reduction is in reverse):
\begin{align*}
s   &~\Ra~ \beta =_{(A+B) \wdot X} \alpha ~;_\alpha s \\[1ex]
    &~\Ra~ \beta\<z\> \cdot \alpha\left[\case{x}{\sigma_0(x)}{y}
        {\sigma_1(y)} z\right] \cdot \beta =_X \alpha ~;_\alpha s \\[1ex]
    &~\Ra~ \beta\<z\> \cdot
       \case{x}{\alpha[\sigma_0(x)] \cdot \beta =_X \alpha ~;_\alpha s}
            {y}{\alpha[\sigma_1(y)] \cdot \beta =_X \alpha ~;_\alpha s} \\[1ex]
    &~\Ra~ \beta\<z\> \cdot \case{x}{\alpha[\sigma_0(x)] \cdot s}
                               {y}{\alpha[\sigma_1(y)] \cdot s} 
\end{align*}
This shows how the term $s$ can be decomposed as a cotuple of the composites
with the injections. Thus $(A+B) \wdot X$ is the coproduct of $A \wdot X$ and
$B \wdot X$.

This completes our discussion of the completeness of the logic. We have
shown that:

\begin{proposition}
The terms of the message passing logic between single types with cut as
composition form a linear additive actegory.
\end{proposition}

In fact, if $\A$ is any monoidal category this shows how we can construct
a linearly distributive category from $\A$ and the empty poly-$\A$-actegory
(i.e., the initial linear $\A$-actegory). The result is, of course,
definitely a non-empty linear actegory: the units $\top$ and $\bot$ must
always be present which implies that $A \circ \bot$ and $A \bullet \top$ for
all $A \in \A$ are non-trivial objects. It is also not hard to see that
this will be a $*$-autonomous category as, inductively $(A \wdot X)^* =
A \bdot X^*$ and $(A \bdot Y)^* = A \wdot Y^*$, with the base case the
units.

\section{Representability and soundness}\label{sec-rep}

In the theory of polycategories representability plays a crucial role 
in getting between the purely categorical (object-to-object) view and the 
circuit (polycategorical) view of the maps. In fact, to be the category of
maps\footnote{In this section and in the sequel we will call a polymap or
multimap with singleton domain and singleton codomain simply a \emph{map}.}
of a representable polycategory is precisely to be a linearly distributive
category. Similarly to be the category of maps of a representable
multicategory is precisely to be a monoidal category.

For the message passing logic an exactly analogous correspondence holds:

\begin{theorem} \label{thm-main}
To be the category of maps of a representable (additive) poly-actegory is
precisely to be a linear (additive) actegory.
\end{theorem}

Instead of a polycategory we must start with a polycategory with a
multicategorical action: a \emph{poly-actegory}. This makes the polymaps have
a type which match the sequent structure of the message passing logic:
composition in a poly-actegory is given by the two sorts of cut and the
evident associativity and interchange laws must hold. Representability then
entails the representability of all the features of the logic: the tensor,
the par, and the two actions. To be a linear actegory is then precisely to
be the maps in a representable poly-actegory.  

The soundness of the term calculus is therefore determined by its soundness 
in any representable poly-actegory. However, the equations of the term 
calculus were developed from the proof theory in Section~\ref{sec-logic-mc} 
and Section~\ref{sec-logic-mp}.  Thus with the identification of the 
proof theory with poly-actegories we have the soundness by construction. 

The purpose of this section is to introduce poly-actegories as a formulation
of the proof theory of the message passing logic and to prove
Theorem~\ref{thm-main}.

\subsection{Poly-actegories and circuit representation}

The proof theory for the two-sided cut rule lies in polycategories.  
Essentially they are the natural deduction style proofs for the system. 
Polycategories have an extremely intuitive representation using circuits
(see~\cite{BCST}). This gives further indication that the proof theory and 
its categorical formulation is capturing a very natural phenomenon. 

We wish to show how the proof theory of message passing, as represented 
by poly-actegories, may also be represented using circuits.
Every proof in multiplicative linear logic has a presentation as a 
circuit~\cite{BCST}. To extend this to the message passing
logic we must describe how the two action rules are to be interpreted
using circuits. 

A multicategory $\M$ has objects $\ob(\M)$ and multimaps $x:\Phi \ra B$
where $\Phi \subset \ob(\M)$. A polycategory $\P$ has objects $\ob(\P)$ and
polymaps $u:\Gamma \ra \Delta$ where $\Gamma,\Delta \subset \ob(\P)$.
Multimaps and polymaps are represented in the circuit notation, where
suppose $\Phi = A_1,\ldots,A_m$, $\Gamma = X_1,\ldots,X_m$, and
$\Delta = Y_1,\ldots,Y_n$, as
\[
\scalebox{0.85}{\rb{-\height/2}{\ig[scale=\scale]{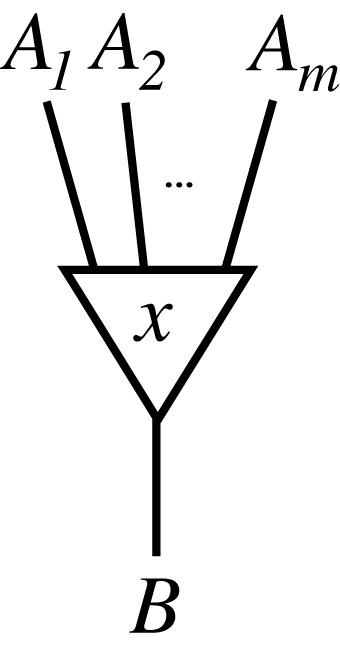}}}
\qquad \text{and} \qquad
\scalebox{0.85}{\rb{-\height/2}{\ig[scale=\scale]{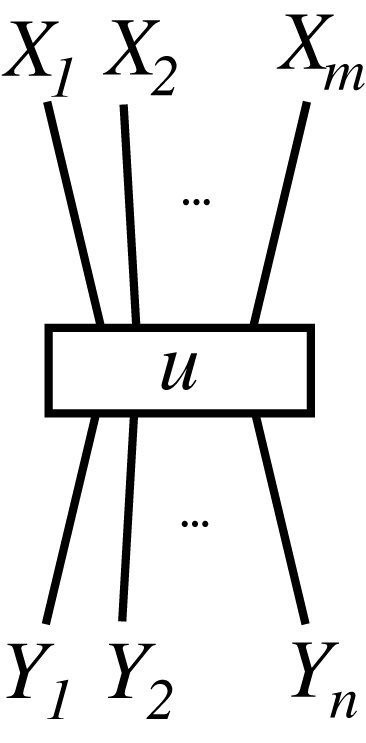}}}
\]
respectively. Since we consider symmetric multicategories and symmetric
polycategories this allows the ``wires'' to cross. Often, to simplify the
circuit notation, we will use a double line to indicate a possibly empty
subset of objects as in 
\[
    \scalebox{0.85}{\rb{-\height/2}{\ig[scale=\scale]{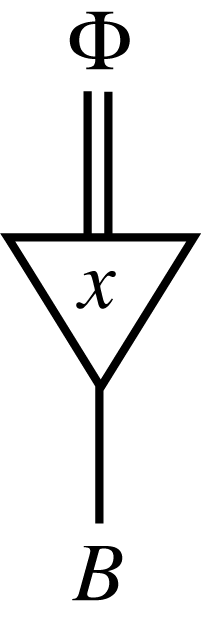}}}
    \qquad \text{and} \qquad
    \scalebox{0.85}{\rb{-\height/2}{\ig[scale=\scale]{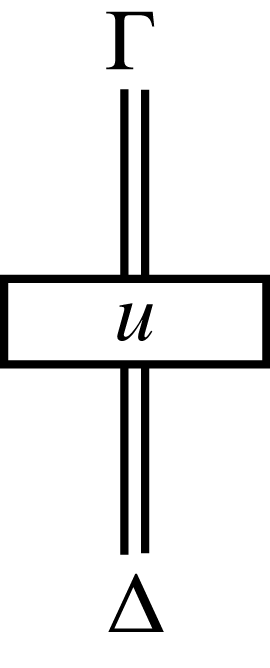}}}~.
\]
A \emph{poly-actegory} $\X$ is a polycategory $\P$ acted on by a multicategory
$\M$. What this means is that the polymaps in $\X$, instead of the usual
polymaps as in $\P$ above, will contain inputs with types coming from both
$\P$ and $\M$. A typical example is $s: \Phi \mid \Gamma \ra \Delta$. These
sorts of polymaps may also be represented using the circuit notation as
\[
\scalebox{0.85}{\rb{-\height/2}{\ig[scale=\scale]{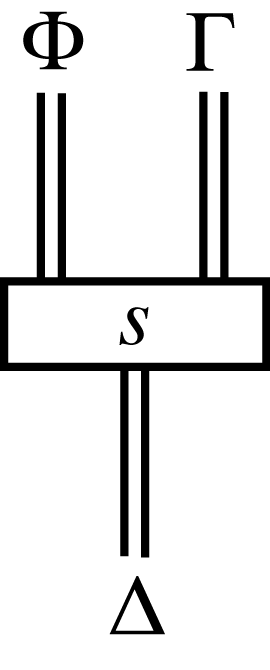}}}
\]
where the type of a wire indicates whether it is from $\M$ or $\P$.

There are two types of composition in a poly-actegory. The first composes a
multimap in $\M$ with a polymap in $\X$. Given a multimap $f:\Phi \ra A$
in $\M$ and a polymap $s:\Psi,A,\Psi' \mid \Gamma \ra \Delta$ in $\X$.
Composing (on $A$) results in a polymap in $\X$ of type
\[
    f \cdot s :\Psi,\Phi,\Psi' \mid \Gamma \ra \Delta.
\]
where the centre dot notation ``\,$\cdot$\,'' represents cutting a multi\-map
into a poly\-map. In circuit notation this is represented as
\[
    \scalebox{0.85}{\rb{-\height/2}{\ig[scale=\scale]{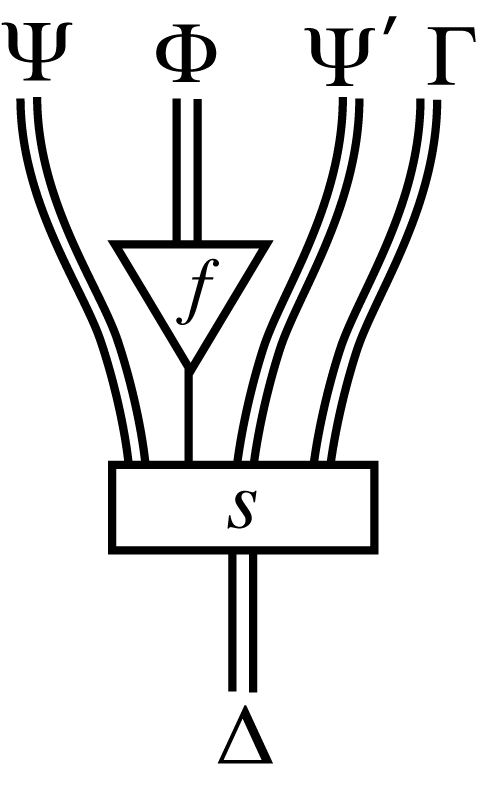}}}~.
\]
The second type of composition is between polymaps in $\X$. Suppose
$s:\Phi \mid \Gamma \ra \Delta_1,X,\Delta_2$ and $t:\Psi \mid
\Gamma_1,X,\Gamma_2 \ra \Delta$ are polymaps in $\X$. Composing (on $X$)
results in a polymap in $\X$ of type
\[
s\, ;\, t :\Phi,\Psi \mid \Gamma_1,\Gamma,\Gamma_2 \ra \Delta_1,\Delta,\Delta_2.
\]
One possible simplified circuit diagram (in which $\Gamma_2 = \Delta_1 =
\emptyset$ so that there are no crossings) for this cut is given by
\[
    \scalebox{0.85}{\rb{-\height/2}{\ig[scale=\scale]{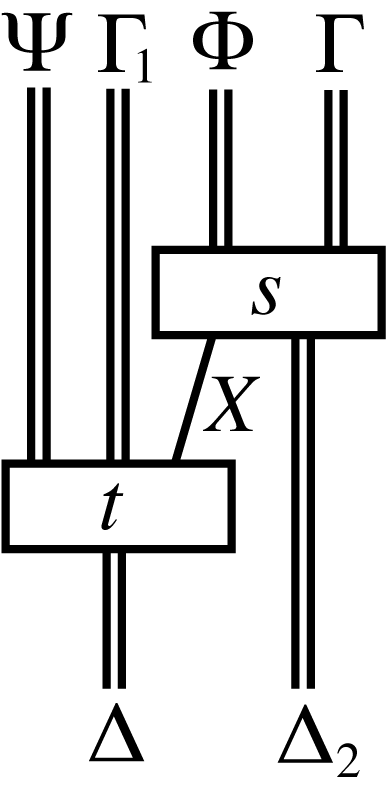}}}~.
\]

All of the evident associativity and interchange laws must, of course, hold
in any poly-actegory.

\begin{figure}
\begin{center}
\framebox[\boxwidth]{
\begin{tabular}{cc}
    ($\wdot~I$) \scalebox{0.85}{\rb{-\height/2}{\ig[scale=\scale]{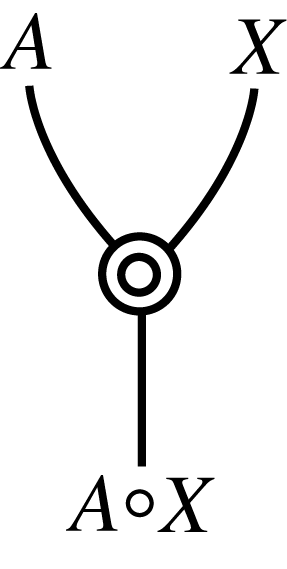}}}
    \qquad\qquad & \qquad\qquad
    ($\wdot~E$) \scalebox{0.85}{\rb{-\height/2}{\ig[scale=\scale]{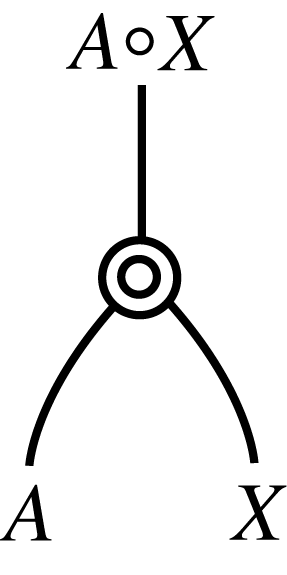}}}
\\[10ex]
    ($\bdot~I$) \scalebox{0.85}{\rb{-\height/2}{\ig[scale=\scale]{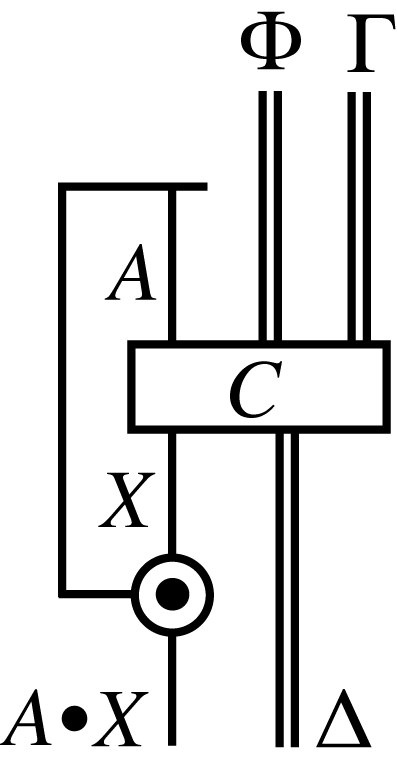}}}
    \qquad\qquad & \qquad\qquad
    ($\bdot~E$) \scalebox{0.85}{\rb{-\height/2}{\ig[scale=\scale]{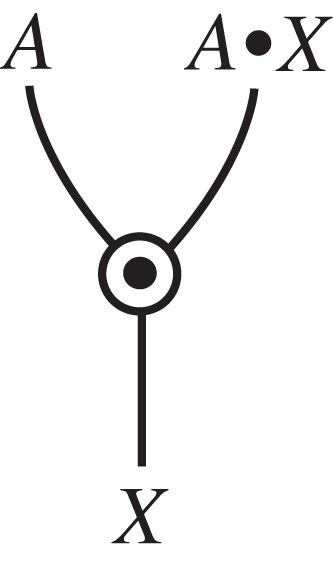}}}
\\[12ex]
\multicolumn{2}{c}{where $C$ in ($\bdot~I$) is a circuit.}
\end{tabular}
}\end{center}
\caption{$\wdot$ and $\bdot$ circuit introduction and elimination rules}
\label{fig-circuit}
\end{figure}

There are circuits corresponding to introduction and elimination rules for
each of the connectives $\ox$, $\oplus$, $\top$, and $\bot$. We refer the
interested reader to~\cite{BCST} for the full story. Here we are interested
in the action rules $\wdot$ and $\bdot$. The introduction and elimination
circuit diagrams for $\wdot$, which we label by ($\wdot~I$) and ($\wdot~E$),
are given in Figure~\ref{fig-circuit}. Notice the similarity to the
introduction and elimination rules of the $\ox$ and $\oplus$
connectives~\cite{BCST}.

The $\bdot$ rule is a binding rule in the sense that the introduction rule
must involve a ``scope box''~\cite{fill} so that they are only applicable to
the situation where one has a subcircuit $C$ to attach the link to. It
replaces a derivation $A, \Phi \mid \Gamma \Vd X,\Delta$ with a derivation
$\Phi \mid \Gamma \Vd A \bdot X, \Delta$. The circuit diagrams, labeled by
($\bdot~I$) and ($\bdot~E$), are also given in Figure~\ref{fig-circuit}.
With this rule notice the similarity to the introduction and elimination rules
(including the scope box) for $\multimap$, the linear implication~\cite{fill}.

\begin{figure}
\begin{center}
\framebox[\boxwidth]{
\scalebox{0.85}{
\begin{tabular}{cc}
    \scalebox{0.85}{\rb{-\height/2}{\ig[scale=\scale]{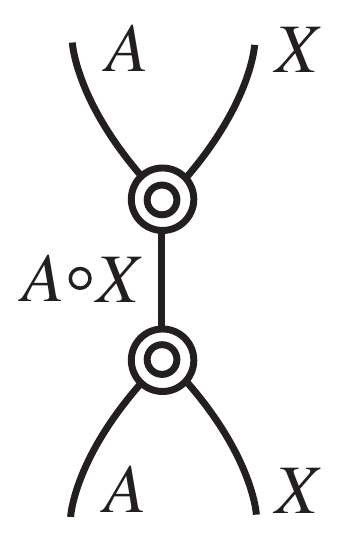}}} \Ra~
    \scalebox{0.85}{\rb{-\height/2}{\ig[scale=\scale]{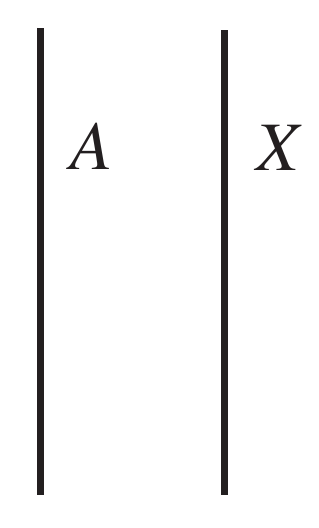}}}
    \qquad & \qquad
    \scalebox{0.85}{\rb{-\height/2}{\ig[scale=\scale]{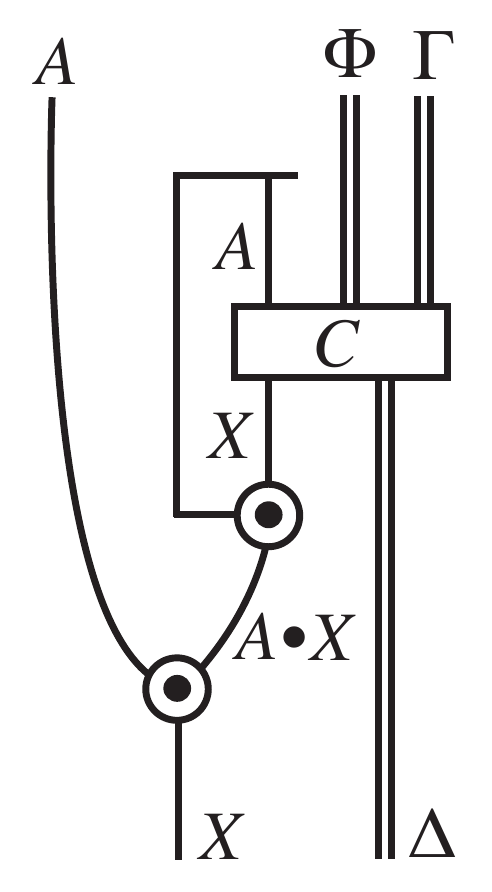}}} \Ra~
    \scalebox{0.85}{\rb{-\height/2}{\ig[scale=\scale]{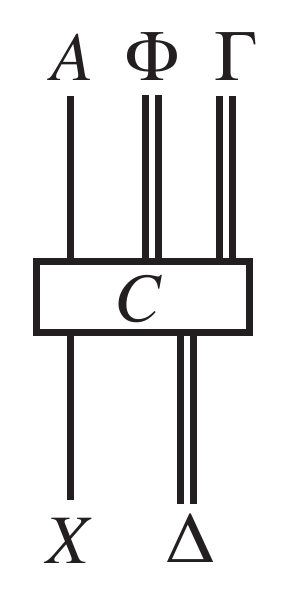}}}
\\
    \scalebox{0.85}{\rb{-\height/2}{\ig[scale=\scale]{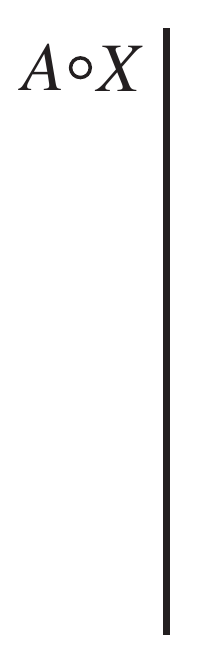}}} ~\Ra~
    \scalebox{0.85}{\rb{-\height/2}{\ig[scale=\scale]{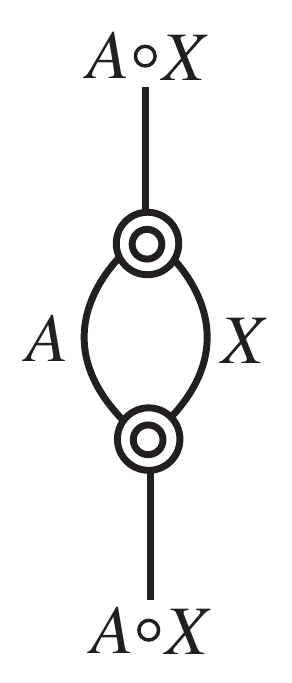}}}
    \qquad & \qquad
    \scalebox{0.85}{\rb{-\height/2}{\ig[scale=\scale]{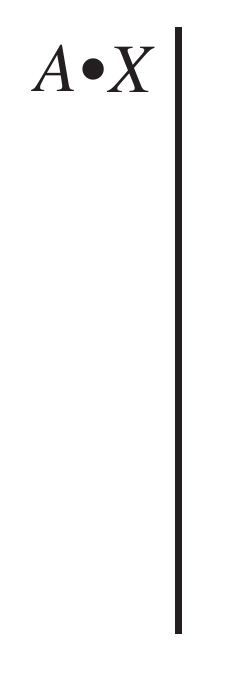}}} ~\Ra~
    \scalebox{0.85}{\rb{-\height/2}{\ig[scale=\scale]{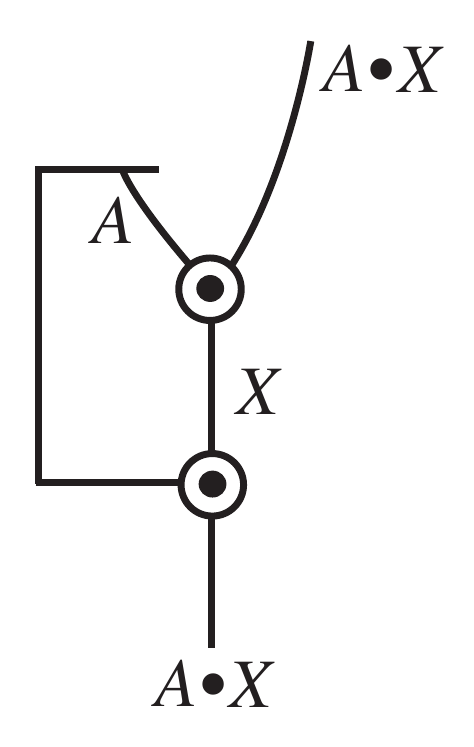}}}
\end{tabular}
}}\end{center}
\caption{Circuit reduction and expansion rules}
\label{fig-circ-re}
\end{figure}

For any connective there are two types of circuit rewrites: \emph{reductions}
which allow one to simplify a circuit which involves an elimination rule
immediately after an introduction rule, and \emph{expansions} which
``split'' a wire carrying a compound formula into ``simpler'' wires, and
ultimately to atomic wires.  The reductions and expansions for $\wdot$ and
$\bdot$ are given by in Figure~\ref{fig-circ-re}.

\begin{figure}
\begin{center}
\framebox[\boxwidth]{
\scalebox{0.75}{
\rb{-\height/2}{\ig[scale=\scale]{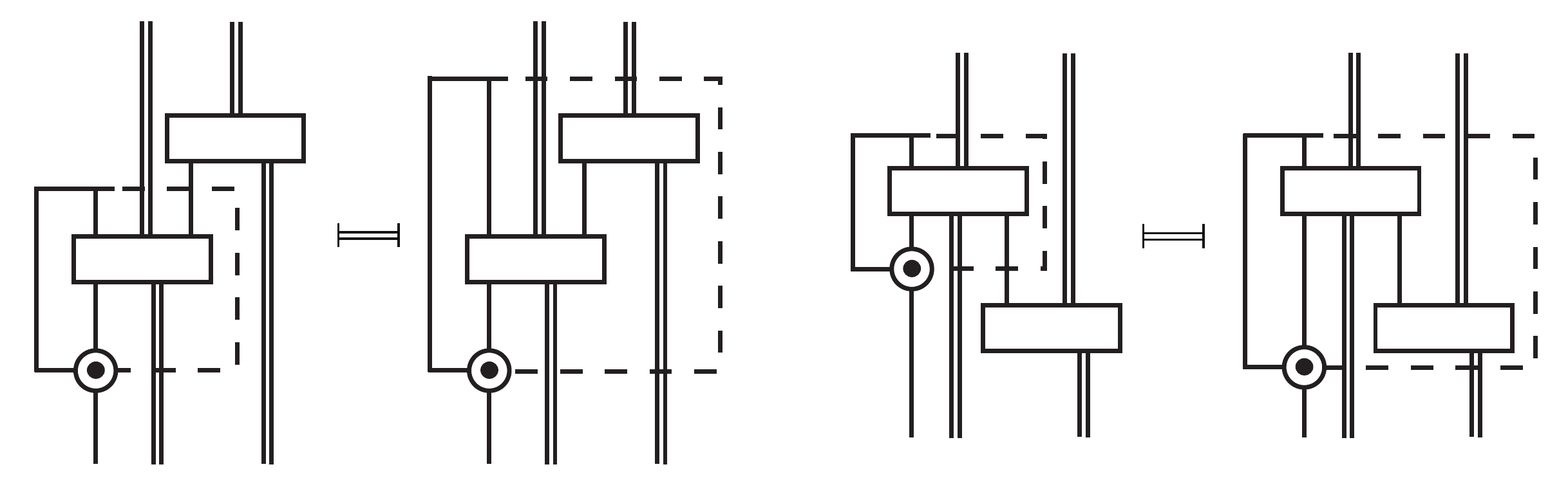}}
}}\end{center}
\caption{$\bdot$ scope rules}
\label{fig-circ-scope}
\end{figure}

In circuits, links are allowed to slide up or down corresponding to cut
elimination steps or equations. For example, in simplified circuit notation,
sliding of the $\wdot$ node
\[
\scalebox{0.85}{\rb{-\height/2}{\ig[scale=\scale]{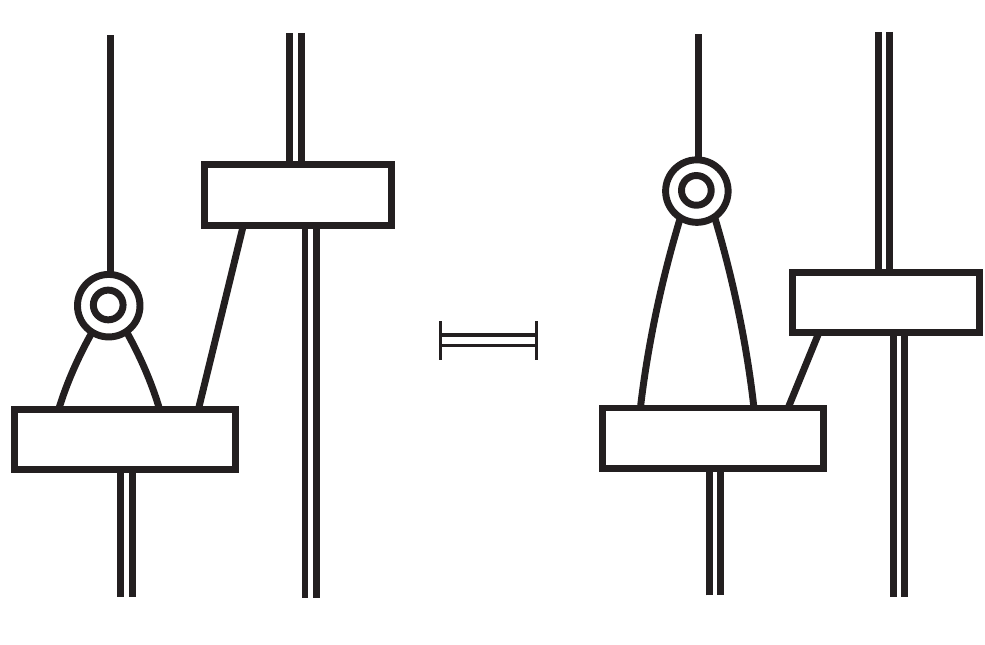}}}
\]
corresponds to the [sequent-$\wdot_l$] cut elimination rule. This should be
familiar, however, the scoping rules for $\bdot$ may not be. We present the
equations (where we have again slightly simplified the circuit notation) in
Figure~\ref{fig-circ-scope}. 

\subsection{Representability of the actions}

There are natural bijections
\[
    \AxiomC{$\Phi,A \mid \Gamma,X\Vd \Delta$}
    \doubleLine
    \UnaryInfC{$\Phi \mid \Gamma,A \wdot X \Vd \Delta$}
    \DisplayProof
\qquad \text{and} \qquad
    \AxiomC{$\Phi,A \mid \Gamma \Vd Y, \Delta$}
    \doubleLine
    \UnaryInfC{$\Phi \mid \Gamma \Vd A \bdot Y, \Delta$}
    \DisplayProof
\]
obtained in the top-to-bottom direction respectively by the $\wdot_l$ and
$\bdot_r$ inference rules, and in the bottom-to-top direction by cutting
respectively with the derivations
\[
    \AxiomC{$\text{id}_A$}
    \UnaryInfC{$A \vd A$}
    \AxiomC{$\text{id}_X$}
    \UnaryInfC{$\emptyset \mid X \Vd X$}
    \BinaryInfC{$A \mid X \Vd A \wdot X$}
    \DisplayProof
\qquad \text{and} \qquad
    \AxiomC{$\text{id}_A$}
    \UnaryInfC{$A \vd A$}
    \AxiomC{$\text{id}_X$}
    \UnaryInfC{$\emptyset \mid X \Vd X$}
    \BinaryInfC{$A \mid A \bdot X \Vd X$}
    \DisplayProof~.
\]
These are called the \emph{representing polymaps} for the covariant and
contravariant actions respectively. They are clearly natural as when
represented they are the identity maps.

The representing polymap of the covariant action may be described as a term:
\[
\alpha[x] \cdot \gamma =_X \alpha ~~::~~ x:A \mid \gamma:X \Vd
\alpha:A \wdot X.
\]
which allows one to see the bijection at the level of the term logic. Suppose
\[
s::\Phi,x:A \mid \Gamma,\alpha:X \Vd \Delta.
\]
The calculation
\begin{align*}
\alpha[x] \cdot \gamma =_X \alpha ~~;_\alpha~ \alpha\<y\> \cdot s
    & \Ra \gamma =_X \alpha ~~;_\alpha~ s \\
    & \Ra s,
\end{align*}
where the first rewrite is the [$\wdot_r$-$\wdot_l$] cut elimination step,
establishes one direction of the bijection. Now suppose
\[
    t::\Phi \mid \Gamma,\alpha:A \wdot X \Vd \Delta.
\]
The other direction is given by
\begin{align*}
\gamma\<x\> \cdot (\alpha[x] \cdot \gamma =_X \alpha ~;_\alpha~ t)
    & ~\Ra~ (\gamma\<x\>\cdot \alpha[x] \cdot \gamma =_X \alpha) ~;_\alpha~ t \\
    & ~\Ra~ \gamma =_{A \wdot X} \alpha ~;_\alpha~ t \\
    & ~\Ra~ t
\end{align*}
where the first rewrite is the [$\wdot_l$-sequent] cut elimination step in
reverse, and the second uses that $\gamma\<x\> \cdot \alpha[x]
\cdot \gamma =_X \alpha$ is the identity on $A \wdot X$.~

\subsection{Representability for poly-actegories}

The proof theory of the message passing logic is a poly-actegory.
In this setting each polymap, besides having multiple inputs and outputs
can, in addition, have inputs from a multicategorical world. Although the
inputs may be typed to come from different worlds, once this
distinction is erased one is simply left with a polycategory in which
certain types are served by multimaps alone. This view determines the
requirements on the poly-actegory composition. Demanding representability of
the (multicategorical) tensor $*$ and its unit $I$, the (polycategorical)
tensor $\ox$ and its unit $\top$, the par $\parr$ and its unit $\bot$,
the covariant action $\wdot$, and the contravariant action $\bdot$, then
forces, we claim, all the proof equivalences discussed in
Section~\ref{sec-msg-equiv} for the message passing logic.

To prove this would be stretching the patience of the reader and is,
besides, relatively standard categorical proof theory. Instead we shall
focus on how a linear actegory arises from the maps of these settings. 
Given the circuit representation it is actually very straightforward to
verify that the required coherence diagrams are satisfied. Thus, the main
objective of this section is to show how the data of a linear actegory
arises.

A significant feature of a linear actegory is the parameterised adjunction
between the covariant and contravariant action. This arises directly from the
representability by the following two-way series of inferences:
\[
\AxiomC{$A \wdot X \Vd Y$}
\doubleLine
\UnaryInfC{$A \mid X \Vd Y$}
\doubleLine
\UnaryInfC{$X \Vd A \bdot Y$}
\DisplayProof
\]
We can also derive all the coherence isomorphisms using representability.
Here we give the derivation of such for the binary connectives (the units
coherences are derived in a similar manner):
\[
   a^*_\wdot \quad : \quad
   \AxiomC{$(A * B) \wdot X \Vd Y$}
   \doubleLine
   \UnaryInfC{$A * B \mid X \Vd Y$}
   \doubleLine
   \UnaryInfC{$A,B \mid X \Vd Y$}
   \doubleLine
   \UnaryInfC{$A \mid B \wdot X \Vd Y$}
   \doubleLine
   \UnaryInfC{$A \wdot (B \wdot X) \Vd Y$}
   \DisplayProof
\qquad\qquad\qquad
   a^*_\bdot \quad : \quad
   \AxiomC{$X \Vd (A * B) \bdot Y$}
   \doubleLine
   \UnaryInfC{$A * B \mid X \Vd Y$}
   \doubleLine
   \UnaryInfC{$A,B \mid X \Vd Y$}
   \doubleLine
   \UnaryInfC{$A \mid X \Vd B \bdot Y$}
   \doubleLine
   \UnaryInfC{$X \Vd A \bdot (B \bdot Y)$}
   \DisplayProof
\]
\[
   a^\wdot_\ox \quad : \quad
   \AxiomC{$A \wdot (X \ox Y) \Vd Z$}
   \doubleLine
   \UnaryInfC{$A \mid X \ox Y \Vd Z$}
   \doubleLine
   \UnaryInfC{$A \mid X,Y \Vd Z$}
   \doubleLine
   \UnaryInfC{$A \wdot X,Y \Vd Z$}
   \doubleLine
   \UnaryInfC{$(A \wdot X) \ox Y \Vd Z$}
   \DisplayProof
\qquad\qquad\qquad
   a^\bdot_\parr \quad : \quad
   \AxiomC{$X \Vd A \bdot (Y \parr Z)$}
   \doubleLine
   \UnaryInfC{$A \mid X \Vd Y \parr Z$}
   \doubleLine
   \UnaryInfC{$A \mid X \Vd Y,Z$}
   \doubleLine
   \UnaryInfC{$X \Vd A \bdot Y, Z$}
   \doubleLine
   \UnaryInfC{$X \Vd (A \bdot X) \parr Z$}
   \DisplayProof
\]

To derive the distributions we have to use the poly-actegorical composition
together with the representing polymaps. The derivations are below, however,
note that these are not derivations of the logic, but in a representable
poly-actegory. In this setting we have just the poly-actegorical
composition, the representing polymaps, and the equivalences. The binary
inference in each of the inferences below is the poly-actegorical
composition. 
\[
   d^\wdot_\parr :\!\!
   \AxiomC{$X \parr Y \Vd X,Y \!\!\!\!\!$}
   \AxiomC{$A \mid X \Vd A \wdot X$}
   \BinaryInfC{$A \mid X \parr Y \Vd A \wdot X,Y$}
   \doubleLine
   \UnaryInfC{$A \wdot (X \parr Y) \Vd (A \wdot X) \parr Y$}
   \DisplayProof 
\quad
   d^\bdot_\ox :\!\!
   \AxiomC{$A \mid A \bdot X \Vd X \!\!\!\!\!$}
   \AxiomC{$X,Y \Vd X \ox Y$}
   \BinaryInfC{$A \mid A \bdot X, Y \Vd X \ox Y$}
   \doubleLine
   \UnaryInfC{$(A \bdot X) \ox Y \Vd A \wdot (X \ox Y)$}
   \DisplayProof 
\]
\medskip
\[ 
   d^\wdot_\bdot ~:~
   \AxiomC{$B \mid B \bdot X \Vd X$}
   \AxiomC{$A \mid X \Vd A \wdot X$}
   \BinaryInfC{$A,B \mid B \bdot X \Vd A \wdot X$}
   \doubleLine
   \UnaryInfC{$A \wdot (B \bdot X) \Vd B \bdot (A \wdot X)$}
   \DisplayProof
\]

\begin{figure}
\begin{center}
\framebox[\boxwidth]{
\scalebox{0.85}{\rb{-\height/2}{\ig[scale=\scale]{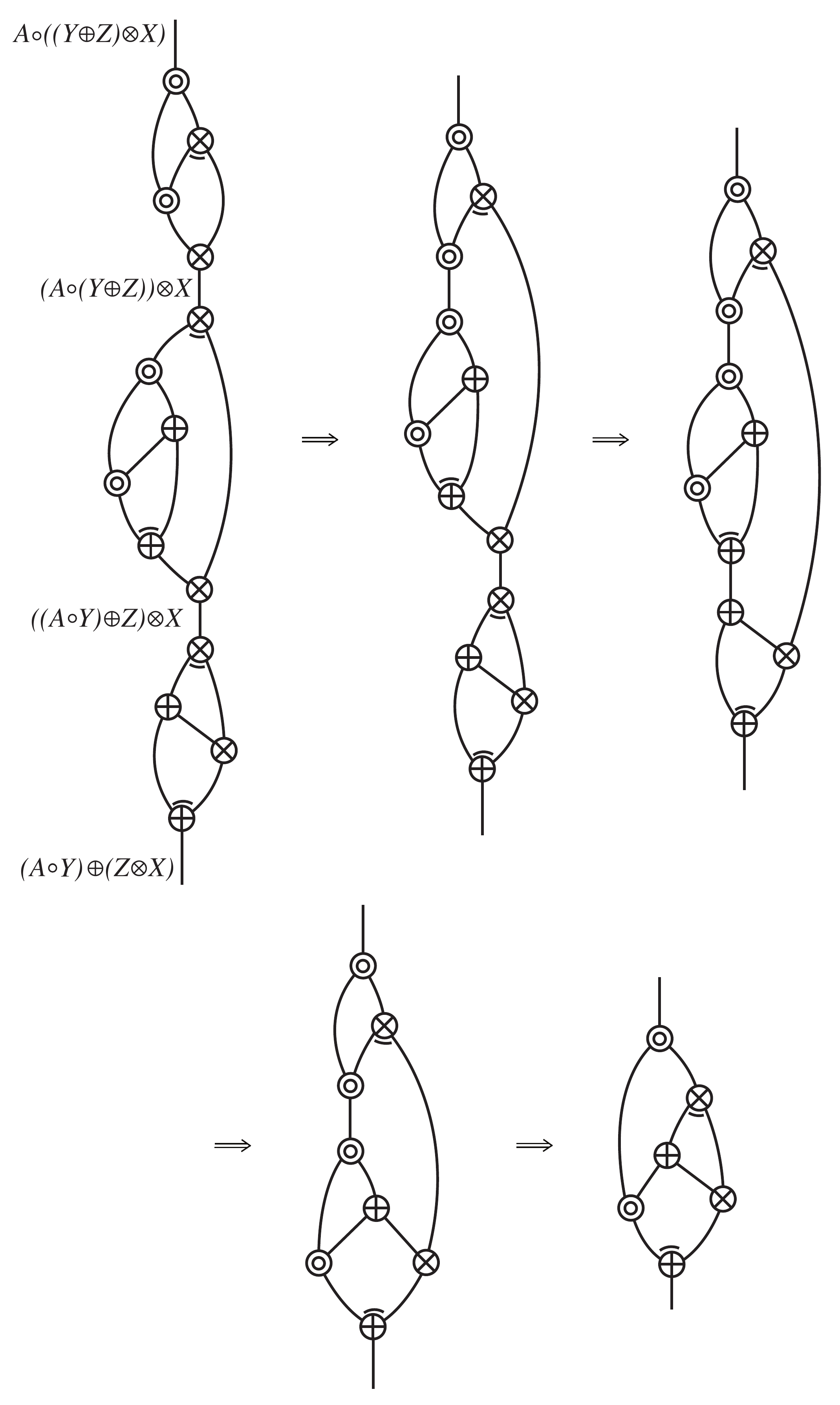}}}
}\end{center}
\caption{Circuit coherence calculation: top route}\label{fig-circ-coh1}
\end{figure}

\begin{figure}
\begin{center}
\framebox[\boxwidth]{
\scalebox{0.85}{\rb{-\height/2}{\ig[scale=\scale]{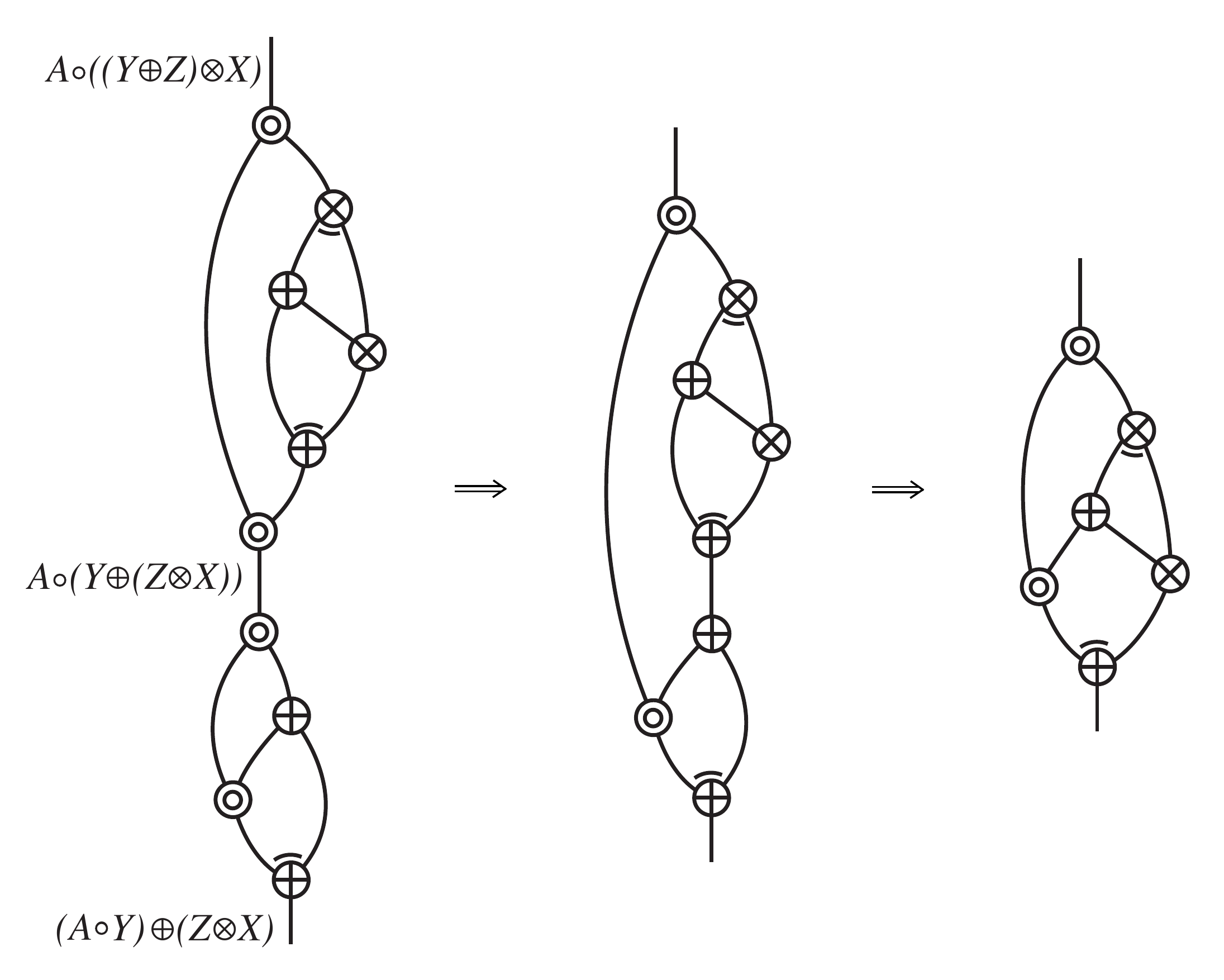}}}
}\end{center}
\caption{Circuit coherence calculation: bottom route}\label{fig-circ-coh2}
\end{figure}

To prove soundness of the message passing logic we have presented a recipe
which relies on the fact (established in the next section) that to be the
category of maps of a representable poly-actegory is precisely to be a
linear actegory (Theorem~\ref{thm-main}). To check that the proof
equivalences of the message passing logic will hold in a representable
poly-actegory is then straightforward. As an example, consider the
coherence diagram (19) which we recall here:
\[
\xygraph{ {A \wdot ((Y \parr Z) \ox X)} (
    :[u(1.7)r(1.4)] {(A \wdot (Y \parr Z)) \ox X} ^-{a^\wdot_\ox}
    :[dr] {((A \wdot Y) \parr Z) \ox X} ^-{d^\wdot_\parr \ox X}
    :[d(1.4)] {(A \wdot Y) \parr (Z \ox X)}="e" ^-{d^\parr_\ox} 
    )
    :[d(1.7)r(1.4)] {A \wdot (Y \parr (Z \ox X))} _-{A \wdot d^\parr_\ox}
    :"e" _-{d^\wdot_\parr}}
\]
It is an easy circuit calculation, which we show in
Figures~\ref{fig-circ-coh1} and~\ref{fig-circ-coh2}, to see that both routes
of the coherence diagram are equivalent.

We conclude:

\begin{proposition}
The category of maps of a representable poly-actegory form a linear
actegory.
\end{proposition}

\subsection{Soundness} \label{sec-soundness}

This section is devoted to showing that all the reductions and equations in
a representable poly-actegory $\P$ (e.g., $\P=\PMsg$) hold in any linear
actegory $\X$. That is, we wish to show that given a linear actegory one may
build from it a representable poly-actegory. Once again we concentrate on the 
action rules and refer the reader to the start of~\cite{CKS} for a description 
of representability, and to~\cite{CS}, for soundness for 
representable polycategories. For simplicity, we will 
ignore most instances of associativity.

The first step is to show how the two cut rules of a poly-actegory arise
from the data of a linear actegory. The ``action cut'' of a multimap into
a polymap will be described first. To this end suppose that we have a
multimap and polymap
\[
f : \Phi \ra A \qquad \text{and} \qquad
s : A,\Psi \mid \Gamma_1,X,\Gamma_2 \ra \Delta
\]
respectively that are represented by the maps (between singletons)
\[
\hat{f}:\hat{\Phi} \ra A \qquad \text{and} \qquad
\hat{s}:\hat{\Gamma}_1 \ox A \wdot X \ox \hat{\Gamma}_2 \ra \hat{\Delta}
\]
in $\X$ (the tildes over $\Gamma$'s and $\Phi$'s indicating the results of
representing). The composite of $f$ and $s$ on $A$ should be a polymap in $\P$ as 
\begin{equation}\label{fdots}\tag{X}
f \cdot s : \Phi,\Psi \mid \Gamma_1,X,\Gamma_2 \ra \Delta.
\end{equation}
In $\X$ we may form the composite 
\[
\xymatrix{\hat{\Gamma}_1 \ox \hat{\Phi} \wdot X \ox \hat{\Gamma}_2
    \ar[rr]^-{1 \ox \hat{f} \wdot X \ox 1} &&
    \hat{\Gamma}_1 \ox A \wdot X \ox \hat{\Gamma}_2 \ar[r]^-{\hat{s}} &
    \hat{\Delta}}
\]
which, by reversing representability, gives the desired polymap~\eqref{fdots}.

The composite of two polymaps can be defined in a similar manner. Suppose
that 
\[
s : \Phi \mid \Gamma \ra \Delta_1,X,\Delta_2 \qquad \text{and} \qquad
t : \Psi \mid \Gamma_1,X,\Gamma_2 \ra \Delta
\]
are polymaps represented by
\[
\hat{s} : \hat{\Gamma} \ra (\hat{\Delta}_1 \parr X) \parr \hat{\Delta}_2
\qquad \text{and} \qquad
\hat{t} : \hat{\Gamma}_1 \ox (X \ox \hat{\Gamma}_2) \ra \hat{\Delta}~.
\]
The composite of $s$ and $t$ on $X$ should be a polymap in $\P$ as
\begin{equation}\label{st}\tag{Y}
s ; t : \Phi,\Psi \mid \Gamma_1,\Gamma,\Gamma_2 \ra
\Delta_1,\Delta,\Delta_2.
\end{equation}
In $\X$ we may form the composite 
\begin{align*}
& \xygraph{
{\hat{\Gamma}_1 \ox \hat{\Gamma} \ox \hat{\Gamma}_2} 
    :[r(4)] {\hat{\Gamma}_1 \ox ((\hat{\Delta}_1 \parr X) \parr \hat{\Delta}_2) \ox \hat{\Gamma}_2} ^-{1 \ox \hat{s} \ox 1}
}
\\
& \xygraph{
    :[r(3.2)] {\hat{\Gamma}_1 \ox ((\hat{\Delta}_1 \parr X)  \ox \hat{\Gamma}_2) \parr \hat{\Delta}_2} ^-{1 \ox d^{\parr'}_\ox}
    :[r(5.1)] {\hat{\Gamma}_1 \ox (\hat{\Delta}_1 \parr (X  \ox
\hat{\Gamma}_2)) \parr \hat{\Delta}_2} ^-{1 \ox d^\parr_\ox \parr 1}
}
\\
& \xygraph{
    :[r(3.2)] {\hat{\Delta}_1 \parr (\hat{\Gamma}_1 \ox (X  \ox
\hat{\Gamma}_2)) \parr \hat{\Delta}_2} ^-{d^{\ox'}_\parr \parr 1}
    :[r(4)] {\hat{\Delta}_1 \parr \hat{\Delta} \parr \hat{\Delta}_2} ^-{1 \parr \hat{t} \parr 1}
}
\end{align*}
which, by reversing representability, gives the desired polymap~\eqref{st}.

This shows how the two cuts arise in the representable poly-actegory built
from a linear actegory. Associativity and the interchange laws follow 
from a combination of functoriality, naturality, and associativity in $\X$.
As an example suppose there are multimaps and a polymap
\[
f: \Phi \ra A, ~~ f':\Phi' \ra A', ~~\text{and}~~
s: \Psi_1,A,\Psi_2,A',\Psi_3 \mid \Gamma,X,\Gamma_2,X',\Gamma_3 \ra \Delta
\]
which are represented as
\[
\hat{f}: \hat{\Phi} \ra A, \quad \hat{f}':\hat{\Phi}' \ra A', \quad
\text{and} \quad
\hat{s}: \hat{\Gamma} \ox A \wdot X \ox \hat{\Gamma}_2 \ox A' \wdot X' \ox
         \hat{\Gamma}_3 \ra \hat{\Delta}.
\]
If the interchange law is to hold then
\[
f \cdot (f' \cdot s) = f' \cdot (f \cdot s).
\]
From the definition the left-hand side and right-hand side respectively are
given by reversing representability in the composites below.
\[
\xymatrix@R=7ex{
\hat{\Gamma} \ox \hat{\Phi} \wdot X \ox \hat{\Gamma}_2 \ox \hat{\Phi}' \wdot X' \ox \hat{\Gamma}_3
    \ar[d]_-{1 \ox \hat{f} \wdot X \ox 1 \ox 1 \ox 1} \\
\hat{\Gamma} \ox A \wdot X \ox \hat{\Gamma}_2 \ox \hat{\Phi}' \wdot X' \ox \hat{\Gamma}_3
    \ar[d]_-{1 \ox 1 \ox 1 \ox \hat{f}' \wdot X' \ox 1} \\
\hat{\Gamma} \ox A \wdot X \ox \hat{\Gamma}_2 \ox A' \wdot X' \ox \hat{\Gamma}_3
    \ar[d]_-{\hat{s}} \\
\hat{\Delta}}
\qquad\qquad
\xymatrix@R=7ex{
\hat{\Gamma} \ox \hat{\Phi} \wdot X \ox \hat{\Gamma}_2 \ox \hat{\Phi}' \wdot X' \ox \hat{\Gamma}_3
    \ar[d]^-{1 \ox 1 \ox 1 \ox \hat{f}' \wdot X' \ox 1} \\
\hat{\Gamma} \ox \hat{\Phi} \wdot X \ox \hat{\Gamma}_2 \ox A' \wdot X' \ox \hat{\Gamma}_3
    \ar[d]^-{1 \ox \hat{f} \wdot X \ox 1 \ox 1 \ox 1} \\
\hat{\Gamma} \ox A \wdot X \ox \hat{\Gamma}_2 \ox A' \wdot X' \ox \hat{\Gamma}_3
    \ar[d]^-{\hat{s}} \\
\hat{\Delta}}
\]
That both composites are equal follows from functoriality of $\ox$.

\begin{proposition}
Every linear actegory is the category of maps of some representable
poly-actegory.
\end{proposition}

This now proves Theorem~\ref{thm-main}, that to be a linear actegory is
precisely to be the category of maps of a representable poly-actegory.

Coproducts in the polycategorical setting are discussed extensively
in~\cite{P}. They have not been mentioned the discussion of this section 
as their presence or absence is completely orthogonal to the main result.  

\section{Conclusion}

We have now completed the tour of the diagram connecting proof theory, 
categorical semantics, and term calculus as promised:

\[
\xymatrix@C=5ex@R=5ex{
{\begin{array}{c} \text{proof theory} \\ \hline \text{poly-actegory}
 \end{array}}~~ \ar@{->}[dr]_{\S \ref{sec-mess-term}} \ar@{<-}[rr]^{\S \ref{sec-soundness}} &&
~~{\begin{array}{c} \text{category theory} \\ \hline \text{linear actegory}
\end{array}} \ar@{<-}[dl]^{\S \ref{semantics}}
\\
& {\begin{array}{c} \text{term calculus} \\ \hline \text{message passing}
\end{array}}}
\]

A reasonable question to ask is whether we really have a complete set of 
coherence diagrams and of equations. Our response to this is, of course, that
we did try to be reasonably complete. However, we are happy to admit that, in
a system of this size, it is quite possible that we have overlooked something.
That said, however, as we now have at least three different ways to view the 
subject matter all of which agree, we are confident that the basic story is 
complete and that these ideas, insofar as they are not completely fleshed 
out here, can be.

The aim of the paper was to show that message passing could be accommodated
in the proof theoretic framework for concurrency provided by the two-sided
proof theory of cut elimination (which is a fragment of linear logic). In
particular we feel that by providing a categorical semantics we have
anchored this correspondence in a way which will facilitate the exploration
of semantic models. Thus, we feel that we have now lain out the story of how
message passing can be modeled proof theoretically and categorically in
sufficient detail to establish the viability of this perspective. There
remains a lot to be done.

In the process of writing the paper the proof system underwent a number of
downsizing changes in an attempt to make it more manageable. For example,
the additives at the message passing level were sacrificed for this reason.
From the programming perspective, features such as (the initial and final)
datatypes are desirable \emph{at both levels} and this is an aspect to which
we would like to return. Particularly, at the message passing level
datatypes are of significant interest as they allow the expression of
communication protocols in a formal manner.



\end{document}